\documentclass[english,11pt]{article}%
\usepackage{bbm,ulem}
\usepackage{graphicx}
\usepackage{subcaption}
\usepackage{makeidx}
\usepackage{amssymb}

\usepackage{tikz-cd}

\usepackage[utf8]{inputenc}
\usepackage{amsfonts}
\usepackage{amsmath}
\usepackage{amssymb}
\usepackage{amsthm}
\usepackage{url}
\usepackage[colorlinks,citecolor=blue]{hyperref}
\usepackage{mathabx}
\usepackage{wrapfig}
\usepackage{dsfont}
\usepackage{resizegather}
\usepackage{yfonts}
\usepackage{authblk}
\usepackage{geometry}

\DeclareMathOperator*{\argmin}{arg\,min}

\providecommand{\U}[1]{\protect\rule{.1in}{.1in}}

\usepackage{xcolor}

\newtheorem{prop}{Proposition}[section]
\newtheorem{cor}[prop]{Corollary}

\newtheorem{rmk}[prop]{Remark}
\newtheorem{lem}[prop]{Lemma}

\newtheorem{theo}[prop]{Theorem}
\newtheorem{examp}[prop]{Example}

\newcommand{\tr}{\mbox{\rm Tr}}

\newcommand{\LL}{\mathbb{L}}

\newcommand{\PP}{\mathbb{P}}

\newcommand{\RR}{\mathbb{R}}

\newcommand{\TT}{\mathbb{T}}
\newcommand{\UU}{\mathbb{U}}
\newcommand{\VV}{\mathbb{V}}
\newcommand{\XX}{\mathbb{X}}

\newcommand{\YY}{\mathbb{Y}}

\newcommand{\Ba}{ {\cal B }}

\newcommand{\Da}{ {\cal D }}

\newcommand{\Ka}{ {\cal K }}

\newcommand{\Sa}{ {\cal S }}

\newcommand{\Va}{ {\cal V }}
\newcommand{\Ua}{ {\cal U }}

\newcommand{\Ga}{ {\cal G }}

\newcommand{\Ha}{ {\cal H }}
\newcommand{\Ja}{ {\cal J }}
\newcommand{\Pa}{ {\cal P }}
\newcommand{\Za}{ {\cal Z }}

\newcommand{\Wa}{ {\cal W }}
\newcommand{\Qb}{ {\bf Q}}

\newcommand{\Pb}{ {\bf P }}

\newcommand{\point}{\mbox{\LARGE .}}
\newcommand{\cqfd}{\hfill\blbx \\}
\def\blbx{\hbox{\vrule height 5pt width 5pt depth 0pt}\medskip}

\def \PP{\mathbb{P}}
\def \RR{\mathbb{R}}

\def \LL{\mathbb{L}}
\def \JJ{\mathbb{J}}

\newcommand{\cchi}{\protect\raisebox{2pt}{$\chi$}}

\numberwithin{equation}{section}

\newcommand{\vertiii}[1]{{\left\vert\kern-0.25ex\left\vert\kern-0.25ex\left\vert #1
    \right\vert\kern-0.25ex\right\vert\kern-0.25ex\right\vert}}

\usepackage{amsopn}
\begin{document}

  \title{New Trends in the Stability of Sinkhorn Semigroups }

\author{P. Del Moral\thanks{Centre de Recherche Inria Bordeaux Sud-Ouest, Talence, 33405, France. {\footnotesize E-Mail:\,} \texttt{\footnotesize pierre.del-moral@inria.fr}} \& A. Jasra\thanks{School of Data Science, The Chinese University of Hong Kong, Shenzhen,  Shenzhen, China. {\footnotesize E-Mail:\,} \texttt{\footnotesize ajayjasra@cuhk.edu.cn}}}

\maketitle
  \begin{abstract}    
  \noindent Entropic optimal transport problems play an increasingly important role in machine learning and generative modelling.  In contrast with optimal transport maps which often have limited applicability in high dimensions,  Schr\"odinger bridges 
 can be solved using the celebrated Sinkhorn's algorithm, a.k.a.~the iterative proportional fitting procedure. The stability properties of Sinkhorn bridges when the number of iterations
 tends to infinity is a very active research area in applied probability and machine learning.
 Traditional proofs of convergence are mainly based on nonlinear versions of Perron-Frobenius theory and related Hilbert projective metric techniques, gradient descent, Bregman divergence techniques and Hamilton-Jacobi-Bellman equations, including propagation of convexity profiles based on coupling diffusions by reflection methods.
  The objective of this review article is to present, in a self-contained manner,  recently developed
 Sinkhorn/Gibbs-type semigroup analysis based upon contraction coefficients and
 Lyapunov-type operator-theoretic techniques. These  powerful,  off-the-shelf semigroup methods
are based upon transportation cost inequalities (e.g.~log-Sobolev, Talagrand quadratic inequality, curvature estimates),  $\phi$-divergences,  Kan\-to\-rovich-type criteria and
Dobrushin contraction-type coefficients on  weighted Banach spaces as well as Wasserstein distances. 
This novel semigroup analysis allows one to unify and simplify many arguments in the stability of Sinkhorn algorithm. It also yields new contraction estimates w.r.t.~generalized $\phi$-entropies, as well as weighted total variation norms, Kantorovich criteria and Wasserstein distances.
\end{abstract}

\newpage
{\scriptsize
\tableofcontents
}

\newpage

\section{Entropic optimal transport}
\subsection{Schr\" odinger and Sinkhorn bridges}
Let $\mu$ and $\eta$ be two probability measures 
on measurable spaces $\XX$ and $\YY$ and  
let $\Pi(\mu,\eta)$ be the set of couplings of $\mu$ and $\eta$ on  the product space $(\XX\times\YY)$, with 
prescribed first and second coordinate marginals $(\mu,\eta)$.  The (static) Schr\"odinger bridge $P_{\mu,\eta}$ from $\mu$ to $\eta$ with respect to a reference probability measure $P$ on  $(\XX\times\YY)$ is defined by 
  \begin{equation}\label{def-entropy-pb-v2}
 P_{\mu,\eta}:=\argmin_{Q\,\in\, \Pi(\mu,\eta)}\Ha(Q|P)
\end{equation}
where $\Ha(Q|P)$ is the relative entropy of $Q$ w.r.t.~$P$;  $\Ha(Q|P)=\int \log\left(\tfrac{dQ}{dP}\right)dQ$,
$\tfrac{dQ}{dP}$ the Radon-Nikodym derivative of $Q$ with respect to $P$.
The bridge map $P\mapsto P_{\mu,\eta}$  is well defined  
as soon as  
$\Ha(Q|P)<\infty$ for  some $Q\in \Pi(\mu,\eta)$,  see for instance,  the seminal article by Csisz\'ar~\cite{csiszar-2}, \cite[Section 6]{nutz} and the survey article by L\'eonard~\cite{leonard} including the references therein. 
In contrast with conventional optimal transport maps, with limited applicability in high dimensions,  Schr\"odinger bridges (\ref{def-entropy-pb-v2})
 can be solved using the celebrated Sinkhorn's algorithm, a.k.a.~iterative proportional fitting procedure or bi-proportional fitting in the statistics literature~\cite{cuturi,genevay-cuturi-2,genevay,sinkhorn,sinkhorn-2,sinkhorn-3}, as well as the relative state attraction and matrix scaling~\cite{faharat,idel,schoen,schoen2} approaches  in computer sciences.  These ideas provide a concrete roadmap to solve optimal transport problems and,  as mentioned,  detail numerous applications in computer science and statistics.
 
We now briefly outline the basic principles of the algorithm; Section~\ref{state-space-conditioning-sec}  presents a more detailed summary.  Given a  Markov transition $K(x,dy)$ from $\XX$ into $\YY$ we set
\begin{eqnarray}
 (\mu K)(dy)&:=&\int_{\XX}\mu(dx)\, K(x,dy)\quad\mbox{\rm and}\quad
 (\mu\times K)(d(x,y)):=\mu(dx)\, K(x,dy).
 \label{muK}
\end{eqnarray}
Given a probability measure $R(d(y,x))$ on $(\YY\times\XX)$, 
 we denote by $ R^{\flat}$ the probability measure on $(\XX\times\YY)$ defined by
 $$
 R^{\flat}(d(x,y))= R(d(y,x)).
 $$
Observe that
\begin{equation}\label{def-entropy-note}
 P:=\mu\times K
\Longrightarrow
P=P_{\mu,\mu K}\in\Pi(\mu,\mu K)
\quad \mbox{\rm and}\quad
(P_{\mu,\eta})^{\flat}=P^{\flat}_{\eta,\mu}:=(P^{\flat})_{\eta,\mu}.
\end{equation}
The Sinkhorn iterations are defined sequentially for any $n\geq 0$ by a collection of probability distributions
\begin{equation}\label{def-Pa-n}
 \Pa_{2n}=\mu\times\Ka_{2n}
\quad \mbox{\rm and}\quad
 \Pa_{2n+1}=(\eta\times\Ka_{2n+1})^{\flat}
\end{equation}
starting from a reference measure of the form $\Pa_0=P:=(\mu\times K)$ at rank $n=0$.
For $n \ge 1$, the  Markov transitions $\Ka_n$ in (\ref{def-Pa-n}) {are defined} sequentially by the conditioning Bayes-type  formulae
\begin{equation}
\left\{\begin{array}{l}
(\pi_{2n}\times\Ka_{2n+1})^{\flat}=\mu\times\Ka_{2n}\quad \mbox{\rm and}\quad
\pi_{2n+1}\times\Ka_{2(n+1)}=(\eta\times\Ka_{2n+1})^{\flat}\\
\\
\mbox{\rm with the distributions}\quad \pi_{2n}:=\mu \Ka_{2n}\quad\mbox{\rm and}\quad \pi_{2n+1}:=\eta \Ka_{2n+1}.
\end{array}\right.
\label{s-2}
\end{equation}
When $n\rightarrow\infty$, Sinkhorn bridges $\Pa_{n}$ converge towards the Schr\"odinger bridge $P_{\mu,\eta}$ and the marginal measures $(\pi_{2n+1},\pi_{2n})$
converge towards $(\mu,\eta)$. In a mode of convergence that we define later,  one has
$$
\Pa_{n}~\longrightarrow_{n\rightarrow\infty}~P_{\mu,\eta}\quad \mbox{\rm and}\quad
(\pi_{2n+1},\pi_{2n})~\longrightarrow_{n\rightarrow\infty}~(\mu,\eta).
$$
The limiting properties as $n$ grows are what we call, loosely,  stability and this term also includes rates of convergence.

\subsection{Sinkhorn semigroups}\label{sinkhorn-sg-intro-sec}

For finite state spaces, the Sinkhorn algorithm is a well known simple procedure that repeatedly applies  column and row normalization on matrix spaces~\cite{deming,fortet,kruithof,sinkhorn-2,sinkhorn-3}. 
In problems of low-to-moderate dimension,  whilst the recursions can be computed with ease,  it is challenging to obtain explicit rates of convergence. 
For Gaussian models, Sinkhorn equations can also be solved analytically using sequential Bayes conjugacy principles~\cite{adm-24,cramer,cramer-2,cramer-3,janadi}.  As expected, this recursive formulation is closely related to the Kalman filter and related Riccati matrix difference equations, and it yields algorithms that can be implemented in practical settings without further approximations. A self contained analysis of Schr\" odinger bridges and the stability properties of Sinkhorn bridges for a general class of linear-Gaussian models is provided in Section~\ref{Gauss-review-sec}. It is noted that in this Gaussian setting,  the stability of the recursions is related to the stability of the Riccati equation;  this concept is made concrete later in the article.
In more general situations, while Sinkhorn iterations as presented in (\ref{s-2}) may look appealing {and easy to implement},  one should note that they are based on nonlinear conditional/conjugate transformations with generally no finite-dimensional recursive solutions and, therefore, they do not lead to a practical algorithm. In this context, 
 stability analysis of Sinkhorn bridges is crucial in the convergence analysis of any numerical approximations, including time discretization of diffusion bridges, neural networks  and related
 score matching approximation techniques~\cite{schen1,tchen,bortoli,doucet-bortoli,sohl,song1,song2}.

As underlined in~\cite{adm-24}, the key observation is that
Sinkhorn semigroups belong to the class of time varying Markov processes. 
Indeed, integrating (\ref{s-2}) we check the collection of fixed point equations
\begin{equation}\label{pre-gibbs}
\mu=\pi_{2n}\Ka_{2n+1}=\mu \Sa_{2n+1}
\quad\mbox{\rm and}\quad
\eta= \pi_{2n+1}\Ka_{2(n+1)}=
\eta \Sa_{2(n+1)}
\end{equation}
with the forward-backward Markov transitions
\begin{equation}\label{fb-sinkhorn}
\begin{array}{l}
\Sa_{2n+1}:=
\Ka_{2n}\Ka_{2n+1}\quad\mbox{\rm and}\quad
\Sa_{2(n+1)}:=
\Ka_{2n+1}\Ka_{2(n+1)}\\
\\
\Longrightarrow
 \pi_{2n-1}\Sa_{2n+1}=\pi_{2n+1}
\quad\mbox{\rm and}\quad
 \pi_{2n}\Sa_{2(n+1)}=\pi_{2(n+1)}
\end{array}
\end{equation}
where the product of transitions $\Ka_{n}\Ka_{n+1}$ is a shorthand notation for the composition $\Ka_{n}\circ \Ka_{n+1}$ of the right or left-action operators defined respectively in (\ref{ract-ref}) (operation on functions) and in (\ref{muK}). In view of the fixed point equations (\ref{pre-gibbs}) along with \eqref{fb-sinkhorn}, when $n\rightarrow\infty$, we expect Sinkhorn distributions $\pi_{2n}$ to converge towards $\eta$ while $\pi_{2n+1}$ converge towards $\mu$.

A schematic picture of the time varying forward-backward semigroup of the distribution flow $\pi_{2n}$ is given below:
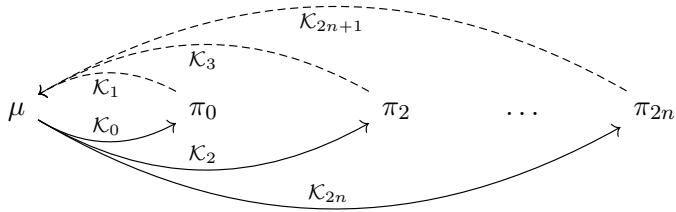
\[
\begin{tikzcd}
\mu   \ar[rrrrrr,bend right,"\Ka_{2n}"]   \ar[rrrr,bend right,"\Ka_2"] \ar[rr,bend right,"\Ka_0"] &&  \ar[ll, "\Ka_1", bend right,  dashrightarrow]  \pi_0&& \ar[llll,bend right,"\Ka_3", dashrightarrow] \pi_2&\ldots& \ar[llllll,bend right,"\Ka_{2n+1}", dashrightarrow] \pi_{2n}
\end{tikzcd}
\]
Note that the product $\Sa_{2n+1}=\Ka_{2n}\Ka_{2n+1}$ can be interpreted as the $2$-block Gibbs transition with target measure $\Pa_{2n}$ on $(\XX\times\YY)$. The transition $\Ka_{2n}(x,dy)$ is the $\Pa_{2n}$-distribution of the coordinate $y$ given $x$; while
$\Ka_{2n+1}(y,dx)$ is the $\Pa_{2n}$-distribution of the coordinate $x$ given $y$.

In addition, whenever the sequence of entropies of Sinkhorn distributions (\ref{fb-sinkhorn}) w.r.t.~the fixed point measures (\ref{pre-gibbs}) are  summable we have
the equivalence
\begin{equation}\label{kdec2n-intro}
\lim_{n\rightarrow\infty}\Ha(P_{\mu,\eta}|\Pa_{2n})=0\Longleftrightarrow
\Ha(P_{\mu,\eta}|\Pa_{2n})=\sum_{l\geq n}\left[\Ha(\eta|\pi_{2l})+\Ha(\mu|\pi_{2l+1})\right].
\end{equation}
A detailed proof of the above assertion is provided in  Section~\ref{enproj-sec}, see also Proposition~\ref{prop-decrease2n01} for equivalent conditions. This crucial observation shows that the convergence of Sinkhorn bridges
is exactly characterized by the convergence of Sinkhorn semigroups toward their limiting distributions.

As shown in (\ref{pre-gibbs}) and (\ref{fb-sinkhorn}) Sinkhorn transition semigroups $\Sa_n$ belong to the class of time varying Markov processes  sharing a common invariant target measure at each individual step. 
These classes of time varying Markov semigroups arise in a variety of areas in applied mathematics including nonlinear filtering, physics and molecular chemistry; see for instance~\cite{ado-24,dm-04,dg-01,dhj,saloff-zuniga,saloff-zuniga-2,saloff-zuniga-3}  and references therein. 
The stability analysis of time inhomogenous Markov processes is a broad and notoriously difficult subject, 
far more challenging than their time homogeneous counterparts; this is because the operators may drastically change during the semigroup evolution.
The articles~\cite{douc,saloff-zuniga,saloff-zuniga-2} provide several useful probabilistic and analytical  tools, including coupling methods, spectral analysis as well as functional inequalities such as Nash or log-Sobolev inequalities.  Nevertheless, Sinkhorn semigroups are defined sequentially by nonlinear conjugate transformations and hence the application of these refined analytical techniques in this context is rather challenging. 

The nonlinear Sinkhorn/Gibb-type semigroup approach taken in this article differs significantly from the conventional approaches developed in the optimal transport literature; this is discussed in some detail in Section~\ref{lit-review}. This novel framework allows one to unify and simplify many arguments in the stability of Sinkhorn bridges, as well as to improve upon some existing results. 
We shall survey, in a self contained manner,  some powerful off-the-shelf semigroup techniques recently developed in the series of articles~\cite{adm-24,adm-25,dp-25}, themselves based on the contraction operator-theoretic analysis developed in~\cite{ado-24,dhj,dg-01,dm-03,dmg-25,dm-penev-2017}.  The approaches that we focus upon are:
\begin{itemize}
\item{Transportation cost inequalities including local log-Sobolev and Talagrand quadratic cost inequalities as well as Fisher-Lipschitz-type inqualities and a novel entropic stability theorem.}
\item{An operator-theoretic framework based on contraction coefficients, including Lyapunov methodologies on Weighted Banach spaces,  Kantorovich criteria and Wasserstein distances. }
\end{itemize}

In the context of Sinkhorn and Schr\"odinger  bridges, the article~\cite{dp-25} presents a  semigroup approach
combining transportation cost inequalities including log-Sobolev and Talagrand quadratic inequalities 
  with the theory of Riccati matrix difference equations arising in filtering and optimal control theory.  In the present article, we illustrate these inequalities in the context of log-concave-at-infinity models.  Other types of models, such as  sub-Gaussian or log-concave,  are discussed in~\cite{dp-25}.  The analysis of marginal models with lower/upper bounded curvature is closely related to Gaussian models.
  Using Brascamp-Lieb and Cramer-Rao inequalities,  the analysis of these strongly convex models also relies on discrete time algebraic Riccati equations (see for instance Appendix D in~\cite{dp-25} and Section~\ref{mr-stconv-sec} in the present article). When applied to the stability of Sinkhorn semigroups, they also yield
 a series of novel contraction estimates in terms of fixed point of Riccati equations.   
As mentioned above, a detailed analysis of Sinkhorn bridges for linear-Gaussian models is presented in \cite{adm-24}, including finite-dimensional recursive formulation of the Sinkhorn algorithm and sharp convergence rates in terms of Riccati matrix difference equations. 

The  article~\cite{adm-25} introduces an operator-theoretic framework yielding  new exponential decays of Sinkhorn iterates towards Schr\"odinger bridges with respect to general classes of $\phi$-divergences and Kantorovich-type criteria, including the relative entropy,  squared Hellinger integrals, $\alpha$-divergences,
$\LL_p$-norms,  Havrda-Charvat entropies,
Jensen-Shannon divergence, as well as  Jeffreys and R\'enyi divergences. \cite{adm-25} also develops a novel 
  Lyapunov contraction approach under minimal regularity conditions. This methodology  
   provides quantitative exponential stability estimates for a large class of Sinkhorn semigroups
  on weighted Banach spaces.  Contraction estimates are expressed in terms of  Kantorovich semi-distances, Wasserstein distances as well as Dobrushin-type contraction coefficients and weighted total variation norms.
   This novel Lyapunov methodology is applied to  a variety of situations, ranging from polynomial growth potentials and heavy tailed marginals on general normed spaces,  to more sophisticated boundary state space models, including semi-circle transitions, Beta, Weibull, exponential marginals as well as  semi-compact models. This approach also allows one to consider statistical finite mixture of the above models, including kernel-type density estimators of complex data distributions arising in generative modeling.

As a result of the methods we consider, in this article,  we can often obtain quantitative and refined estimates. A
precise description of the convergence rates and estimation constants are provided in the series of articles~\cite{adm-24,adm-25,ado-24,dp-25}.  We remark that the articles~\cite{adm-24,adm-25,dp-25} are partly based on well-established techniques in Riccati difference equations and in the stability theory of Markov semigroups. 
 The exponential stability of Sinkhorn bridges in~\cite{adm-24,dp-25} rely on
 the stability analysis of Riccati matrix difference equations presented in~\cite{dh-22}.
 The contraction properties of Sinkhorn semigroups presented in~\cite{adm-25} combine  the contraction
  analysis of Markov semigroups  developed in~\cite{dg-01,dm-03} with  the Lyapunov operator theoretic framework on weighted Banach spaces  presented in
 ~\cite{dhj,dm-penev-2017},
and further developed in the articles
 ~\cite{ado-24,dmg-25}.  In the context of time-homogeneous models, we also refer to the pioneering articles by Hasminskii~\cite{hasminskii} and Meyn and Tweedie~\cite{meyn-1}, see also the book~\cite{meyn-2} and the more recent seminal article by Hairer and J. C. Mattingly~\cite{hairer-mattingly}.

\subsection{Literature review}\label{lit-review}

The stability properties of Schr\" odinger and Sinkhorn bridges when the number of iterations
 tends to $\infty$ is a very active research area in applied probability and machine learning.
 Most of the literature on the exponential convergence of Sinkhorn equations is concerned with
finite state spaces~\cite{borwein,fienberg-1970,sinkhorn-2,sinkhorn-3,soules} as well as compact state spaces or bounded cost functions using nonlinear versions of Perron-Frobenius theory~\cite{brualdi,idel} and related
Hilbert projective metrics techniques~\cite{chen,deligiannidis,franklin,marino}, based on the seminal article by 
Birkhoff~\cite{birkoff}. See also~\cite{greco2023non}
for a refined analysis of Sinkhorn iterates and their gradient on the $d$-dimensional torus as well as the more recent article~\cite{chizat2025sharper} on compact spaces and bounded cost functions.  
 
 Sublinear rates have been developed in the articles~\cite{alts-2017,chak-2018,dvu-2017}.
 Linear rates  on non-necessarily compact spaces were first obtained by L\'eger in~\cite{leger} using elegant gradient descent and Bregman divergence techniques (see also the recent articles~\cite{doucet-bortoli,karimi-2024}). Refined convergence rates at least one order faster for a general class of models including sub-Gaussian target marginals has also been developed in~\cite{promit-2022} under some exponential moments estimates on Schr\" odinger potentials.

 The first exponential stability analysis of unbounded models and not necessarily compactly supported marginals is presented in the article~\cite{durmus}. In this study, the authors investigate quantitative contraction rates for target marginal distributions with an asymptotically positive log-concavity profile and quadratic symmetric cost functions  associated with a sufficiently large regularization parameter. The approach taken in~\cite{durmus} relies on well known links between gradient and Hessian of Sinkhorn and Schr\" odinger potentials along Hamilton-Jacobi-Bellman equations with conditional expectations and covariances~\cite{chewi2023entropic,pooladian2021entropic}, see also~\cite{del2025entropic}. The main conditions in~\cite{durmus} are expressed in terms of the integrated convexity profiles of marginal potentials. Extending ideas from~\cite{conforti2024weak}  the proofs are based on the propagation of these convexity profiles based on coupling diffusion  by reflection techniques introduced in~\cite{eberle2016reflection}.  The exponential entropic decays presented in~\cite{durmus} were further developed and recently refined in~\cite{chiarini} to  apply to all values of the regularization parameter and some classes of perturbation of log-concave models. However this stability approach requires one to estimate the semi-concavity properties of generally unknown dual Schr\"odinger bridge transitions and it only applies to entropically regularized transport problems with symmetric quadratic-type costs. Thus, this framework excludes the important cases when the reference transition is associated with linear Gaussian transitions arising in Ornstein-Uhlenbeck-type diffusion generative models
and denoising diffusions. 

The article~\cite{eckstein} is also based on an extension of Hilbert’s projective metric for spaces of integrable functions of bounded growth. Extending results in~\cite{chen,deligiannidis} to unbounded cost functions, the recent article ~\cite{eckstein} presents exponential decays of these new classes of
Hilbert’s projective metrics with respect to some cone of functions
 in settings where the marginal distributions have sufficiently light tails compared to the growth of the cost function.   These results also imply the exponential convergence of the total variation distance. On normed finite dimensional spaces, the main condition in ~\cite{eckstein} is that  the log of the reference transition density
grows faster that a polynomial of order $p$ and the tails of both marginals decay faster that $e^{-r^{p+q}}$ for some $q>0$. This condition does not cover Gaussian models.
  It is out of the scope of this article to present more details of
 these different approaches.  
  Despite the numerous references given above, a complete literature review is not possible.  
We  refer to reader to the books~\cite{nutz,peyre} as well as the survey articles~\cite{chen-2,garg,huynh,idel,sejourne,wang,zhang},  recent article articles~\cite{adm-24,adm-25,chiarini,durmus,dp-25,gt} and the numerous references therein.
 We underline that all the methodologies discussed above differ from the recently developed Sinkhorn semigroup contraction analysis 
 and the Lyapunov operator-theoretic framework discussed in Section~\ref{sinkhorn-sg-intro-sec}.  
 
\subsection{Article focus and structure}

The objective of this survey article is to detail,  in a pedagogical manner,  the subject of Sinkhorn semigroup techniques, along with quantitative stability estimates.  Although,  as mentioned above,  we do discuss some existing results,  we do not survey the entire literature with all of the results that have been derived.  One of the main motivations of this paper is to equip the reader with many of the tools that are needed to analyze Sinkhorn semigroups,  without having to go beyond the technical methods that have been employed in the analysis of non-homogenous Markov chains (e.g.~\cite{dhj,douc}).  This in turn leads to short and elegant proofs of existing results,  along with,  in some cases,  a substantial improvment over those currently published.   We again remark that,  other than the finite state-space or Gaussian case,  Sinkhorn algorithm recursions are numerically intractable.
One an implication is that an important ingredient in the mathematical analysis of numerical approximations is the time stability of the previously mentioned recursions.  Indeed an understanding of time stability could lead to a motivation for \emph{new numerical methods}.
This idea is analogous to the stability of the non-linear filter which has been studied in some depth for more than 30 years (e.g.~\cite{dg-01} or \cite{ba-cr-08,dm-04} and the references therein).  Once one understands whether a non-linear filter is stable,  one may hope to derive numerical algorithms which approximate the filter and can be stable themselves (we refer the reader to \cite{ba-cr-08,dm-04} to understand what stability is in non-linear filtering).

This article is structured as follows.   In the following subsection \ref{sec:notation} we provide the main notations that are adopted in the paper.  Section \ref{sec:prelim} presents some background analysis on Schr\" odinger and Sinkhorn bridges.   In section \ref{sec:stab_markov} we give a summary on the stability of Markov operators which is the approach that we adopt.  In section \ref{SS-bridges} a brief review of Sinkhorn bridges is given. 
Finally in section \ref{sec:stability_analysis} we detail our stability analysis. 

\subsection{Some basic notation}\label{sec:notation}
\subsubsection*{Integral operators}
Given a  Markov transition $K(x,dy)$ from $\XX$ into $\YY$ and a function $f$ on $\YY$ we set
\begin{equation}\label{ract-ref}
 K(f)(x):=\int_{\YY}~K(x,dy)~f(y).
\end{equation}
We are using $K$ to denote both the left-action operator on measure (\ref{muK}) and the right-action operator on functions (\ref{ract-ref}).
Denote by
$ \Pi(\mu,\point)$ be the set of probability measures on $(\XX\times\YY)$ with $\XX$-marginal $\mu$; and 
 $ \Pi(\point,\eta)$  is the set of probability measures on $(\XX\times\YY)$ with $\YY$-marginal $\eta$.  Note that
 $$
( \mu\otimes \eta)(d(x,y)):= \mu(dx)~\eta(dy)\Longrightarrow
(\mu\otimes \eta)\in \Pi(\mu,\eta).
 $$
The Dirac delta measure located at $x$ is denote by $\delta_x$.
 We also denote by $\vert\nu \vert:=\nu_++\nu_-$  the total variation measure associated with
 a Hahn-Jordan decomposition $\nu=\nu_+-\nu_-$ of a bounded signed measure $\nu$.
For indicator functions $f=1_{A}$ of a measurable subset $A\subset \XX$ and a given measure $\mu$ on $\XX$ sometimes we write $\mu(A)$ instead of $\mu(1_A)$. 

 Given a probability measure $P(d(x,y))$ on a product space $(\XX\times\YY)$ we denote by
 $ P^{\flat}$ the probability measure on $(\YY\times\XX)$ defined by
 $$
 P^{\flat}(d(y,x))= P(d(x,y)).
 $$
When $R(d(y,x))$ is a probability measure $(\YY\times\XX)$, 
 $ R^{\flat}$ is the probability measure on $(\XX\times\YY)$ defined by
 $$
 R^{\flat}(d(x,y))= R(d(y,x)).
 $$
To illustrate the notation, for any Markov transition $L(y,dx)$ from $\YY$ into $\XX$ we have
$$
R=\eta \times L\Longleftrightarrow
R(d(y,x))=(\eta\times L)(d(y,x))=\eta(dy)\, L(y,dx).
$$
In this situation, we have
\begin{eqnarray*}
 R^{\flat}(d(x,y))=(\eta \times L)^{\flat}(d(x,y)):=\eta(dy)\, L(y,dx)&\Longrightarrow&
(\eta \times L)^{\flat}\in\Pi(\eta L,\eta).
\end{eqnarray*}

 We let $\Ba(\XX)$ denote the set of measurable functions on $\XX$  and  $\Vert f\Vert:=\sup_{x\in \XX}|f(x)|$  denotes the uniform norm of a function $f\in \Ba(\XX)$.
Denote by ${\sf Prob}(\XX)$ be the set of  probability measures on $\XX$. 
Next, we let $g\in\Ba(\XX)$ be a  lower semi-continuous (abbreviated l.s.c.) function, lower bounded away from zero, and  we denote by $\Ba_g(\XX)$   the subspace of functions $f\in \Ba(\XX)$ such that $\Vert f/g\Vert<\infty$.  
Denote by ${\sf Prob}_g(\XX)$ be the set of  probability measures $\mu$ on $\XX$ such that $\mu(g)<\infty$.   
$\Ba_{\infty}(\XX)$ is the collection of locally upper bounded and uniformly positive functions  
 that growth at infinity (in the sense that  their sub-level sets are compact).  $\Ba_0(\XX)= \{\theta~:~1/\theta \in B_{\infty}(\XX)\}$ is the sub-algebra of (bounded) positive functions $\theta(x)$ which are locally lower bounded and vanish at infinity  in the sense that, for any $0<\epsilon<\Vert \theta \Vert$, the set $\{\theta\geq \epsilon\}$ is a non-empty compact set. 

\subsubsection*{Differential operators}
We let $\nabla f(x)=\left[\partial_{x_i}f(x)\right]_{1\leq i\leq d}$ be the gradient column vector associated with some smooth function $f(x)$ from $\RR^d$ into $\RR$.
 Given some smooth function $F(x)$ from $\RR^d$ into $\RR^p$
we denote by $\nabla h=\left[\nabla F^1,\ldots,\nabla F^p\right]\in \RR^{n\times p}$ the gradient matrix associated with the column vector
 function $F=(F^i)_{1\leq i\leq p}$, for some $p\geq 1$. For instance, for a given matrix 
 $A\in\RR^{q\times p}$ we have $Ah(x)\in\RR^q$ and therefore
\begin{equation}\label{rule-ah}
 \nabla (A F(x))=(\nabla F(x))~A^{\prime}
\end{equation}
where $A^{\prime}$ is the transpose
of the matrix $A$.
Given a smooth  function $f(x,y)$ from $\RR^d\times\RR^d$ into $\RR$ and $i=1,2$
  we denote by $\nabla_i$ the gradient operator defined by
   $$
 \nabla_if=\left[
 \begin{array}{c}
 (\nabla_if)^1\\
 \vdots\\
  (\nabla_if)^d
 \end{array}
 \right] \quad\mbox{\rm with}\quad
  (\nabla_1f)^j:= \partial_{x^j}f \quad\mbox{\rm and}\quad
  (\nabla_2f)^j:= \partial_{y^j}f.
 $$
For smooth $p$-column valued functions $F(x,y)=(F^j(x,y))_{1\leq j\leq p}$ we also set
$$
  \nabla_i F=\left[ \nabla_i F^1,\ldots,
  \nabla_iF^p\right].
   $$ 
In this notation, the second order differential operators  $\nabla_i \nabla_j$ on function $f(x,y)$ from $\RR^d\times\RR^d$ into $\RR$  are defined by $(\nabla_i \nabla_j) (f):=  \nabla_i (\nabla_j(f))$ and we set $
  \nabla_i^2:=  \nabla_i \nabla_i$ when $i=j$.
  
\subsubsection*{Conditional covariances}

Let $\Ga l_d$ denote the general linear group of $(d\times d)$-invertible matrices, and $\Sa^+_d\subset \Ga l_d$ the subset 
 of positive definite matrices. We sometimes use the L\" owner partial ordering notation $v_1\geq v_2$ to mean that a symmetric matrix $v_1-v_2$ is positive semi-definite (equivalently, $v_2 - v_1$ is negative semi-definite), and $v_1>v_2$ when $v_1-v_2$ is positive definite (equivalently, $v_2 - v_1$ is negative definite). Given $v\in\Sa_d^+$ we denote by $v^{1/2}$  the principal (unique) symmetric square root.

  Consider bounded positive integral operators from $\RR^d$ into itself of the form
\begin{equation}\label{def-KV}
  K(x,dy)=e^{-W(x,y)}~dy\quad \mbox{\rm and}\quad
  K_V(x,dy):=\frac{  K(x,dy)~ e^{-V(y)}}{K(e^{-V})(x)}
 \end{equation}
  for some real valued functions $W(x,y)$ and $V(y)$. Denote by 
  $W_x$ and $W_{x,x}$ the gradient and the Hessian of the transition potential
  $W(x,y)$ with respect to the first coordinate; that is 
\begin{equation}\label{ref-grad-hessian-W}
W_x(y):=\nabla_1W(x,y)\quad \mbox{\rm and}\quad
W_{x,x}(y):=\nabla_1^2W(x,y).
\end{equation}
The gradient and the Hessian formulae are
\begin{equation}\label{ref-grad-hessian}
  \begin{array}{l}
 f(x):= \log{K(e^{-V})(x)}\\
 \\
 \Longrightarrow
 \nabla f(x)=-K_V(W_x)(x)
 \quad \mbox{\rm and}\quad
  \nabla^2 f(x)=-K_V(W_{x,x})(x)+  \mbox{\rm cov}_{K_V,W}(x)
 \end{array} \end{equation}
 with the conditional covariance
\begin{equation}\label{ref-cond-cov}
  \mbox{\rm cov}_{K_V,W}(x):=\frac{1}{2}\int~
  K_{V}(x,dy_1)~  K_{V}(x,dy_2)~(W_x(y_1)-W_x(y_2))~(W_x(y_1)-W_x(y_2))^{\prime}
\end{equation}
where we recall $v^{\prime}$ is the transpose
of a matrix or vector $v$. Denote by $\mbox{\rm mean}_{K}(x)$ an,d $\mbox{\rm cov}_{K}(x)$ the conditional mean and covariance functions
\begin{equation}\label{cov-mat-W}
\begin{array}{rcl}
\mbox{\rm mean}_{K}(x)&:=&\displaystyle\int~K(x,dy)~y\\
&&\\
\mbox{\rm cov}_{K}(x)&:=&\displaystyle\frac{1}{2}\int~
  K(x,dy_1)K(x,dy_2)~(y_1-y_2)~(y_1-y_2)^{\prime}.
  \end{array}
\end{equation}
Assume the log-density function $W$ satisfy the strong convexity condition
\begin{equation}\label{strong-convex}
\tau^{-1}\leq \nabla_2^2W(x,y)\leq\overline{\tau}^{-1}\quad \mbox{\rm for some $\tau,\overline{\tau}\in \Sa_d^+$.}
\end{equation} 
In this context, the Brascamp-Lieb inequality ensures that
$$
\mbox{\rm cov}_{K}(x)\leq \int K(x,dy)~\left(\nabla^2_2W(x,y)\right)^{-1}\leq \tau.
  $$
Conversely,  the Cramer-Rao inequality also ensures that
$$
\mbox{\rm cov}_{K}(x)^{-1}\leq \int K(x,dy)~\nabla^2_2W(x,y)\leq\overline{\tau}^{-1}.
$$
We summarize the above discussion with  the covariance estimate
\begin{equation}\label{bl-cr}
(\ref{strong-convex})\Longrightarrow
\overline{\tau}\leq \mbox{\rm cov}_{K}(x)\leq \tau.
\end{equation}

\subsubsection*{Riccati difference equations}
The Frobenius matrix norm of a given matrix $v\in \RR^{d_1\times d_2}$ is defined by
$
\left\Vert v\right\Vert_{F}^2=\tr(v^{\prime}v)
$, 
with the trace operator $\tr(\cdot)$.  The spectral norm is defined by $\Vert v\Vert_2=\sqrt{\lambda_{\sf max}(v^{\prime}v)}$ and $\lambda_{\sf max}(v^{\prime}v)$ is the maximum eigenvalue of $v^{\prime}v$.
Note that, for $x\in\RR^{d \times 1}=\RR^d$, the Frobenius and the spectral  norm $\Vert x\Vert_F=\sqrt{x^{\prime}x}=\Vert x\Vert_2$ coincides with the Euclidean norm. 
When there is no possible confusion, sometimes
 we use the notation $\Vert\cdot\Vert$ for any equivalent matrix or vector norm.

 For any $u,v\in\Sa_{d}^+$ we recall the Ando-Hemmen inequality
\begin{equation}\label{square-root-key-estimate}
\Vert u^{1/2}-v^{1/2}\Vert \leq \left[\lambda^{1/2}_{\rm min}(u)+\lambda^{1/2}_{\rm min}(v)\right]^{-1}~\Vert u-v\Vert
\end{equation}
that holds for any unitary invariant matrix norm $\Vert\cdot \Vert$, including the spectral and the Frobenius norms, see for instance Theorem 6.2 on page 135 in~\cite{higham}, as well as Proposition 3.2  in~\cite{hemmen}.

Let us associate with some  $\varpi\in\Sa^+_d$  the increasing map $\mbox{\rm Ricc}_{\varpi}$ from $\Sa^0_d$ into  {$\Sa^+_d$} defined by
\begin{equation}\label{ricc-maps-def}
\begin{array}{lccl}
\mbox{\rm Ricc}_{\varpi}: &\Sa^0_d &\mapsto &\Sa^+_d\\
&v &\leadsto &\mbox{\rm Ricc}_{\varpi}(v):=(I+(\varpi+v)^{-1})^{-1}.\\ 
\end{array}
\end{equation}
A refined stability analysis of Riccati matrix differences $v_{n+1}:=\mbox{\rm Ricc}_{\varpi}(v_n)$ and the limiting stationary matrices $r=\mbox{\rm Ricc}_{\varpi}(r)$ associated with these maps is provided in~\cite{adm-24}.
Notice that
\begin{eqnarray}
r=\mbox{\rm Ricc}_{\varpi}(r)&\Longleftrightarrow&
r^{-1}=I+(\varpi+r)^{-1}
\Longleftrightarrow\varpi r^{-1}+I=
(\varpi+r)r^{-1}=(\varpi+r)+I\nonumber\\
&\Longleftrightarrow&\varpi r^{-1}=\varpi+r
\Longleftrightarrow r^{-1}=I+\varpi^{-1}r\Longleftrightarrow
I=r+r\,\varpi^{-1}r.\label{equiv-r}
\end{eqnarray}
Also note the increasing property
$$
v_1\leq v_2\Longrightarrow \mbox{\rm Ricc}_{\varpi}(v_1)\leq \mbox{\rm Ricc}_{\varpi}(v_2)
$$
For any $n\geq p\geq 1$ and $v\in\Sa^d_0$ we also have the estimates
\begin{equation}\label{ricc-maps-incr}
 \mbox{\rm Ricc}_{\varpi}^p(0)\leq \mbox{\rm Ricc}_{\varpi}^n(0)\leq \mbox{\rm Ricc}_{\varpi}^n(v)
\leq \mbox{\rm Ricc}_{\varpi}^{n-1}(I)\leq  \mbox{\rm Ricc}_{\varpi}^{p-1}(I)\leq I
\end{equation}
These matrix equations belong to the class of discrete algebraic Riccati equations (DARE), {and} no analytical solutions are available for general models.  The unique positive definite fixed point {of the Riccati differences} is given by 
\begin{equation}\label{def-fix-ricc-1}
(I+\varpi^{-1})^{-1}\preceq r:=-\frac{\varpi}{2}+\left(\varpi+\left(\frac{\varpi}{2}\right)^2\right)^{1/2}\preceq I.
\end{equation}
 In addition, there exists some constant $c_{\varpi}$
such that
\begin{equation}\label{cv-ricc-intro}
 \Vert v_{n}-r\Vert_2\leq c_{\varpi}~(1+\lambda_{\text{\rm min}}(\varpi+r))^{-2n}~
\Vert v_0-r\Vert_2
\end{equation}
{The contraction rates in \eqref{cv-ricc-intro} are based on Floquet-type representation of Riccati flows and they are sharp} (see Theorem 1.3 in ~\cite{dh-22}). For further discussion see ~\cite{adm-24,dh-22}, and the references therein.

   The positive and negative part of $a\in\RR$ are defined respectively by $a_+=\max(a,0)$
and  $a_-=\max(-a,0)$. Given $a,b\in\RR$ we also set $a\vee b=\max(a,b)$ and
$a\wedge b=\min(a,b)$. 
  
  \section{Some background on bridges}\label{sec:prelim}

\subsection{Overview}
  
 The theory of entropic optimal transport  
 is at the carrefour of statistics and applied probability, optimal control theory,
 machine learning and artificial intelligence. Entropic bridges have emerged as powerful techniques in a various fields of physics, economics, 
 engineering and computer sciences. Every community often has a different way
to interpret and formulate Schr\"odinger and Sinkhorn bridges and every formulation has its own advantages and limitations. To guide the reader, this section provides a comprehensive overview of some of the main formulations used in this field.
  
  Section~\ref{conjugate-formulation-sec} presents a probabilistic formulation of Sinkhorn bridges in terms of conjugate and dual integral operators. Equivalent Bayes' type formulations in terms of state space models are discussed in Section~\ref{state-space-conditioning-sec}. Section~\ref{dyn-statis-bridges-sec}  underlines the links between static and dynamic Schr\" odinger bridges in the context of discrete generation Markov processes. Forward-backward formulations for continuous time generative-type diffusion 
  models are discussed in~\ref{cont-discrete-sec}. Section~\ref{ent-sec-proj} provides an alternative analytical interpretation of Sinkhorn bridges in terms of entropic projections.

  \subsection{Conjugate and dual transitions}\label{conjugate-formulation-sec}
When
 $\delta_{x}K\ll \mu K$ for any $x\in\XX$, we have the Bayes' formula
\begin{equation}\label{Bayes1}
\mu\times K=\left((\mu K)\times K^{\ast}_{\mu}\right)^{\flat}\quad\mbox{\rm with}\quad
K^{\ast}_{\mu}(y,dx):=\mu(dx)~\frac{d\delta_x K}{d\mu K}(y).
\end{equation}
A more general definition of the dual transition $K^{\ast}_{\mu}$ is provided in Section 2 in the article~\cite{dm-03}.
Note the reversible property 
$$
\mu(dx) K(x,dy)K^{\ast}_{\mu}(y,d\overline{x})=\mu(d\overline{x}) K(\overline{x},dy)K^{\ast}_{\mu}(y,dx)
\quad\mbox{\rm and therefore}\quad
\mu K K^{\ast}_{\mu}=\mu.
$$
By symmetry arguments, when
 $\delta_{y}L\ll \eta L$ for any $y\in\YY$, we have the Bayes' formula
\begin{equation}\label{Bayes2}
\left(\eta\times L\right)^{\flat}=(\eta L)\times L_{\eta}^{\ast}
\quad\mbox{\rm with}\quad
L^{\ast}_{\eta}(x,dy):=\eta(dy)~\frac{d\delta_y L}{d\eta L}(x).
\end{equation}
Again, we have the reversible property
$$
\eta(dy) L(y,dx)L^{\ast}_{\eta}(x,d\overline{y})=\eta(d\overline{y}) L(\overline{y},dx)L^{\ast}_{\eta}(x,dy)
\quad\mbox{\rm and therefore}\quad
\eta L L^{\ast}_{\eta}=\eta.
$$

Consider some reference probability measure 
$P\in \Pi(\mu,\point)
$. We have disintegration formula $P=(\mu\times K)$, for some
Markov transition $K$ from $\XX$ into $\YY$.
Sinkhorn bridges $\Pa_n$ are defined sequentially starting from a reference probability
$$
\Pa_0:=P\Longleftrightarrow \left(\Pa_0:=\mu\times \Ka_0
\quad\mbox{\rm with}\quad
\Ka_0:=K\right).
$$
Rewritten in a slightly different way, the Bayes' rule (\ref{Bayes1}) takes the form
$$
 \Pa_0:=\mu\times \Ka_0=\left(\pi_0\times \Ka_1\right)^{\flat}
\in \Pi\left(\mu,\pi_0\right)
\quad\mbox{\rm with}\quad \pi_0:=\mu \Ka_0\quad\mbox{\rm and}\quad
\Ka_1:=(\Ka_0)^{\ast}_{\mu}.
$$
In the same vein, Bayes' rule (\ref{Bayes2}) takes the form
$$
\Pa_1:=\left(\eta\times \Ka_1\right)^{\flat}=\pi_1\times \Ka_2\in \Pi\left(\pi_1,\eta\right)
\quad\mbox{\rm with}\quad \pi_1:=\eta \Ka_1\quad\mbox{\rm and}\quad
\Ka_2:=(\Ka_1)^{\ast}_{\eta}.
$$
Iterating the above procedure we define for even indices a collection of bridges 
\begin{equation}\label{bayes2n}
 \Pa_{2n}:=\mu\times \Ka_{2n}=\left(\pi_{2n}\times \Ka_{2n+1}\right)^{\flat}
\in \Pi\left(\mu,\pi_{2n}\right)
\end{equation}
with the marginal distributions
$$
\pi_{2n}:=\mu \Ka_{2n}\quad\mbox{\rm and the transitions}\quad
\Ka_{2n+1}:=(\Ka_{2n})^{\ast}_{\mu}.
$$
Sinkhorn bridges for odd indices are given by
\begin{equation}\label{bayes2n1}
\Pa_{2n+1}:=\left(\eta\times \Ka_{2n+1}\right)^{\flat}=\pi_{2n+1}\times \Ka_{2(n+1)}\in \Pi\left(\pi_{2n+1},\eta\right)
\end{equation}
with the marginal distributions
$$
 \pi_{2n+1}:=\eta \Ka_{2n+1}\quad\mbox{\rm and the transitions}\quad
\Ka_{2(n+1)}:=(\Ka_{2n+1})^{\ast}_{\eta}.
$$

\subsection{State space models}\label{state-space-conditioning-sec}
In terms of random variables $(X,Y)$ with distribution $P=\mu\times K$
we have
$$
P(d(x,y))=\PP((X,Y)\in d(x,y))=
\PP(X\in dx)~\PP(Y\in dy~|~X=x)=\mu(dx)~K(x,dy)
$$
with the distributions
$$
\begin{array}{l}
\mu(dx):=\PP(X\in dx)\quad\mbox{and}\quad K(x,dy):=\PP(Y\in dy~|~X=x)\\
\\
\Longrightarrow
(\mu K)(dy)=\PP(Y\in dy)\quad\mbox{and}\quad
K^{\ast}_{\mu}(y,dx)=\PP(X\in dx~|~Y=y).
\end{array}
$$
For continuous random variables on $\XX=\RR^p=\YY$   with conditional densities
with respect to the Lebesgue measure we have
$$
\begin{array}{l}
\mu(dx)=p_X(x)dx \quad\mbox{and}\quad K(x,dy)=p_{Y|X}(y|x)dy\\
\\
\Longrightarrow
P(d(x,y))=p_{X,Y}(x,y)~dx dy:=p_X(x)p_{Y|X}(y|x)~dxdy
\end{array}$$
where $p_{X,Y}(x,y)$ is the joint density of $(X,Y)$, 
$p_X(x)$ the marginal density of $X$ and $p_{Y|X}(y|x)$ is the conditional density of 
$Y$ given $X=x$.  In this context, Bayes' rule (\ref{Bayes1}) takes the form
$$
p_{X,Y}(x,y)=p_X(x)p_{Y|X}(y|x)=p_Y(y)p_{X|Y}(x|y)
\quad\mbox{and}\quad
K^{\ast}_{\mu}(y,dx)=p_{X|Y}(x|y)~dx.
$$
 \subsubsection*{Sinkhorn bridges}
Consider random variables $(X_n,Y_n)$ such that
$$
\Pa_{n}(d(x,y))=\PP((X_n,Y_n)\in d(x,y))
$$
with the Sinkhorn bridges $\Pa_n$ defined in (\ref{def-Pa-n}).
For even indices we have
$$
\begin{array}{l}
 \Pa_{2n}:=\mu\times \Ka_{2n}\\
 \\
 \Longleftrightarrow
\mu(dx)= \PP(X_{2n}\in dx)\quad \mbox{\rm and}\quad
 \Ka_{2n}(x,dy)=\PP(Y_{2n}\in dy~|~X_{2n}=x).
\end{array}$$
In this context, Bayes' rule (\ref{bayes2n}) holds with 
\begin{equation}\label{bayes2nb}
 \pi_{2n}(dy)=\PP(Y_{2n}\in dy)
\quad \mbox{\rm and}\quad
  \Ka_{2n+1}(y,dx)=\PP(X_{2n}\in dx~|~Y_{2n}=y).
\end{equation}
For odd indices we have
$$
\begin{array}{l}
 \Pa_{2n+1}:=\left(\eta\times \Ka_{2n+1}\right)^{\flat}\\
 \\
 \Longleftrightarrow
\eta(dy)= \PP(Y_{2n+1}\in dy)\quad \mbox{\rm and}\quad
 \Ka_{2n+1}(y,dx)=\PP(X_{2n+1}\in dx~|~Y_{2n+1}=y).
\end{array}
$$
Again, 
Bayes' rule (\ref{bayes2n1}) holds with 
\begin{equation}\label{bayes2n1b}
 \pi_{2n+1}(dx)=\PP(X_{2n+1}\in dx)
\quad \mbox{\rm and}\quad
  \Ka_{2(n+1)}(x,dy)=\PP(Y_{2n+1}\in dy~|~X_{2n+1}=x).
\end{equation}

 \subsubsection*{Density models}
 
Consider locally compact Polish spaces $(\XX,\lambda)$ and $(\YY,\nu)$ endowed with some locally bounded positive measures $\lambda$ and $\nu$. Also assume the marginal probability measures are given by
\begin{equation}\label{ref-intro-UV}
\mu(dx)=\lambda_U(dx):=e^{-U(x)}~\lambda(dx)\quad \mbox{\rm and}\quad
\eta(dy)=\nu_V(dy):=e^{-V(y)}~\nu(dy)
\end{equation}
for some normalized potential functions $U:x\in\XX\mapsto U(x)\in \RR$ and $V:y\in\XX\mapsto V(y)\in \RR$.
We have
$$
\PP((X_n,Y_n)\in d(x,y))=p_{(X_n,Y_n)}(x,y)~\gamma(d(x,y))
\quad \mbox{\rm with}\quad
 \gamma:=\lambda\otimes\nu.
$$
 In addition
\begin{equation}\label{rec-bayes-marg}
\begin{array}{rclcrcl}
p_{(X_{2n},Y_{2n})}(x,y)&=&p_{X_{2n}}(x)~p_{Y_{2n}|X_{2n}}(y|x)& \mbox{\rm with}&
p_{X_{2n}}(x)&=&e^{-U(x)}\\
&&&&&&\\
p_{(X_{2n+1},Y_{2n+1})}(x,y)&=&p_{Y_{2n+1}}(y)~p_{X_{2n+1}|Y_{2n+1}}(x|y)& \mbox{\rm with}&
p_{Y_{2n+1}}(y)&=&e^{-V(y)}.
\end{array}
\end{equation}
In this scenario, the recursions in  (\ref{bayes2nb})  and  (\ref{bayes2n1b}) take the form
\begin{equation}\label{rec-bayes}
p_{X_{2n+1}|Y_{2n+1}}:=p_{X_{2n}|Y_{2n}}\quad
\mbox{\rm and}\quad
p_{Y_{2(n+1)}|X_{2(n+1)}}:=p_{Y_{2n+1}|X_{2n+1}}.
\end{equation}

 \subsection{Dynamic and static bridges}\label{dyn-statis-bridges-sec}
 Consider a Markov chain $Z_k$ on $\XX$ with transition semigroup
 $ M_{l,k}$ with $l\leq k$ defined by
 $$
 M_{l,k}(z_{l},dz_k):=\PP(Z_k\in dz_k~|~Z_{l}=z_{l}).
 $$
We fix the time horizon $n$ and observe that 
$$
\begin{array}{l}
\mu(dz_0):=\PP(Z_0\in dz_0)\\
\\
\Longrightarrow
P(d(z_0,z_n)):=\PP((Z_0,Z_n)\in d(z_0,z_n))=(\mu \times K)(d(z_0,z_n))\quad \mbox{\rm with}\quad K=M_{0,n}.
\end{array}
$$

 \subsubsection*{Pinned  process}
 
The conditional distribution of the random path $(Z_1,\ldots,Z_{n-1})$ given the initial and terminal states $(Z_0,Z_n)$ (we call this a pinned Markov chain)
is given by $$
\begin{array}{l}
\displaystyle \PP\left( (Z_1,\ldots,Z_{n-1})\in
(d(z_1,\ldots,z_{n-1}))~|~(Z_0,Z_n)=(z_0,z_n)\right)\\
\\
\displaystyle= P_{0\rightarrow n}^{(z_0,z_n)}(d(z_1,\ldots,z_{n-1})):=\prod_{1\leq k<n}M_{k-1,k}^{(n,z_n)}(z_{k-1},dz_{k}).
\end{array}
$$
with the collection of Markov transitions
$$
M_{k-1,k}^{(n,z_n)}(z_{k-1},dz_{k}):=
\frac{M_{k-1,k}(z_{k-1},dz_k)~M_{k,n}(z_k,dz_n)}{
M_{k-1,n}(z_{k-1},dz_n)}.
$$
Note that the above transitions do not depend on $\mu$ and we have
\begin{eqnarray*}
\Pb(d(z_0,\ldots,z_{n}))&:=&\PP\left( (Z_0,\ldots,Z_{n})\in
(d(z_0,\ldots,z_{n}))\right)\\
&=&P(d(z_0,z_n))\times P_{0\rightarrow n}^{(z_0,z_n)}(d(z_1,\ldots,z_{n-1})).
\end{eqnarray*}

Let $\Pi_n(\mu,\eta)$ be the set of probability measures on the path space $\XX^{n+1}:=\XX\times\XX^n$
with marginals $\mu$ and $\eta$ at time $k=0$ and $k=n$.
The dynamic Sch\"odinger bridge between $\mu$ and $\eta$ with reference distribution $\Pb$ is given by the bridge map
$$
\Pb_{\mu,\eta}:=\argmin_{\Qb\,\in\, \Pi_n(\mu,\eta)}\Ha(\Qb~|~\Pb).
$$
Any probability measure $\Qb\in \Pi_n(\mu,\eta)$
can be disintegrated w.r.t. the marginal distribution $Q(d(z_0,z_n))$ of the initial and terminal states, namely 
$$
\Qb(d(z_0,\ldots,z_{n}))
=Q(d(z_0,z_n))\times Q_{0\rightarrow n}^{(z_0,z_n)}(d(z_1,\ldots,z_{n-1})).
$$
This yields the entropy factorization
$$
\Ha(\Qb~|~\Pb)=
\Ha(Q~|~P)+\int~Q_n(d(z_0,z_n))~ \Ha(Q_{0\rightarrow n}^{(z_0,z_n)}|P_{0\rightarrow n}^{(z_0,z_n)})
$$
from which we conclude~\cite{adm-24,chen-phd,essid,follmer,leonard} that the static Schr\" odinger  bridge $P_{\mu,\eta}$ and dynamic Schr\" odinger bridge $\Pb_{\mu,\eta}$ are connected by
the formulae
$$
\Pb_{\mu,\eta}(d(z_0,\ldots,z_{n})))=P_{\mu,\eta}(d(z_0,z_n))\times P_{0\rightarrow n}^{(z_0,z_n)}(d(z_1,\ldots,z_{n-1})).
$$

 \subsubsection*{Backward conditioning formulae}
To simplify notation we consider a time homogenous Markov transition $M_{k-1,k}=M$. In this context, 
we have the backward conditioning formulae
$$
\PP\left( (Z_0,\ldots,Z_{n-1})\in
d(z_0,\ldots,z_{n-1})~|~Z_n=z_n\right)=\prod_{n>k\geq 0}~M^{\ast}_{\mu_k}(z_{k+1},dz_{k})
$$
the dual transitions
$$
M^{\ast}_{\mu_k}(z_{k+1},dz_{k}):=\frac{\mu_{k}(dz_{k})~M(z_{k},dz_{k+1})}{(\mu_{k}M)(dz_{k+1})}\quad \mbox{\rm and}\quad
\mu_k(dz)=\PP(Z_k\in dz)=(\mu_{k-1}M)(dz).
$$
This yields the conjugate formula
$$
(\mu,K):=(\mu_0,M_{0,n})\Longrightarrow
K^{\ast}_{\mu}(z_n,dz_0)=
\PP\left( Z_0\in
dz_0~|~Z_n=z_n\right)=
\frac{\mu(dz_{0})~K(z_0,dz_{n})}{(\mu K)(dz_0)}.
$$

 \subsection{Continuous vs discrete time models}\label{cont-discrete-sec}
We illustrate the conjugate formulae presented in Section~\ref{conjugate-formulation-sec} and Section~\ref{state-space-conditioning-sec} in the context of diffusion flows, due to the practical relevance of such methods. 

Consider the stochastic flow
 $s\in [0,t]\mapsto Y^{(0)}_{s}(x)$ on some {\it fixed time horizon $t$} starting at $Y^{(0)}_{0}(x)=x\in\XX=\RR^d=\YY$ defined by the diffusion
\begin{equation}\label{def-diff}
dZ^{(0)}_{s}(x)=b^{(0)}_s(Z^{(0)}_{s}(x))ds+\sigma_s(Z^{(0)}_{s}(x))^{1/2}~dB_s
\end{equation}
for some drift function $b_s(x)$, some diffusion function $\sigma_s(x)\in\Sa^+_d$ and some $d$-dimensional Brownian motion $B_s$. Formally, we may slice the time interval $[0,t]_{\Delta}:=\{s_0,\ldots,s_{n}\}$ via some time mesh $s_{i+1}=s_i+\Delta$ from $s_0=0$ to $s_n=t=n\Delta$. In this notation, for $s=s_k\in [0,t]_{\Delta}$ with $k<n$ we compute $Z^{(0)}_{s+\Delta}(x)$ from $Z^{(0)}_{s}(x)$ using the formula
$$
Z^{(0)}_{s+\Delta}(x)-Z^{(0)}_{s}(x)=b^{(0)}_s(Z^{(0)}_{s}(x))~\Delta+\sigma_s(Z^{(0)}_{s}(x))^{1/2}~(B_{s+\Delta}-B_s)
$$
i.e.~via Euler-Maruyama time discretization.

Given some initial random variable $X_0$ with distribution $\mu$ independent of the Brownian motion $B_s$ we set
$$
Y_0:=Z^{(0)}_{t}(X_0)\Longrightarrow \Pa_0(d(x,y)):=
\PP((X_0,Y_0)\in d(x,y))=\mu(dx)~\Ka_0(x,dy)
$$
with the Markov transition
$$
\Ka_0(x,dy):=\PP(Z^{(0)}_{t}(x)\in dy):=q^{(0)}_t(x,y)~dy.
$$
Observe that
\begin{eqnarray*}
(\mu \Ka_0)(dy)&=&\PP(Y_0\in dy)\\
&=&\PP(Z^{(0)}_{t}(X_{0})\in dy)=p^{(0)}_t(y)~dy
\quad\mbox{\rm with $\forall s\in [0,t]$}\quad p^{(0)}_{s}(y):=\int \mu(dx) ~q^{(0)}_s(x,y).
\end{eqnarray*}
The conjugate transition $\Ka_1$ is associated with the diffusion flow $s\in [0,t]\mapsto Z^{(1)}_{s}(y)$ starting at $Z^{(1)}_{t}(y)=y$   defined by the backward stochastic differential equation
$$
-d_s Z^{(1)}_{s}(y)=b^{(1)}_{s}(Z^{(1)}_{s}(y))~ds+\sigma_s(Z^{(1)}_{s}(y))^{1/2}~dB_s
$$
 with the drift function
$$
b^{(1)}_{s}:=-b^{(0)}_s+\nabla_{\sigma_s} \log p^{(0)}_s
\quad\mbox{\rm and}\quad 
~(\nabla_{\sigma_s} \log p^{(0)}_s)^j:=\frac{1}{p^{(0)}_s}~\sum_{1\leq i\leq d}\partial_{x_i}\left(\sigma_s^{j,i} p^{(0)}_s\right).
$$
Formally, slicing the interval as above, for $s=s_k\in [0,t]_{\Delta}$ with $k>0$ we compute $Z^{(0)}_{s-\Delta}(x)$ from $Z^{(0)}_{s}(x)$ using the backward formula
$$
Z^{(1)}_{s-\Delta}(y)-Z^{(1)}_{s}(y)
=b^{(1)}_{s}(Z^{(1)}_{s}(y))~\Delta+\sigma_s(Z^{(1)}_{s}(y))^{1/2}~(B_{s}-B_{s-\Delta})
$$
We  have the dual formulae
$$
\begin{array}{l}
\mu(dx)~\Ka_0(x,dy)=(\mu \Ka_0)(dy)~\Ka_1(y,dx)
\quad\mbox{\rm and}\quad \Ka_1(y,dx):=\PP(Z^{(1)}_{0}(y)\in dx).
\end{array}
$$
A detailed proof of the above backward conjugate formula can be found in~\cite{anderson,anderson-dp,haussman}. 
At the next iteration, we consider
some initial random variable $Y_1$ with distribution $\eta$ independent of the Brownian motion $B_s$ and we set
$$
X_1:=Z^{(1)}_{0}(Y_1)\Longrightarrow \Pa_1(d(x,y)):=
\PP((X_1,Y_1)\in d(x,y))=\eta(dy)~\Ka_1(y,dx).
$$

\subsection{Some half-bridge entropy formulae}\label{ent-sec-proj}
 
We have the entropy formulae
$$
\Ha((\eta\times L)^{\flat}~|~\mu\times K)=\Ha(\eta\times L~|~(\mu K)\times K^{\ast}_{\mu})=
\Ha((\eta L)\times L_{\eta}^{\ast}~|~\mu\times K).
$$
This yields the variational formulae
\begin{equation}\label{Bayes1ee}
 (\eta\times K^{\ast}_{\mu})^{\flat}=\argmin_{Q\in \Pi(\point,\eta)} \Ha(Q~|~\mu\times K)
\quad \mbox{\rm and}\quad
\mu\times L_{\eta}^{\ast}=\argmin_{Q\in \Pi(\mu,\point)} \Ha((\eta\times L)^{\flat}~|~Q).
\end{equation}
By symmetry arguments, we also have the entropy formulae
$$
\Ha(\mu\times K~|~(\eta\times L)^{\flat})=
\Ha((\mu K)\times K^{\ast}_{\mu}~|~\eta\times L)=
\Ha(\mu\times K~|~(\eta L)\times L_{\eta}^{\ast})
$$
This yields the variational   formulae
\begin{equation}\label{Bayes2ee}
(\eta\times K^{\ast}_{\mu})^{\flat}=
\argmin_{Q\in \Pi(\point,\eta)}\Ha(\mu\times K~|~Q)
\quad\mbox{\rm and}\quad
\mu\times L_{\eta}^{\ast}=
\argmin_{Q\in \Pi(\mu,\point)}\Ha(Q~|~(\eta\times L)^{\flat}).
\end{equation}

Applying the left hand side (l.h.s.) of (\ref{Bayes1ee}) and (\ref{Bayes2ee}) to $K=\Ka_{2n}$
 we have the half-bridge equivalent formulation
$$
\Pa_{2n+1}:= \argmin_{Q\in \Pi(\point,\eta)}\Ha(Q~|~\Pa_{2n})=\argmin_{Q\in \Pi(\point,\eta)}\Ha(\Pa_{2n}~|~Q).
$$
Similarly, applying the right hand side (r.h.s.) of (\ref{Bayes1ee}) and (\ref{Bayes2ee}) to $L=\Ka_{2n+1}$
 we  also have
$$
\Pa_{2(n+1)}:= \argmin_{Q\in \Pi(\mu,\point)}\Ha(\Pa_{2n+1}~|~Q)=\argmin_{Q\in \Pi(\mu,\point)}\Ha(Q~|~\Pa_{2n+1}).
$$

\section{Stability of Markov operators}\label{sec:stab_markov}

\subsection{Overview}
  
The time varying Gibbs-type semigroup formulation of Sinkhorn bridges 
(\ref{fb-sinkhorn}) and the entropy formula (\ref{kdec2n-intro}) show that the convergence analysis of Sinkhorn bridges rely on the stability analysis of a class of time varying Markov processes  sharing a common invariant target measure at each individual step.  This section gives a summary on the stability analysis of Markov operators. Section~\ref{sec-criteria} and Section~\ref{kanto-sect} dicuss the different classes of entropy, divergences and the Kantorovich criteria used in the article.
  Section~\ref{contract-coef-sec} review an operator-theoretic framework based on contraction coefficients and Lyapunov principles. Section~\ref{sect-transport-ineq} is dedicated to transportation cost inequalities including local log-Sobolev and Talagrand quadratic cost inequalities.

  \subsection{Entropy and divergence criteria}\label{sec-criteria}
A  measure $\mu_1$ on $\XX$ is said to be absolutely continuous with respect to another measure $\mu_2$ on $\XX$, and we write $\mu_1\ll \mu_2$, if $\mu_2(A)=0$ implies $\mu_1(A)=0$ for any measurable subset $A\subset \XX$. Whenever $\mu_1\ll \mu_2$, we denote by ${d\mu_1}/{d\mu_2}$ the Radon-Nikodym derivative of $\mu_1$ w.r.t. $\mu_2$.
Consider a convex function 
\begin{equation}\label{Phi-ref}
\begin{array}{cccl}
    \Phi:&\RR_+^2 &\mapsto &\RR\cup\{+\infty\}\\
         &(u,v)   &\leadsto&\Phi(u,v)\quad\mbox{\rm s.t.  $\forall a\in\RR_+$~~~ $\Phi( a u,a v)=a~ \Phi(u,v)$ and  $\Phi(1,1)=0$.}
\end{array}
\end{equation}
The $\Phi$-entropy between two probability measures $\mu_1,\mu_2$ on $\XX$ is defined for any dominating measure $\gamma$ (such that $\gamma\gg \mu_1$ and $\gamma\gg\mu_2$) by the formula
\begin{equation}\label{def-DPhi}
\Ha_{\Phi}(\mu_1,\mu_2):=\int~\Phi\left(\frac{d\mu_1}{d\gamma}(x),\frac{d\mu_2}{d\gamma}(x)\right)~\gamma(dx).
\end{equation}
By homogeneity arguments, the definition (\ref{def-DPhi}) does not depend on the choice of the dominating measure $\gamma$. Note that
 \begin{equation}\label{tv-osc}
 \text{if} \quad \Phi(u,v)=\Phi_0(u,v):=\frac{1}{2}~\vert u-v\vert \quad\text{then}\quad \Ha_{\Phi_0}(\mu_1,\mu_2)=\Vert \mu_1-\mu_2\Vert_{\sf tv}.
 \end{equation}
In the above display, $\Vert \mu_1-\mu_2\Vert_{\sf tv}$ stands for the total variation distance between the probability measures $\mu_1$ and $\mu_2$ defined by
\begin{eqnarray}
\Vert \mu_1-\mu_2\Vert_{\sf tv}&:=&\frac{1}{2}~\sup{\left\{ \vert(\mu_1-\mu_2)(f)\vert~:~f\in\Ba(\XX)~;~\Vert f\Vert \leq 1 \right\}}\nonumber\\
&=&\sup{\left\{ \vert(\mu_1-\mu_2)(f)\vert~:~f\in\Ba(\XX)~;~\mbox{\rm osc}(f) \leq 1 \right\}}\label{tv-osc-norm}
\end{eqnarray}
with the oscillation $\mbox{\rm osc}(f)$ of the function $f\in \Ba(\XX)$ defined by
$$
\mbox{\rm osc}(f):=\sup{\left\{|f(x_1)-f(x_2)|~:~(x_1,x_2)\in \XX^2 \right\}}.
 $$

Whenever $\Phi(1,0)=\infty$, the $\Phi$-entropy coincides with the $\phi$-divergence in the sense of Csisz\'ar~\cite{csiszar}, i.e., we also have
 $$
 \Ha_{\Phi}(\mu_1,\mu_2)=\int~\phi\left(\frac{d\mu_1}{d\mu_2}(x)\right)~\mu_2(dx)\quad \text{where} \quad  \phi(u):=\Phi(u,1)
 $$
 and we have used the convention $\Ha_{\Phi}(\mu_1,\mu_2)=\infty$ when $\nu_1\not\ll \nu_2$. Therefore, the $\Phi$-entropy generalizes the classical relative entropy (a.k.a.~Kullback-Leibler divergence or $I$-divergence). More precisely, choosing the functions
    $$
     \Phi(u,v)=u\log{(u/v)}\quad \mbox{\rm and}\quad \phi(u)= u\log{u}\
    $$
we recover the classical formulae
 $$
 \Ha_{\Phi_1}(\mu_1,\mu_2)=
 \Ha(\mu_1~|~\mu_2):=\mu_1\left(\log{\left(\frac{d\mu_1}{d\mu_2}\right)}\right).
$$
The $\phi$-divergence criteria discussed above also includes $\LL_p$-norms,  Havrda-Charvat entropies,
$\alpha$-divergences, Jensen-Shannon divergence, squared He\-llin\-ger integrals and Kakutani-Hel\-inger,  as well as Jeffreys and R\'enyi divergences, see for instance~\cite{csiszar,dm-03,majernik,nielsen} and references therein.

  \subsection{Kantorovich semi-distances}\label{kanto-sect}
  
Consider a  semi-distance $\varphi$ on a Polish space $\XX$; a semi-distance is essentially a distance function that is l.s.c.~except that the triangle inequality may be violated. 
Some examples of semi-distances that we consider include the discrete metric $\varphi_d(x_1,x_2):=1_{x_1\neq x_2}$ and, for a given
 l.s.c.~function $g$, lower bounded away from zero,  the weighted discrete metric defined by the formulae
 \begin{align}\label{refvarpiV}
\varphi_g(x_1,x_2):=\varphi_d(x_1,x_2)~w_g(x_1,x_2)\quad \mbox{\rm with}\quad w_g(x_1,x_2)=g(x_1)+g(x_2).
 \end{align}

The Kantorovich semi-distance between  two measures $\mu_1,\mu_2\in{\sf Prob}(\XX)$ is defined by
 \begin{equation}\label{ref-cost}
D_{\varphi}(\mu_1,\mu_2):=\inf_{Q\in \Pi(\mu_1,\mu_2)}{Q(\varphi)}.
\end{equation}
The assumption that  $\XX$ is a complete and separable metric space and the fact that $\varphi$ is l.s.c.~and non-negative ensure the existence of an optimal transport plan \cite[Theorem 4.1]{villani-2}.  Moreover, the assumption that $\varphi(x_1,x_2)=0\Leftrightarrow x_1=x_2$ ensure that $D_{\varphi}(\mu_1,\mu_2)=0$ if and only if $\mu_1=\mu_2$.   
The Kantorovich semi-distance   is related to several well-known measures of distance between probability distributions. Specifically, when $(\XX,\varphi)$ is a Polish space
  the quantity   $D_{\varphi}(\mu_1,\mu_2)$ is the  Kantorovitch distance (or Wasserstein 1-distance) between $\mu_1$ and $\mu_2$,  whilst for $p\geq 1$ the Wasserstein $p$-distance $\Da_{\varphi,p}(\mu_1,\mu_2)$    is defined by  $\Da_{\varphi,p}(\mu_1,\mu_2)=D_{\varphi^p}(\mu_1,\mu_2)^p$.  
 We also have the formulae 
  $$
  D_{\varphi_d}(\mu_1,\mu_2)=\|\mu_1-\mu_2\|_{\sf tv}\quad \mbox{\rm as well as}\quad
  D_{\varphi_g}(\mu_1,\mu_2)=\vertiii{ \mu_1-\mu_2}_g
  $$
  with 
  the $g$-weighted total variation norm 
\begin{eqnarray}
\vertiii{ \mu_1-\mu_2}_g&:=&\vert \mu_1-\mu_2\vert(g)\nonumber\\
&=&\sup{\left\{ \vert(\mu_1-\mu_2)(f)\vert~:~f\in\Ba_g(\XX)~;~\Vert f\Vert_g \leq 1 \right\}}\nonumber\\
&=&\sup{\left\{ \vert(\mu_1-\mu_2)(f)\vert~:~f\in\Ba_g(\XX)~;~\mbox{\rm osc}_g(f) \leq 1 \right\}},\label{kr2-dual}
\end{eqnarray}
where $\mbox{\rm osc}_g(f)$ stands for the $\varphi_g$-oscillation  of the function $f\in \Ba_g(\XX)$ defined by the Lipschitz constant
$$
\mbox{\rm osc}_g(f):=\sup{\left\{|f(x_1)-f(x_2)|/\varphi_g(x_1,x_2)~:~(x_1,x_2)\in \XX^2~;~x_1\not=x_2 \right\}}.
 $$
 
  When $\XX=\RR^d$ and $\varphi(x_1,x_2)=\Vert x_1-x_2\Vert_2$ is the Euclidian norm we simplify notation and we write $\Da_{p}$ instead of $\Da_{\varphi,p}$; that is
\begin{equation}\label{def-c2}
\Da_{p}(\mu_1,\mu_2):= \inf_{Q \in \Pi(\mu_1,\mu_2)}~\left(\int \Vert x_1-x_2\Vert_{2}^p~Q(d(x_1,x_2))\right)^{1/p}.
\end{equation}
 For a more thorough discussion on Kantorovich criteria  and 
weighted total variation norms we refer to ~\cite{dmg-25,dm-penev-2017} and references therein.

  \subsection{Contraction coefficients}\label{contract-coef-sec}  

  \subsubsection{Dobrushin coefficient}  
Given a convex function  $\Phi$ satisfying (\ref{Phi-ref}), the Dobrushin $\Phi$-contraction coefficient ${\sf dob}_{\Phi}(\Ka)$ of  some Markov transition $K(x,dy)$  from $\XX$ into $\YY$  is defined by 
\begin{equation}\label{defi-beta-Phi}
{\sf dob}_{\Phi}(K):=\sup_{}\frac{\Ha_{\Phi}(\mu_1K,\mu_2K)}{\Ha_{\Phi}(\mu_1,\mu_2)},
\end{equation}
where the supremum is taken over all probability measures 
$\mu_1,\mu_2$ such that
$\Ha_{\Phi}(\mu_1,\mu_2)>0$. 
 Here, again, the terminology $\Phi$-Dobrushin coefficient comes from the fact that choosing $\Phi(u,v)=\Phi_0(u,v)=\frac{1}{2}|u-v|$ we recover the standard Dobrushin contraction coefficient, i.e., 
\begin{eqnarray}
 {\sf dob}_{\Phi_0}(K)= {\sf dob}(K)&=&\sup_{(x_1,x_2)\in\XX^2}\Vert \delta_{x_1}K-\delta_{x_2}K\Vert_{\sf tv}\label{defi-dob-beta-Phi0}\\
 &=&\sup{\{\mbox{\rm osc}(K(f))~:~f\in\Ba(\YY)~;~\mbox{\rm osc}(f) \leq 1 \}}. \nonumber
\end{eqnarray} 
 We also recall that the Dobrushin coefficient 
 is a universal contraction coefficient, that is we have
\begin{equation}\label{defi-beta-Phi-Phi0}
  {\sf dob}_{\Phi}(K)\leq   {\sf dob}(K):={\sf dob}_{\Phi_0}(K).
\end{equation}
 A proof of the above estimate is provided in~\cite{cohen-0} for finite state spaces and in~\cite{dm-03} for general measurable spaces. For further details on the Dobrushin contraction coefficient we refer the reader to  \cite{dg-01,dm-03,dob-56} and \cite[Chapter 4]{dm-04}.

  \subsubsection{Kantorovich Lipschitz norm}  
 Consider a Markov integral operator $K$ and some other  l.s.c.~function  $h$,  lower bounded away from zero on $\YY$ and such that $\Vert K(h)\Vert_g=\Vert K(h)/g\Vert<\infty$.  This condition ensures that the right-action linear operator
$f\in \Ba_{h}(\YY)\mapsto K(f)\in \Ba_{g}(\XX)$ is continuous with operator 
norm defined by
\begin{equation}\label{P-VW-norm}
\vertiii{ K}_{\sf g,h}:=
\Vert K(h)/g\Vert=\sup{\left\{ \Vert K(f)\Vert_g~:~f\in\Ba_h(\YY)~;~\Vert f\Vert_h\leq 1\right\}}<\infty.
\end{equation}
The Lipschitz constant or Kantorovich norm \cite{dob-96} of the integral left-action operator 
 $\mu\in {\sf Prob}_g(\XX)\mapsto
 \mu K\in {\sf Prob}_h(\YY)$ is defined by
\begin{align}\label{V2W-ref}
{\sf lip}_{g,h}(K):=\sup
\frac{D_{\varphi_h}(\mu_1K,\mu_2K)}{D_{\varphi_g}(\mu_1,\mu_2)}=\sup
\frac{\vertiii{\mu_1K-\mu_2K}_h}{\vertiii{\mu_1-\mu_2}_g}
\end{align}
where the  supremum is taken over all $ (\mu_1,\mu_2)\in {\sf Prob}_g(\XX)^2$ such that $\vertiii{\mu_1-\mu_2}_g>0$. We have the equivalent dual formulation
\begin{equation}\label{P-V-norm-osc-equiv}
{\sf lip}_{g,h}(K)=\sup{\{\mbox{\rm osc}_{g}(K(f))~:~f\in\Ba_h(\YY)~;~\mbox{\rm osc}_{h}(f) \leq 1 \}}=\sup
\frac{D_{\varphi_h}(\delta_{x_1}K,\delta_{x_2}K)}{D_{\varphi_g}(\delta_{x_1},\delta_{x_2})}
\end{equation}
where the  supremum is taken over all states $ (x_1,x_2)\in \XX^2$ such that $\varphi_g(x_1,x_2)>0$. For a more detailed discussion on these equivalent formulations see \cite{dmg-25} or \cite[Section 8.2]{dm-penev-2017}. 

  \subsubsection{Lyapunov methodologies}  
Consider a Markov transition $L(y,dx)$ from $\YY$ into $\XX$. Assume that
\begin{equation}\label{Lyap-KL}
K(h)\leq \epsilon~ g+c\quad \mbox{\rm and}\quad
L(g)\leq \epsilon~ h+c
\end{equation}
for some constants $\epsilon\in ]0,1[$ and $c>0 $ and some l.s.c.~functions $g,h$  lower bounded away from zero. 
Also assume there exists  a map $\iota:l\in [l_0,\infty[\mapsto \iota(l)\in ]0,1[$, for some  $l_0>0$, such that for any $w_g(x_1,x_2)\leq l$ and $w_h(y_1,y_2)\leq l$ we have
\begin{equation}\label{loc-contracKL}
\Vert \delta_{x_1} K-\delta_{x_2} K\Vert_{\sf tv}\vee \Vert \delta_{y_1} L-\delta_{y_2} L\Vert_{\sf tv}\leq 1-\iota(l).
\end{equation}

Condition (\ref{loc-contracKL}) is equivalent to the local minorization condition
\begin{equation}\label{loc-min}
\delta_{x_1}K\wedge\delta_{x_2}K\geq \iota(l)~\eta\quad \mbox{\rm and}\quad 
\delta_{y_1}L\wedge  \delta_{y_2}L\geq \iota(l)~\mu
\end{equation}
for some probability measures $(\mu,\eta)$ on  $(\XX,\YY)$ that may depend on the parameters $l$ and the states $(x_1,x_2)$ and $(y_1,y_2)$  in the sub-level sets
$w_g(x_1,x_2)\leq l$ and $w_h(y_1,y_2)\leq l$ (cf. for instance~\cite{ado-24,dhj,dm-penev-2017}).
We have
$$
\{w_g\leq l\}\subset \{g\leq l\}^2\quad \mbox{\rm and}\quad
\{w_h\leq l\}\subset \{h\leq l\}^2.
$$
Thus, (\ref{loc-min}) is met as soon as
\begin{equation}\label{loc-min2}
\forall (x,y)\in (\{g\leq l\}\times \{h\leq l\})\qquad
\delta_{x}K\geq \iota(l)~\eta\quad \mbox{\rm and}\quad 
\delta_{y}L\geq \iota(l)~\mu.
\end{equation}

We recall the notation $\Ba_{\infty}(\XX)$ and $\Ba_{0}(\XX)$ from section \ref{sec:notation}.
It should be noted that, in (\ref{Lyap-KL}), the functions $g,h$ are not necessarily proper Lyapunov functions,  in the sense that $g$ and $h$ are not assumed to belong to  the sub-algebrae $\Ba_{\infty}(\XX)$ and  $\Ba_{\infty}(\YY)$ 
(recall the collection of locally upper bounded and uniformly positive functions  
 that growth at infinity (in the sense that  their sub-level sets are compact)). 
When $g,h$ are bounded (\ref{loc-min}) are minorization conditions on the whole state spaces.  Conversely, when $(g,h)\in (\Ba_{\infty}(\XX)\times\Ba_{\infty}(\YY))$ condition (\ref{loc-min}) are local minorization conditions on compact sub-level sets.  

We note that if a function $h>0$ satisfies
\begin{equation}\label{vanish-0}
K(h)/g\leq \theta \in \Ba_0(\XX)
\end{equation}
then,  for any $0<\epsilon<\Vert\theta \Vert\wedge 1$, we have
\begin{align}\label{eq:Lap}
K(h)\leq \epsilon~1_{\{\theta< \epsilon\}}~g+1_{\{\theta\geq \epsilon\}}~\theta~g\leq \epsilon~g+c\quad\mbox{\rm with}\quad 
c:=\Vert \theta\Vert~\sup_{x\in\{\theta\geq \epsilon\}}g(x).
\end{align}
Consequently, the l.h.s.~drift condition in (\ref{Lyap-KL}) is satisfied.

Note the rescaling properties
\begin{equation}\label{over-gh}
\begin{array}{l}
\displaystyle\overline{g}:=\frac{1}{2}+\frac{\epsilon}{2c}~g\quad
\quad \mbox{\rm and}\quad
\overline{h}:=\frac{1}{2}+\frac{\epsilon}{2c}~h\\
\\
\displaystyle\Longrightarrow\quad
K(\overline{h})\leq \epsilon~ \overline{g}+\frac{1}{2}\quad \mbox{\rm and}\quad
L(\overline{g})\leq \epsilon~ \overline{h}+\frac{1}{2}.
\end{array}
\end{equation}
In addition, for any $a\geq 0$, we also have
$$
\frac{1}{2}+a~ \overline{g}=\frac{a+1}{2}\left(1+~\frac{\epsilon}{c}~\frac{a}{a+1}~g\right)
\quad \mbox{\rm and}\quad
\frac{1}{2}+a~ \overline{h}=\frac{a+1}{2}\left(1+~\frac{\epsilon}{c}~\frac{a}{a+1}~h\right).
$$
The following theorem is a slight variation of contraction estimates presented in~\cite{ado-24,dmg-25} (itself based on the weighted norm contraction analysis presented in~\cite{dhj,dm-penev-2017}).  
\begin{theo}\label{theo-lipa}
Assume (\ref{Lyap-KL}) and (\ref{loc-contracKL}) are satisfied.  Then there exists some $a>0$ and $\varrho\in ]0,1[$ that only depends on the parameters $(\epsilon,c)$ and the function $\iota(r)$ 
such that
$$
\begin{array}{l}
\displaystyle
{\sf lip}_{g_a,h_a}(K)\vee {\sf lip}_{h_a,g_a}(L)\leq \varrho
\quad
\mbox{with}\quad
 g_{a}:=\frac{1}{2}+a~ g\quad
\mbox{and}\quad
h_{a}:=\frac{1}{2}+a~ h,
\end{array}$$
and $( \overline{g}, \overline{h})$ as in (\ref{over-gh}). In addition, we have
$$
{\sf lip}_{g_a,g_a}(KL)\vee
{\sf lip}_{h_a,h_a}(LK)
\leq {\sf lip}_{g_a,h_a}(K)~ {\sf lip}_{h_a,g_a}(L)\leq \varrho^2.
$$
\end{theo}

Observe that
$$
\varphi_{g_a}(x_1,x_2):=\varphi_d(x_1,x_2)~\left(1+a~w_g(x_1,x_2)\right).
$$
Assume there exists a semi-distance $ \psi$ on $\XX$ such that
$$
\begin{array}{l}
\psi(x_1,x_2)\leq w_g(x_1,x_2)
\Longrightarrow a~ \psi(x_1,x_2)\leq \varphi_{g_a}(x_1,x_2)
\Longrightarrow a~D_{\psi}(\mu_1,\mu_2)\leq \Vert \mu_1-\mu_2\Vert_{g_a}.
\end{array}
$$
In this situation,  Theorem~\ref{theo-lipa} yields the 
Kantorovich semi-distance exponential decays
\begin{equation}\label{wn2Dpsi}
a~D_{\psi}(\mu_1 (KL)^n,\mu_2 (KL)^n)\leq \varrho^{2n} ~\Vert \mu_1-\mu_2\Vert_{g_a}.
\end{equation}
When $(\XX,\psi)$ is a Polish space, the above estimates yield the exponential decays w.r.t.~the Wasserstein 1-distance $\Da_{\psi}$. By \cite[Theorem 2.3]{dmg-25},  under the assumptions of Theorem~\ref{theo-lipa}, there exists some $b>0$ such that
$$
\psi_{a,b}:=
\varphi_{g_a}+b~\psi\Longrightarrow
{\sf lip}_{\psi_{a,b},\psi_{a,b}}(KL)\leq 1-\left(1-{\sf lip}_{g_a,g_a}(KL)\right)^2.
$$
\begin{examp}\label{g-lyap-ill}
For Lyapunov functions on $\XX=\RR^d$ of the form
$$
g(x):=c_1~\exp{(c_2\Vert x\Vert^p)}
$$
 for constants $c_1,c_2>0$ and some parameter $p>0$ we have
 \begin{eqnarray*}
w_g(x_1,x_2)&\geq& 2c_1 e^{c_2\Vert x_1\Vert^p/2}
e^{c_2\Vert x_2\Vert^p/2} \geq 2c_1 \left(\exp{\left(\frac{c_2}{2}~(\Vert x_1\Vert^p+\Vert x_2\Vert^p)\right)}-1\right)\\
&\geq &\psi_{\sf exp}(x_1,x_2):=2c_1 \left(
\exp{ \left(\frac{c_2}{2^{1\vee p}} ~\Vert x_1-x_2\Vert^p\right)}-1\right).
\end{eqnarray*}
The first estimate comes from the fact that   $ u^2+v^2\geq 2uv$ for any $u,v\in\RR$.
The last estimate follows from the fact that
$$
\Vert x_1\Vert^p+\Vert x_2\Vert^p\geq 2^{-(p-1)_+}~\Vert x_1-x_2\Vert^p.
$$
Also note that for any $q\geq 1$ we have
$$
\psi_{\sf exp}(x_1,x_2)\geq \psi_{\sf p,q}(x_1,x_2):= c_{p,q}
\Vert x_1-x_2\Vert^{pq}\quad \mbox{\rm with}\quad
c_{p,q}:= \frac{2c_1}{q!}~\left(\frac{c_2}{2^{1\vee p}}\right)^q.
$$
In this case, in terms of the Wasserstein metric (\ref{def-c2}) we have
$$
a~c_{p,q}~\Da_{pq}(\mu_1,\mu_2)^{pq}\leq \Vert \mu_1-\mu_2\Vert_{g_a}.
$$
\end{examp}

  \subsection{Transportation cost inequality}\label{sect-transport-ineq}

  \subsubsection{Log-Sobolev constants}
We say that a probability measure $\mu_2$ satisfies the quadratic transportation cost inequality $\TT_2(\rho)$ with constant $\rho>0$ if for any probability $\mu_1$ we have
\begin{eqnarray}
\frac{1}{2}~\Da_{2}(\mu_1,\mu_2)^2&\leq& \rho~ \Ha(\mu_1~|~\mu_2).\label{K-talagrand}
\end{eqnarray}
We say that a Markov transition $K$ from $\XX$ into $\YY$ satisfies the quadratic transportation cost inequality $\TT_2(\rho)$ when $\delta_x K$ satisfies a quadratic transportation cost inequality $\TT_2(\rho)$ for any $x\in\XX$.

\begin{lem}[\cite{dp-25,del2025entropic}]\label{lem-Int-log-sob}
The Markov transition $K$ satisfies the $\TT_2(\rho)$ inequality  if and only if for any probability distribution $\mu$ 
and any Markov transition $L$  we  have 
\begin{equation}\label{lee-0}
\Da_{2}(\mu L,\mu K)^2\leq \int\mu(dx)~\Da_2\left(\delta_x L,\delta_xK\right)^2 \leq 2\rho ~\Ha(\mu\times L~|~\mu\times K).
 \end{equation} 
 In addition, we have the bias  estimate
 $$
\int~\mu(dx)~  \Vert\mbox{\rm mean}_{K}(x)-\mbox{\rm mean}_{L}(x)\Vert_{2}^2\leq 
 2\rho~\Ha\left(\mu\times L~|~\mu\times K\right)
 $$
 as well as the covariance entropic bounds
 $$
 \begin{array}{l}
\displaystyle \int~\mu(dx)~
 \Vert \mbox{\rm cov}_{K}(x)-\mbox{\rm cov}_{L}(x)\Vert_F\\
 \\
\displaystyle\leq 4\rho~\Ha\left(\mu\times L~|~\mu\times K\right)+c_{\mu,L}~\left(8\rho~\Ha\left(\mu\times L~|~\mu\times K\right)\right)^{1/2}\quad\mbox{with}\quad
c_{\mu,L}^2:=\mu\left(\tr(\mbox{\rm cov}_{L})\right)
 \end{array}
 $$

 \end{lem}

When $\XX=\RR^d$ and  $d\mu_1/d\mu_2$ is differentiable   
the (relative) Fisher information of $\mu_1$ with respect to $\mu_2$ is defined by
$$
\Ja(\mu_1~|~\mu_2)= \mu_1\left(\Vert\nabla \log{(d\mu_1/d\mu_2)}\Vert_2^2\right)
$$
We also use the convention $  \Ja(\mu_1~|~\mu_2)=\infty$  when $\mu_1\not\ll \mu_2$.

A probability measure $\mu_2$ satisfies the log-Sobolev inequality $LS(\rho)$  if for any probability $\mu_1$   we have
\begin{eqnarray}
 \Ha(\mu_1~|~\mu_2)
&\leq &\frac{\rho}{2}~
\Ja(\mu_1~|~\mu_2).\label{K-log-sob}
\end{eqnarray}
In addition,  a Markov transition $K(x,dy)$ from $\RR^d$ into itself  satisfies  the log-Sobolev inequality $LS(\rho)$ when $\delta_xK$ satisfies a   log-Sobolev inequality $LS(\rho)$ for any $x\in\RR^d$. Finally,  a Markov transition $K$ satisfies  the local log-Sobolev inequality $LS_{\sf loc}(\rho)$ if 
for any $x_1,x_2\in\RR^d$ we have the inequalities
\begin{equation}\label{reflocls}
\Ha(\delta_{x_1}K~|~\delta_{x_2}K)\leq \frac{\rho}{2}~\Ja(\delta_{x_1}K~|~\delta_{x_2}K).
\end{equation}
A theorem by Otto and Villani (\cite[Theorem 1]{otto-villani}) ensures that the  log-Sobolev inequality $LS(\rho)$
implies the quadratic transportation cost inequality $\TT_2(\rho)$.

We also say that $K$ satisfies the 
Fisher-Lipschitz inequality  $\JJ(\kappa)$ with constant $\kappa>0$  if 
for any $x_1,x_2\in\RR^d$ we have the inequalities
\begin{equation}\label{reflocllip}
\Ja(\delta_{x_1}K~|~\delta_{x_2}K)\leq 
  \kappa^2~\Vert x_1-x_2\Vert_2^2.
\end{equation}

  \subsubsection{Some illustrations}
Several criteria with explicit estimates of the log-Sobolev constants of Gibbs measures with non uniformly convex potentials are available in the literature under rather mild conditions, including the celebrated Bakry-Emery criteria on Riemannian manifolds~\cite{bakry-emery}. It is clearly out of the scope of the present article to review these estimates, see for instance~\cite{bgl,cattiaux-14,cattiaux,monmarche,otto,villani,villani-2} and references therein. 

Consider a Boltzmann-Gibbs probability measure on $\RR^d$ of the form
$$
\lambda_U(dx)=e^{-U(x)}~dx
$$
for some potential function $U(x)$ such that 
\begin{equation}\label{ex-lg-cg-v2}
\nabla^2 U(x)\geq 
a~1_{\Vert x\Vert_2 \geq  l}~I-b 
~1_{\Vert x\Vert_2< l}~I
\end{equation}
for some $a>0$, $b\geq 0$ and $l\geq 0$.
In this context there exists some $\overline{U}$ with bounded oscillations such that
$\nabla^2 (U+\overline{U})\geq (a/2)~I$. By the Holley-Stroock perturbation lemma~\cite{holley},
$\lambda_U$ satisfies the $LS(\rho)$-inequality (\ref{K-log-sob}) with some parameter $\rho=\rho(a,b,l)$ that depends on $(a,b,l)$.
 For an explicit  construction of such function we refer to~\cite{monmarche-AHL}.
As an example,  the one-dimensional double-well potential function $U(x)=x^4-a~x^2$, with $a>0$ are convex at infinity and thus satisfy 
the $LS(\rho)$-inequality for some $\rho>0$.

The local Hessian condition (\ref{ex-lg-cg-v2}) is often termed ``strongly convexity outside a ball'' or in short "convex-at-infinity". To avoid confusion,  we underline that this condition strongly differs from the strong convexity of $U$ 
on the non convex domain $\{x\in \RR^d:\Vert x\Vert_2\geq r\}$ discussed in~\cite{ma,peters,yan}.  It is also well known (see for instance Corollary 9.3.2 and Corollary 5.7.2 in~\cite{bgl}) that
\begin{equation}\label{ref-convex-UV}
\begin{array}{l}
\exists  \sigma\in  \Sa^+_d\quad\mbox{\rm such that}\quad \forall x\in\RR^d\quad
\nabla^2 U(x)\geq \sigma^{-1}\\
\\
\Longrightarrow\quad \mbox{\rm $\lambda_U$ satisfy the $LS(\rho)$ inequality with $\rho=\Vert\sigma\Vert_2$}.
\end{array}\end{equation}
We also refer to~\cite{monmarche} for a variety of log-Sobolev criteria when the curvature is not positive everywhere based on  Holley-Stroock and Aida-Shigekawa perturbation arguments~\cite{aida-shi,holley}.

We will now illustrate the transportation cost inequalities discussed above in the context of Gaussian models.

Denote by $\nu_{m,\sigma}$ the Gaussian distribution on $\XX=\RR^d$ with mean $m\in\RR^d$ and covariance matrix $\sigma \in \Sa^+_d$. In addition, let ${\sf gauss}_{\sigma}$ denote the probability density function of the distribution $\nu_{0,\sigma}$, with covariance matrix $\sigma \in \Sa^+_d$.  In this context, we have
$$
\nu_{m,\sigma}=\lambda_U\quad \mbox{\rm with}\quad
U(x)=-\log{{\sf gauss}_{\sigma}(x-m)}\Longrightarrow \nabla^2 U(x)=\sigma^{-1}. 
$$
For any parameters
$
(m,\overline{\sigma})\in(\RR^d\times \Sa^+_{d})$ 
we recall that 
 \begin{equation}\label{KL-def}
2~\Ha\left(\nu_{m,\sigma}~|~\nu_{\overline{m},\overline{\sigma}}\right)
 = 
  D(\sigma~|~\overline{\sigma})+\Vert \overline{\sigma}^{-1/2}\left(m-\overline{m}\right)\Vert^2_2
 \end{equation}
with the Burg (or log-det) divergence
\begin{equation}\label{burg-def}
 D_{\sf burg}(\sigma~|~\overline{\sigma}):=\tr\left(\sigma\overline{\sigma}^{-1}-I\right)-\log{\mbox{det}\left(\sigma\overline{\sigma}^{-1}\right)}.
\end{equation}
Let $(\alpha,\beta,\tau)\in  \left(\RR^{d}\times\Ga l_d\times \Sa^+_{d}\right)$,  we have the linear Gaussian transition
\begin{equation}\label{def-W}
\begin{array}{l}
K(x,dy)=e^{-W(x,y)}~dy
\quad \mbox{\rm with}\quad
W(x,y)=-\log{{\sf gauss}_{\tau}(y-(\alpha+\beta x))}\\
\\
\Longrightarrow \nabla_2W(x,y)=\tau^{-1}(y-(\alpha+\beta x))\quad \mbox{\rm and}\quad
 \nabla_2^2W(x,y)=\tau^{-1}\geq \Vert \tau\Vert_2^{-1}~I.
\end{array}
\end{equation}
The curvature estimate in the above display ensures that 
 $K$ satisfies  the  $LS(\rho)$-inequality with $\rho=\Vert \tau\Vert_2$.
Observe that
$$
2~\Ha\left(\delta_{x_1} K~|~\delta_{x_2} K\right)=
(x_1-x_2)^{\prime}\left((\tau^{-1/2}\beta)^{\prime}(\tau^{-1/2}\beta)\right)(x_1-x_2).
$$
On the other hand, we have
\begin{eqnarray*}
\Ja(\delta_{x_1}K~|~\delta_{x_2}K)
&=&\int~K(x_1,dy)~\Vert\nabla_2W(x_1,y)-\nabla_2W(x_2,y)\Vert^2\\
&=&(x_1-x_2)^{\prime}\left((\tau^{-1/2}\beta)^{\prime}\tau^{-1}(\tau^{-1/2}\beta)\right)(x_1-x_2).
\end{eqnarray*}
We recover the fact that  $K$ satisfies  the  $LS_{\sf loc}(\Vert \tau\Vert_2)$-inequality
and the $\JJ(\Vert\tau^{-1}\beta\Vert_2)$ inequality.

\subsubsection{A stability theorem}
For any probability measures $\mu_1$ and $\mu_2$ 
on  $\XX$ and a Markov transition $K(x,dy)$  from $\XX$ into $\YY$
we have
\begin{eqnarray*}
\Ha(\mu_1~|~\mu_2)&=&\Ha(\mu_1\times K~|~\mu_2\times K)=\Ha\left((\mu_1 K)\times K^{\ast}_{\mu_1}~|~(\mu_2 K)\times K^{\ast}_{\mu_2}\right)\nonumber\\
&=& \Ha(\mu_1 K~|~\mu_2 K)+\int (\mu_1 K)(dy)~\Ha\left(\delta_yK^{\ast}_{\mu_1}|\delta_yK^{\ast}_{\mu_2}\right)
\end{eqnarray*}
and hence that
\begin{equation}
\Ha(\mu_1 K~|~\mu_2 K)\leq \Ha(\mu_1~|~\mu_2).\label{contK}
\end{equation}
In information theory, the above inequality is also known as the data-processing inequality.
For any probability measures $\mu$ and $\eta$ 
on  $\XX$ and $\YY$,
from  the convexity of the function $h(u)=u\log u$ we have 
\begin{eqnarray*}
  h\left(\frac{d\mu K }{d\eta}(y)\right)
 &=& h\left(\int~\mu(dx)~\frac{d\delta_{x} K }{d\eta}(y)\right)\leq \int \mu(dx)~ h\left(\frac{d\delta_x K }{d\eta}(y)\right).
\end{eqnarray*}
This implies that
\begin{eqnarray*}
\Ha(\mu K ~|~\eta)&=&\int~\eta(dy) ~ h\left(\frac{d\mu K }{d\eta}(y)\right)\\
&\leq &\int~~\mu(dx) \int~\eta(dy) ~ h\left(\frac{d\delta_x K }{d\eta}(y)\right) =
\int~\mu(dx)~\Ha(\delta_{x} K ~|~\eta).
\end{eqnarray*}
We conclude that
\begin{equation}\label{Hconvx}
\Ha(\mu K ~|~\eta)\leq \int~\mu(dx)~\Ha(\delta_{x} K ~|~\eta).
\end{equation}

Let $\mu$ and $\eta$ be two probability measures 
on  $\XX$  and let $K(x,dy)$ be Markov transition  from $\XX$ into $\YY$.
From the bridge map definition  given in (\ref{def-entropy-pb-v2}) that for any reference measure of the form $P:=(\mu\times K)$ for some given probability measures $\mu$ on $\XX$ we have
$$
P_{\mu,\mu K}=\argmin_{Q\,\in\, \Pi(\mu,\mu K)}\Ha(Q|P)=\mu\times K\quad\mbox{\rm and}\quad
 P_{\mu,\eta K}:=\argmin_{Q\,\in\, \Pi(\mu,\eta K)}\Ha(Q|P).
$$
We are now in position to state and prove the following stability theorem.

\begin{theo}[\cite{dp-25}]\label{key-theo}
Assume that $\XX=\YY=\RR^d$ and the Markov transition
 $K$ satisfies  the local log-Sobolev inequality $LS_{\sf loc}(\rho)$ and the Fisher-Lipschitz inequality $\JJ(\kappa)$ inequality. In this situation, for any reference measure of the form $P:=(\mu\times K)$ for some probability measures $\mu$
on  $\RR^d$, 
and for any probability measure $\eta$ on $\RR^d$ have
\begin{equation}\label{H2W}
\Ha(\eta K|\mu K)\leq \Ha\left(P_{\mu,\eta K}~|~P_{\mu,\mu K}\right)
\leq  \kappa^2\rho ~ \frac{1}{2}~\Da_2(\mu,\eta)^2.
\end{equation}
In addition, we have the following estimates:\\
\noindent $\bullet$ When $\eta$ satisfies the $\TT_2(\rho^{\prime})$ inequality we have
\begin{equation}\label{H2WetaT2}
 \Ha\left(P_{\mu,\eta K}~|~P_{\mu,\mu K}\right)
\leq  \kappa^2\rho\rho^{\prime} ~\Ha(\mu|\eta).
\end{equation}
$\bullet$ When $\mu$ satisfies the $\TT_2(\rho^{\prime})$ inequality we have
\begin{equation}\label{H2WmuT2}
 \Ha\left(P_{\mu,\eta K}~|~P_{\mu,\mu K}\right)
\leq  \kappa^2\rho\rho^{\prime} ~\Ha(\eta|\mu).
\end{equation}
$\bullet$ When $K$ satisfies the $\TT_2(\rho^{\prime})$ inequality we have
\begin{equation}\label{H2WKT2}
\Da_2(\eta K,\mu K)^2\leq 2\rho^{\prime}~\Ha\left(P_{\mu,\eta K}~|~P_{\mu,\mu K}\right)
\leq  \kappa^2\rho\rho^{\prime} ~\Da_2(\mu,\eta)^2.
\end{equation}
\end{theo}

\proof
Applying (\ref{Hconvx}) to a Markov transition $T$ we check that
\begin{eqnarray*}
\Ha\left(\mu\times (TK)|\mu \times K\right)&=&
\int \mu(dz)
\Ha((\delta_{z} T)K ~|~\delta_{z}K)\\
& \leq&
 \int~(\mu\times T)^{\flat}(d(x,z))~\Ha(\delta_{x} K ~|~\delta_{z}K).
\end{eqnarray*}
Choosing $T$ such that $\mu T=\eta$ we check that
$$
(\mu\times T)^{\flat}\in \Pi\left(\eta,\mu\right)\quad \mbox{\rm and}\quad
(\mu\times (TK))\in \Pi\left(\mu ,\eta K\right).
$$
On the other hand, by (\ref{reflocls}) we have
$$
\Ha(\delta_{x_1}K~|~\delta_{x_2}K)\leq \frac{\rho}{2}~  \kappa^2~\Vert x_1-x_2\Vert_2^2.
$$
For any $T$ such that $\mu T=\eta$ we have
$$
P:=\mu\times K\Longrightarrow \Ha\left(P_{\mu,\eta K}~|~P_{\mu,\mu K}\right)
\leq \frac{\rho}{2}~  \kappa^2 \int~\mu(dz)~T(z,dx)~~\Vert z-x\Vert_2^2.
$$
We end the proof of (\ref{H2W}) choosing the optimal $2$-Wasserstein coupling transition.
The estimates (\ref{H2WetaT2}) and (\ref{H2WmuT2}) are immediate consequences of the definition of the
 $\TT_2(\rho^{\prime})$ inequality.
The estimate (\ref{H2WKT2}) is a consequence of Lemma~\ref{lem-Int-log-sob} and
(\ref{H2W}).
\cqfd

\section{A brief review on Sinkhorn bridges}\label{SS-bridges}

\subsection{Overview}

Entropic optimal transport for density models on Polish spaces  is often presented as 
an optimal transport problem with an entropic regularization. The entropic regularization problem is solved by Sinkhorn iterations  for small values of the regularization parameter to solve the initial optimal transport problem. 
A brief presentation of these models is provided in Section~\ref{ert-sec}. 
Recall that the Sinkhorn bridges introduced in (\ref{def-Pa-n}) are defined sequentially by conjugate formulae and adjoint integral operators (cf. Section~\ref{conjugate-formulation-sec}). This probabilistic interpretation allows one to connect the convergence Sinkhorn semigroups with the stability analysis of Markov semigroups. 
In Section~\ref{sec-monotone-Markov} we briefly review some elementary monotone properties of Sinkhorn semigroup based on the non expansive property of Markov operators (\ref{contK}) w.r.t. the relative entropy criteria.
In the context of density models, Schr\"odinger and Sinkhorn bridges
can alternatively be defined by a sequence of entropic-potential functions.
These entropic potentials are defined in terms of a system of recursive equations. 
These dual-type formulations are presented in Section~\ref{sch-pot}. 
Section~\ref{enproj-sec} and Section~\ref{commut-sec} provide some entropic projection properties of bridge maps and present
  some pivotal commutation-type integral transport formulae. Section~\ref{gred-Hess-sec} underlines some links between the gradient and Hessian of entropic-potential functions and the conditional mean and covariance of Sinkhorn transitions.

\subsection{Entropy-regularized transport}\label{ert-sec}
Consider the marginal measures $(\mu,\eta)=(\lambda_U,\nu_{V})$ 
on locally compact Polish spaces $(\XX,\lambda)$ and $(\YY,\nu)$ defined in
(\ref{ref-intro-UV}). 
 Denote $K$ a Markov transition from $\XX$ into $\YY$ of the form 
 \begin{equation}\label{def-Qa}
K(x,dy)=e^{-W(x,y)}~\nu(dy)
\end{equation}
for some transition potential function $W:(x,y)\in(\XX\times\YY)\mapsto W(x,y)\in\RR$.  In this context, given the reference measure $ P:=(\mu\times K)$, for any $Q\,\in\, \Pi(\mu,\eta)$ we check that
 $$
 \begin{array}{l}
 \displaystyle\frac{d Q}{d P}=\frac{d( \mu\otimes \eta)}{dP}\times 
  \frac{d Q}{d( \mu\otimes \eta)}
 \displaystyle\Longrightarrow  \Ha(Q~|~P)= Q(W)+
\Ha(Q~|~\mu\otimes \eta)+\Ha(\eta~|~\nu).
  \end{array}
 $$
 Therefore the problem (\ref{def-entropy-pb-v2}) is equivalent to the optimal entropic transport problem
 \begin{equation}\label{reg-version-0}
P_{\mu,\eta}:=\argmin_{Q\,\in\, \Pi(\mu,\eta)}\Ha(Q|P)=\argmin_{Q\,\in\, \Pi(\mu,\eta)}\left[ Q(W)+
 \Ha\left(Q~|~\mu\otimes \eta\right)\right].
\end{equation}
Choosing a potential transition of the form $W(x,y)=\Wa(x,y)/t+c_t$,  for some parameter $t>0$ and some constant $c_t$, the problem (\ref{def-entropy-pb-v2}) is now equivalent to the optimal entropic transport problem
\begin{equation}\label{reg-version}
P_{\mu,\eta}:=\argmin_{Q\,\in\, \Pi(\mu,\eta)}\left[ Q(\Wa)+t~
 \Ha\left(Q~|~\mu\otimes \eta\right)\right].
\end{equation}
  In the optimal transport literature the function $\Wa$ is sometimes called the cost function and the parameter $t>0$ is interpreted as a regularization parameter.  In this context, the ``unregularized" problem corresponding to $t=0$ coincides with the  Kantorovich optimal transport problem
  with cost function $\Wa$. The quadratic cost model on $\XX=\RR^d=\YY$ with the Lebesgue measure $\nu(dy)=dy$ is defined by
\begin{equation}\label{gibbs-cost}
W(x,y)=\frac{1}{t}~\Wa(x,y)+\frac{d}{2}\log{(2\pi t)}\quad \mbox{\rm with}\quad
\Wa(x,y)=\frac{1}{2}~\Vert x-y\Vert^2
\end{equation}
In this situation, $K(x,dy)$ corresponds to the heat equation transition semigroup on a given time horizon $t>0$.
  
We remark that $  (\mu\otimes \eta)\in \Pi(\mu,\eta)$ and
  $$
\Ha(\mu\otimes \eta|P)<\infty\Longleftrightarrow
(\mu\otimes \eta)(W)<\infty.
  $$
In contrast with Kantorovich criteria or related Wasserstein-type semi-distances,
note that even when $\XX=\YY$ and the same marginal $\mu=\eta$
the identity-map deterministic coupling has infinite entropy; that is we have that
$$
Q(d(x,y)):=\mu(dx)~\delta_x(dy)\in  \Pi(\mu,\eta)\Longrightarrow
 \Ha(Q~|~P)=+\infty.
$$
This establishes that the Schr\"odinger bridge  $P_{\mu,\eta}$  strongly depends on the reference probability distribution $P=(\mu\times K)$, equivalently on the cost function $W$ for density models with reference transition $K$ given by (\ref{def-Qa}).

\subsection{Some monotone properties}\label{sec-monotone-Markov}
Applying (\ref{pre-gibbs}) we have the Markov transport equations 
\begin{equation}\label{decrease2n}
\left(\pi_{2(n+1)},\eta\right)=\left(\mu \Ka_{2(n+1)},\pi_{2n+1} \Ka_{2(n+1)}\right)
\quad\mbox{\rm and}\quad \left(\mu,\pi_{2n+1}\right)=
\left(\pi_{2n} \Ka_{2n+1},\eta \Ka_{2n+1}\right).
\end{equation}
Combining (\ref{decrease2n}) with (\ref{contK}) we check that
$$
\begin{array}{rcccl}
 \Ha\left(\pi_{2(n+1)}|\eta\right)&\leq &\Ha(\mu|\pi_{2n+1})&\leq &
 \Ha\left(\pi_{2n}|\eta\right)\\
 &&&&\\
 \Ha\left(\eta|\pi_{2(n+1)}\right)&\leq &\Ha(\pi_{2n+1}|\mu)&\leq &
 \Ha\left(\eta|\pi_{2n}\right).
\end{array} $$
In the same vein, we have the Markov transport equations 
\begin{equation}\label{decrease2n1}
\left(\pi_{2n+1},\mu\right)=\left(\eta \Ka_{2n+1},\pi_{2n} \Ka_{2n+1}\right)
\quad\mbox{\rm and}\quad \left(\eta,\pi_{2n}\right)=
\left(\pi_{2n-1} \Ka_{2n},\mu \Ka_{2n}\right).
\end{equation}
Combining (\ref{decrease2n1}) with (\ref{contK}) we check that
$$
\begin{array}{rcccl}
\Ha\left(\pi_{2n+1}|\mu\right)&\leq& \Ha(\eta|\pi_{2n})&\leq &
 \Ha\left(\pi_{2n-1}|\mu\right)
 \\
 &&&&\\
  \Ha\left(\mu|\pi_{2n+1}\right)&\leq& \Ha(\pi_{2n}|\eta)&\leq &
 \Ha\left(\mu|\pi_{2n-1}\right).
\end{array} $$
We summarize the above discussion with the following proposition.
\begin{prop}\label{prop-decrease2n01}
For any $n\geq 0$ we have
$$
\begin{array}{rcccl}
 \Ha\left(\eta|\pi_{2(n+1)}\right)&\leq& \Ha\left(\pi_{2n+1}|\mu\right)&\leq &\Ha(\eta|\pi_{2n})
\\
&&&&\\
\Ha\left(\pi_{2(n+1)}|\eta\right)&\leq&  \Ha\left(\mu|\pi_{2n+1}\right)&\leq& \Ha(\pi_{2n}|\eta).
\end{array} $$
\end{prop}

\subsection{Sinkhorn potentials}\label{sch-pot}
Consider the marginal measures $(\mu,\eta)=(\lambda_U,\nu_V)$ and the Markov transition $K$ defined in
(\ref{ref-intro-UV}) and (\ref{def-Qa}). Note that
$$
 K(x,dy)=k(x,y)~\nu(dy)
\quad  \mbox{\rm with the density}\quad 
k(x,y)=e^{-W(x,y)}.
 $$
We introduce the backward-type (non necessarily Markov) integral operator
 \begin{equation}\label{def-Ra}
K^{\flat}(y,dx):=k^{\flat}(y,x)~\lambda(dx)
\end{equation}
with the density
$$
k^{\flat}(y,x)=e^{-W^{\flat}(y,x)}=k(x,y)
\Longleftrightarrow W^{\flat}(y,x)=W(x,y).
$$
Recalling that $\gamma:=\lambda\otimes\nu$, we have
$$
\Pa_0(d(x,y))=p_0(x,y)~\gamma(d(x,y))
$$
with the density
$$
p_0(x,y):=e^{-U_0(x)}~k(x,y)~e^{-V_0(y)}\quad \mbox{\rm and the potential functions}\quad
(U_0,V_0)=(U,0).
$$
Note that
$$
\Ka_0(x,dy)=\frac{K(x,dy) e^{-V_0(y)}}{K(e^{-V_0})(x)}=K(x,dy)
$$
Using (\ref{rec-bayes-marg}) and (\ref{rec-bayes}) we check that
$$
\Pa_1(d(x,y))=p_1(x,y)~\gamma(d(x,y))
$$
with
$$
p_1(x,y)=e^{-V(y)}~\frac{p_0(x,y)}{\int~p_0(\overline{x},y)~\lambda(d\overline{x})}=
e^{-U_0(x)}~k(x,y)~\frac{e^{-V(y)}}{K^{\flat}(e^{-U_0})(y)}.
$$
This implies that
$$
p_1(x,y):=e^{-U_1(x)}~k(x,y)~e^{-V_1(y)}
$$
with
$$
U_1=U_0\quad \mbox{\rm and}\quad V_1=V+\log{K^{\flat}(e^{-U_0})}
\Longrightarrow
\Ka_1(y,dx)=\frac{K^{\flat}(y,dx)~e^{-U_1(x)}}{K^{\flat}(e^{-U_1})(y)}.
$$
In the same vein, we have
$$
\Pa_2(d(x,y))=p_2(x,y)~\gamma(d(x,y))
$$
with
$$
p_{2}(x,y)=e^{-U(x)}~\frac{p_1(x,y)}{\int~p_1(x,\overline{y})~\nu(d\overline{y})}=
\frac{e^{-U(x)}}{K(e^{-V_1})(x)}~k(x,y)~e^{-V_1(y)}.
$$
This implies that
$$
p_2(x,y):=e^{-U_2(x)}~k(x,y)~e^{-V_2(y)}
$$
with
$$
U_2=U+\log{K(e^{-V_1})}\quad \mbox{\rm and}\quad V_2=V_1
\Longrightarrow
\Ka_2(x,dy)=\frac{K(x,dy) e^{-V_2(y)}}{K(e^{-V_2})(x)}.
$$
Iterating the above procedure we check the following proposition.
\begin{prop}\label{PUVn-prop}
Sinkhorn bridges are given for any $n\geq 0$ by
 \begin{equation}\label{PUVn}
\Pa_n(d(x,y))=p_n(x,y)~\gamma(d(x,y))
\quad\mbox{with}\quad
p_n(x,y):=e^{-U_n(x)}~k(x,y)~e^{-V_n(y)}.
\end{equation}
The potentials $(U_n,V_n)$ are defined for any $n\geq 0$ by the recursions
 \begin{equation}\label{UVn}
\begin{array}{rclcrcl}
U_{2n+1}&=&U_{2n}& & V_{2n+1}&=&V+\log{K^{\flat}(e^{-U_{2n}})}\\
U_{2(n+1)}&=&U+\log{K(e^{-V_{2n+1}})}&& V_{2(n+1)}&=&V_{2n+1}
\end{array}
\end{equation}
with the initial condition $(U_0,V_0)=(U,0)$. 
\end{prop}
By Proposition~\ref{PUVn-prop} Sinkhorn transitions take the following form
 \begin{equation}\label{Kan}
 \begin{array}{rcl}
\displaystyle\Ka_{2n}(x,dy)&=&\displaystyle K_{V_{2n}}(x,dy):=\frac{K(x,dy) ~e^{-V_{2n}(y)}}{K(e^{-V_{2n}})(x)}\\
&&\\
\Ka_{2n+1}(y,dx)&=&\displaystyle K^{\flat}_{U_{2n+1}}(y,dx):=\frac{K^{\flat}(y,dx)~e^{-U_{2n+1}(x)}}{K^{\flat}(e^{-U_{2n+1}})(y)}.
\end{array}
\end{equation}
Rewriting (\ref{Kan}) in a slightly different form we check the following proposition.
\begin{prop}
For any $n\geq 0$ we have
$$
\Ka_{2n}(x,dy)=e^{-W_{2n}(x,y)}~\nu(dy)\quad \mbox{\rm and}\quad
\Ka_{2n+1}(y,dx)=e^{-W_{2n+1}^{\flat}(y,x)}~\lambda(dx)
$$
with the transition potential functions
\begin{eqnarray}
U(x)+W_{2n}(x,y)&=&U_{2n}(x)+W(x,y)+V_{2n}(y)\nonumber\\
V(y)+W_{2n+1}^{\flat}(y,x)&=&V_{2n+1}(y)+W^{\flat}(y,x)+U_{2n+1}(x)\label{UV-Wn}
\end{eqnarray}
with the potential functions $(U_n,V_n)$ as in (\ref{UVn}).
\end{prop}
Using (\ref{PUVn}) and recalling that $(U_0,V_0)=(U,0)$ we check that
\begin{equation}\label{ratio-UV}
V_{2(l+1)}=V_{2l}-\log{\left(\frac{d\eta}{d\pi_{2l}}\right)}\quad \mbox{\rm and}\quad
U_{2(l+1)}=U_{2l}-\log{\left(\frac{d\mu}{d\pi_{2l+1}}\right)}
\end{equation}
This yields the following series formulation.
\begin{prop}\label{UVn-prop-series}
Sinkhorn potentials (\ref{UVn}) are given for any $n\geq 0$ by the formulae
$$
V_{2n}=V_{0}-\sum_{0\leq l<n }\log{\left(\frac{d\eta}{d\pi_{2l}}\right)}\quad \mbox{and}\quad
U_{2n}=U_{0}-\sum_{0\leq l<n }\log{\left(\frac{d\mu}{d\pi_{2l+1}}\right)}.
$$
In addition, we have
$$
\eta(V_{2(n+1)})\leq \eta(V_{2n})\leq \eta(V_0)=0
\quad \mbox{and}\quad
\mu(U_{2(n+1)})\leq \mu(U_{2n})\leq \mu(U).
$$
\end{prop}

Condition $\Ha(P_{\mu,\eta}|P)<\infty$ ensures the convergence of the series
(\ref{kdec2n-intro}) and Sinkhorn entropic potentials functions $(U_n,V_n)\rightarrow_{n\rightarrow\infty}(\UU,\VV)$ with the solution of the fixed point equation
associated with the recursion (\ref{UVn}) defined by the Schr\" odinger system
\begin{equation}\label{def-UU-VV-e}
\begin{array}{l}
\UU=U+\log{K(e^{-\VV})}
\quad\mbox{\rm and}\quad
\VV=V+\log{K^{\flat}(e^{-\UU})}\\
\\
\displaystyle \Longleftrightarrow
\frac{e^{-U}}{K(e^{-\VV})}=e^{-\UU}
\quad\mbox{\rm and}\quad \frac{e^{-V}}{K^{\flat}(e^{-\UU})}=e^{-\VV}
\end{array}
\end{equation}

In this context, the bridge takes the form
$$
P_{\mu,\eta}(d(x,y))=\lambda_U(dx)~K_{\VV}(x,dy)=\nu_V(dy)~K^{\flat}_{\UU}(y,dx)
$$
for some entropic potential functions $(\UU,\VV)$ and the transitions
\begin{equation}\label{def-KVV-KUU}
K_{\VV}(x,dy):=\frac{K(x,dy) ~e^{-\VV(y)}}{K(e^{-\VV})(x)}\quad\mbox{\rm and}\quad
K^{\flat}_{\UU}(y,dx):=\frac{K^{\flat}(y,dx)~e^{-\UU(x)}}{K^{\flat}(e^{-\UU})(y)}.
\end{equation}
An alternative description of Schr\"odinger bridges is provided in Theorem 4.2 in~\cite{nutz}, see also Section~\ref{sec-series-UU} and Theorem~\ref{ac-series-theo}  in the present article for a description of the
entropic potential functions $(\UU,\VV)$ in terms of absolutely convergent series. An explicit description of the entropic potential functions in the context of linear-Gaussian models is also provided in the article~\cite{adm-24}.

\subsection{Entropic projections}\label{enproj-sec}
Using (\ref{bayes2n}) and (\ref{bayes2n1})  for any $n>p$ we have the product formulae
\begin{eqnarray}
\frac{d\Pa_{2n}}{d\Pa_{2p}}(x,y)&=&
\left[\prod_{p\leq l< n}\frac{d\Pa_{2l+1}}{d\Pa_{2l}}(x,y)\right]~\left[\prod_{p\leq l< n}\frac{d\Pa_{2(l+1)}}{d\Pa_{2l+1}}(x,y)\right]\nonumber\\
&=&\left[\prod_{p\leq l< n}\frac{d\eta}{d\pi_{2l}}(y)\right]~\left[\prod_{p\leq l< n}\frac{d\mu}{d\pi_{2l+1}}(x)\right].\label{kp-i}
\end{eqnarray}
This yields for any $Q\in \Pi(\mu,\eta)$ and $n>p\geq 0$ the formulae
\begin{equation}\label{kdec}
Q\left(\log{\frac{d\Pa_{2n}}{d\Pa_{2p}}}\right)=\Ha(Q|\Pa_{2p})-
\Ha(Q|\Pa_{2n})=\sum_{p\leq l< n}\left[\Ha(\eta|\pi_{2l})+\Ha(\mu|\pi_{2l+1})\right].
\end{equation}
Applying the above to $p=0$  and
recalling that $\Pa_0=P$, this yields  the formula
\begin{eqnarray}
\Ha(Q|P)
&=&\Ha(Q|\Pa_{2n})+
\sum_{0\leq l< n}\left[\Ha(\eta|\pi_{2l})+\Ha(\mu|\pi_{2l+1})\right].\label{keysum-p}
\end{eqnarray}
Observe that
\begin{eqnarray}
\Ha(Q|\Pa_{2n})&=&\Ha(Q|\Pa_{2n+1})+Q\left(\log{\frac{d\Pa_{2n+1}}{\Pa_{2n}}}\right)
=\Ha(Q|\Pa_{2n+1})+\Ha(\eta|\pi_{2n})
\nonumber\\
\Ha(Q|\Pa_{2n+1})&=&\Ha(Q|\Pa_{2(n+1)})+Q\left(\log{\frac{d\Pa_{2(n+1)}}{\Pa_{2n+1}}}\right)
=\Ha(Q|\Pa_{2(n+1)})+\Ha(\mu|\pi_{2n+1}).~~~~~~~~\label{revHn}
\end{eqnarray}
This implies that
\begin{eqnarray}
\Ha(Q|P)
&=&\Ha(Q|\Pa_{2n+1})+\Ha(\eta|\pi_{2n})+
\sum_{0\leq l< n}\left[\Ha(\eta|\pi_{2l})+\Ha(\mu|\pi_{2l+1})\right].\label{keysum}
\end{eqnarray}

Combining (\ref{keysum-p}) with (\ref{keysum}) and (\ref{def-entropy-pb-v2})
we find that the Schr\" odinger bridge $P_{\mu,\eta}$ 
with reference measure $P$ defined in (\ref{def-entropy-pb-v2}) can be computed from any   Sinkhorn reference bridge $\Pa_{n}$. Equivalently, the bridge map $P\mapsto P_{\mu,\eta}$
is constant along the flow of Sinkhorn bridges $\Pa_{n}$. More formally we have the following proposition.

\begin{prop}
 For any $n\geq 0$ we have 
\begin{equation}\label{bridge-Pn}
(\Pa_{n})_{\mu,\eta}=\argmin_{Q\,\in\, \Pi(\mu,\eta)}\Ha(Q|\Pa_{n})=P_{\mu,\eta}.
\end{equation}
\end{prop} 

Applying (\ref{kdec}) to $Q=P_{\mu,\eta}$ we have
\begin{equation}\label{kdec2}
\Ha(P_{\mu,\eta}|\Pa_{2p})=
\Ha(P_{\mu,\eta}|\Pa_{2n})+\sum_{p\leq l< n}\left[\Ha(\eta|\pi_{2l})+\Ha(\mu|\pi_{2l+1})\right].
\end{equation}
We further assume  the series on the r.h.s.~of (\ref{kdec2}) converges. 
By Proposition~\ref{prop-decrease2n01}, this condition is met if and only if the sequence of entropies $\Ha(\eta|\pi_{2l})$ or 
$\Ha(\mu|\pi_{2l+1})$ is summable.

Whenever
\begin{equation}\label{kdec2b}
\lim_{n\rightarrow\infty}\Ha(P_{\mu,\eta}|\Pa_{2n})=0,
\end{equation}
 letting $n\rightarrow\infty$ in (\ref{kdec2}), for any $p\geq 0$ we have
\begin{equation}\label{kdec2c}
\Ha(P_{\mu,\eta}|\Pa_{2p})=\sum_{l\geq p}\left[\Ha(\eta|\pi_{2l})+\Ha(\mu|\pi_{2l+1})\right].
\end{equation}
Conversely when (\ref{kdec2c}) is satisfied (\ref{kdec2b}) clearly holds. We summarize the above discussion with the following proposition.
\begin{prop}\label{kdec2n-prop}
 We have
the equivalence (\ref{kdec2n-intro}) as soon as  the sequence of entropies $\Ha(\eta|\pi_{2l})$ or 
$\Ha(\mu|\pi_{2l+1})$ is summable.
\end{prop}

We end this section with some elementary monotone properties.
 Applying (\ref{revHn})  to  the bridge $Q=P_{\mu,\eta}$ and using Proposition~\ref{prop-decrease2n01} we check that
\begin{eqnarray}
\Ha(P_{\mu,\eta}|\Pa_{2(n+1)})&=&\Ha(P_{\mu,\eta}|\Pa_{2n+1})-\Ha(\mu|\pi_{2n+1})\label{pi2nP}\\
&=&\Ha(P_{\mu,\eta}|\Pa_{2n})-\left[\Ha(\eta|\pi_{2n})+\Ha(\mu|\pi_{2n+1})\right]\nonumber\\
&\leq &\Ha(P_{\mu,\eta}|\Pa_{2n})-\left[\Ha(\pi_{2n+1}|\mu)+\Ha(\pi_{2(n+1)}|\eta)\right].\nonumber
\end{eqnarray}
Similarly, we have
\begin{eqnarray}
\Ha(P_{\mu,\eta}|\Pa_{2n+1})&=& \Ha(P_{\mu,\eta}~|~\Pa_{2n})-   \Ha(\eta~|~\pi_{2n})\label{pi2n1P}\\
&=&
\Ha(P_{\mu,\eta}|\Pa_{2n-1})-\left[\Ha(\mu|\pi_{2n-1})
+\Ha(\eta|\pi_{2n})\right]\nonumber\\
&\leq& \Ha(P_{\mu,\eta}|\Pa_{2n-1})-\left[\Ha(\pi_{2n}|\eta)
+\Ha(\pi_{2n+1}|\mu)\right].\nonumber
\end{eqnarray}

\subsection{Commutation formulae}\label{commut-sec}
Observe that
$$
\begin{array}{l}
\displaystyle\frac{d\Pa_{2(l+1)}}{d\Pa_{2l}}(x,y)=
\frac{d\Pa_{2(l+1)}}{d\Pa_{2l+1}}(x,y)~\frac{d\Pa_{2l+1}}{d\Pa_{2l}}(x,y)
=\frac{d\mu}{d\pi_{2l+1}}(x)~\frac{d\eta}{d\pi_{2l}}(y)\\
\\
\displaystyle\Longleftrightarrow
\Ka_{2(l+1)}(x,dy)=\frac{d\mu}{d\pi_{2l+1}}(x)~\Ka_{2l}(x,dy)~\frac{d\eta}{d\pi_{2l}}(y).
\displaystyle
\end{array}
$$
Similarly,  we have
$$
\begin{array}{l}
\displaystyle\frac{d\Pa_{2l+1}}{d\Pa_{2l-1}}(x,y)=
\frac{d\Pa_{2l+1}}{d\Pa_{2l}}(x,y)~\frac{d\Pa_{2l}}{d\Pa_{2l-1}}(x,y)
=\frac{d\eta}{d\pi_{2l}}(y)~\frac{d\mu}{d\pi_{2l-1}}(x)~\\
\\
\displaystyle\Longleftrightarrow
\Ka_{2l+1}(y,dx)=\frac{d\eta}{d\pi_{2l}}(y)~\Ka_{2l-1}(y,dx)~
\frac{d\mu}{d\pi_{2l-1}}(x).
\displaystyle
\end{array}
$$
These computations allow us to deduce the following result.
\begin{prop}\label{prop-commute}
For any $l\geq 0$ we have the commutation formulae
$$
\Ka_{2l}\left(\frac{d\eta}{d\pi_{2l}}\right)=\frac{d\pi_{2l+1}}{d\mu}\quad \mbox{and}\quad
\Ka_{2(l+1)}\left(\frac{d\pi_{2l}}{d\eta}\right)=\frac{d\mu}{d\pi_{2l+1}}.
$$
In addition, for any $l\geq 1$ we have
$$
\Ka_{2l-1}\left(
\frac{d\mu}{d\pi_{2l-1}}\right)=\frac{d\pi_{2l}}{d\eta}
\quad \mbox{and}\quad
\Ka_{2l+1}\left(\frac{d\pi_{2l-1}}{d\mu}\right)=\frac{d\eta}{d\pi_{2l}}.
$$
\end{prop}

In terms of the forward-backward Gibbs-type transitions
(\ref{fb-sinkhorn}) we have
$$
\frac{d\pi_{2l}}{d\eta}=\Ka_{2l-1}\left(\frac{d\mu}{d\pi_{2l-1}}\right)=
\Ka_{2l-1}\Ka_{2l}\left(\frac{d\pi_{2(l-1)}}{d\eta}\right)=\Sa_{2l}\left(\frac{d\pi_{2(l-1)}}{d\eta}\right).
$$
In addition, we have
$$
\frac{d\pi_{2l+1}}{d\mu}=\Ka_{2l}\left(\frac{d\eta}{d\pi_{2l}}\right)=\Ka_{2l}\left(\Ka_{2l+1}\left(\frac{d\pi_{2l-1}}{d\mu}\right)\right)=
\Sa_{2l+1}\left(\frac{d\pi_{2l-1}}{d\mu}\right).
$$
This yields the following result.
\begin{prop}\label{prop-commute-S}
For any $l\geq 1$ we have
$$
\Sa_{2l}\left(\frac{d\pi_{2(l-1)}}{d\eta}\right)=\frac{d\pi_{2l}}{d\eta}
\quad \mbox{and}\quad
\Sa_{2l+1}\left(\frac{d\pi_{2l-1}}{d\mu}\right)=\frac{d\pi_{2l+1}}{d\mu}.
$$
\end{prop}

\subsection{Gradient and Hessian formulae}\label{gred-Hess-sec}

Consider the marginal measures $(\mu,\eta)$ and the transition $K$  discussed in Section~\ref{sch-pot} for some potential functions $(U,V)$ and $W$ on the state spaces  $\XX=\RR^d=\YY$ equipped with the Lebesgue measure $\lambda(dz)=dz=\nu(dz)$.  In addition,  assume $(U,V)$ in (\ref{ref-intro-UV}) are twice differentiable and $W$ is the linear-Gaussian transition potential defined in  (\ref{def-W}), for some parameters $(\alpha,\beta,\tau)\in  \left(\RR^{d}\times\Ga l_d\times \Sa^+_{d}\right)$.  In this scenario the gradient and the Hessian defined in (\ref{ref-grad-hessian-W}) reduce to
\begin{equation}\label{ref-hessian-W-f}
\begin{array}{l}
 W_x(y)=\cchi^{\prime}((\alpha+\beta x)-y)\quad\mbox{and}\quad
W^{\flat}_y(x)=\tau^{-1}(y-(\alpha+\beta x)) \\
\\
\Longrightarrow
W_{x,x}(y)= \cchi^{\prime}\beta
\quad\mbox{and}\quad  W^{\flat}_{y,y}(x)=\tau^{-1}
\quad\mbox{with}\quad\cchi:=\tau^{-1}\beta.
\end{array}
\end{equation}
Therefore, the conditional covariances (\ref{ref-cond-cov}) and (\ref{cov-mat-W}) associated with Sinkhorn transitions (\ref{Kan}) take the form
$$
  \mbox{\rm cov}_{\Ka_{2n},W}(x)=\cchi^{\prime}~
  \mbox{\rm cov}_{\Ka_{2n}}(x)~  \cchi
  \quad \mbox{\rm and}\quad
  \mbox{\rm cov}_{\Ka_{2n+1},W^{\flat}}(y)  =
  \cchi~  \mbox{\rm cov}_{\Ka_{2n+1}}(y)~\cchi^{\prime}.
$$
Combining (\ref{ref-grad-hessian}) with (\ref{UVn}) and (\ref{Kan}) one has the following result.
\begin{prop}\label{prop-gh-uv}
For any $n\geq 0$ we have the gradient formulae
 \begin{eqnarray}
\nabla U_{2n}(x)&=&\nabla U(x)-\cchi^{\prime}\left\{(\alpha+\beta x)-\mbox{\rm mean}_{\Ka_{2n}}(x)\right\}\nonumber\\
\nabla V_{2n+1}(y)&=&\nabla V(y)-\tau^{-1}\left\{y-(\alpha+\beta~ \mbox{\rm mean}_{\Ka_{2n+1}}(y))\right\}\label{ref-nablaUVng}
 \end{eqnarray}
with the conditional expectations $\mbox{\rm mean}_{\Ka_{n}}$ defined in (\ref{cov-mat-W}).  In addition, we have 
  \begin{eqnarray}
\nabla^2 U_{2n}(x)&=&\nabla^2 U(x)- \cchi^{\prime}\beta+\cchi^{\prime}~
  \mbox{\rm cov}_{\Ka_{2n}}(x)~  \cchi\nonumber\\
\nabla^2 V_{2n+1}(y)&=&\nabla^2 V(y)-\tau^{-1}+  \cchi~  \mbox{\rm cov}_{\Ka_{2n+1}}(y)~\cchi^{\prime}.\label{ref-nablaUVn1g}
 \end{eqnarray}
\end{prop}

We remark that (\ref{ref-grad-hessian}) and  (\ref{def-UU-VV-e}) yield 
 \begin{eqnarray*}
\nabla \UU(x)&=&\nabla U(x)-\cchi^{\prime}\left\{(\alpha+\beta x)-\mbox{\rm mean}_{K_{\VV}}(x)\right\}\\
\nabla^2 \UU(x)&=&\nabla^2 U(x)- \cchi^{\prime}\beta+\cchi^{\prime}~
  \mbox{\rm cov}_{K_{\VV}}(x)~  \cchi
 \end{eqnarray*}
with the bridge transition $K_{\VV}$ defined in (\ref{def-KVV-KUU}).
Recalling that $\Ka_{2n}=K_{V_{2n}}$ we have
$$
\nabla (\UU- U_{2n})=\cchi^{\prime}(\mbox{\rm mean}_{K_{\VV}}-\mbox{\rm mean}_{K_{V_{2n}}})
\quad\mbox{\rm and}\quad
\nabla^2 (\UU- U_{2n})=\cchi^{\prime}
  \left(\mbox{\rm cov}_{K_{\VV}}- \mbox{\rm cov}_{K_{V_{2n}}}\right)  \cchi.
$$
The convergence of conditional expectations and covariances of Sinkhorn transitions towards the conditional expectations and covariances of  Schr\"odinger bridges
is discussed in~\cite{del2025entropic}, see also Theorem~\ref{theo-cov-cm} in the present article.
We conclude this discussion by noting that applying Proposition~\ref{prop-gh-uv} along with (\ref{UV-Wn}) we obtain the following interplay between the Hessian of Sinkhorn log-densities and conditional covariance functions.
\begin{prop}
For any $n\geq 0$ we have the decompositions
\begin{eqnarray}
\nabla^2_2W_{2n+1}^{\flat}(y,x)&=&\nabla^2U(x)+ \cchi^{\prime}~
  \mbox{\rm cov}_{\Ka_{2n}}(x)~  \cchi\nonumber\\
\nabla^2_2W_{2(n+1)}(x,y)&=&\nabla^2 V(y)+  \cchi~  \mbox{\rm cov}_{\Ka_{2n+1}}(y)~\cchi^{\prime}.\label{HW-UV-0}
\end{eqnarray}
\end{prop}

\section{Stability analysis}\label{sec:stability_analysis}

\subsection{Overview}

In Section~\ref{SS-bridges} we have presented different formulations of Sinkhorn bridges, including 
geometric entropic projection methodologies, operator-theoretic framework with Markov
integral transport formulae
and commutation-type  formulae, as well as
 dual formulations in terms of potential functions.  The stability of Sinkhorn bridges  discussed in this section combines these
  different formulations with the stability theory of Markov operators provided in
 Section~\ref{sec:stab_markov}.

\subsection{Linear decay rates}
 Using (\ref{keysum-p}) and the monotone properties stated in Proposition~\ref{prop-decrease2n01}  for any $Q\in \Pi(\mu,\eta)$  and $n\geq 0$ we have the estimate
\begin{eqnarray*}
\Ha(Q|P)
&\geq &
\sum_{0\leq l\leq  n}\left(\Ha(\eta|\pi_{2l})+\Ha(\mu|\pi_{2l+1})\right)\geq (n+1)~(\Ha(\eta|\pi_{2n})+\Ha(\mu|\pi_{2n+1})).
\end{eqnarray*}
This yields  an elementary and very short proof of the following rather well known theorem.
\begin{theo}\label{th-gibbs-loop-entrop}
For any $n\geq 1$ we have
\begin{eqnarray*}
\Ha(\eta|\pi_{2n})+\Ha(\mu|\pi_{2n+1})\leq \frac{1}{n+1}~\Ha(P_{\mu,\eta}|P)
\end{eqnarray*}
and, in addition,
$$
\lim_{n\rightarrow \infty}n~ \Ha\left(\eta~|~\pi_{2n}\right)=
\lim_{n\rightarrow \infty} n~\Ha\left(\mu~|~\pi_{2n+1}\right)=0.
$$
\end{theo}
The last assertion is a direct consequence of the convergence of the series (\ref{keysum}). 
Sublinear rates have been developed in the articles~\cite{alts-2017,chak-2018,dvu-2017}.
Linear rates with the robust constant $\mbox{\rm Ent}(P_{\mu,\eta}~|~P)$ on possibly non-compact spaces were first obtained by L\'eger in~\cite{leger} using elegant gradient descent and Bregman divergence
 techniques, see also the recent articles~\cite{doucet-bortoli,karimi-2024}. A refined convergence rate at least one order faster has also been developed in~\cite{promit-2022}.

\subsection{Bounded cost functions}

\subsubsection{Universal contraction rates}
The contraction analysis developed in this section is taken from~\cite{adm-25}.
Consider the marginal measures $(\mu,\eta)=(\lambda_U,\nu_V)$ and the Markov transition $K$  defined in
(\ref{ref-intro-UV}) and (\ref{def-Qa}).
We further assume that $\mbox{\footnotesize osc}(W)<\infty$. In this case, using (\ref{Kan}) we check that
\begin{eqnarray}
\frac{d\delta_{x_1}\Ka_{2n}}{d\delta_{x_2}\Ka_{2n}}(y)\frac{d\delta_{x_2}\Ka_{2n}}{d\delta_{x_1}\Ka_{2n}}(z)
&=&\exp{\left((W(x_2,y)-W(x_1,y))+(W(x_1,z)-W(x_2,z))\right)}~~\label{time-hom}\\
&\geq& \epsilon_W:=
\exp{\left(-2~\mbox{\footnotesize osc}(W)\right)}.\nonumber
\end{eqnarray}
The time-homogeneous property (\ref{time-hom}) is sometimes called  ``structure conservation". A first proof of geometric convergence based on these constant cross-product ratios is due to Fienberg in the seminal article~\cite{fienberg-1970}.

The case when one has bounded  transition potentials, $W$, typically arise on finite or compact spaces. For instance, Markov processes on compact manifolds, such as diffusions on bounded manifolds reflected at the boundary,  have Markov transitions of the form (\ref{def-Qa}) with
a uniformly bounded transition potential $W$ that depends on the time parameter~\cite{ado-24,aronson,nash,varopoulos}.  
The convergence analysis of Sinkhorn bridges for uniformly bounded transition potential $W$ is often developed using Hilbert projective metrics~\cite{chen,deligiannidis,franklin,marino}.  As remarked in \cite[Section 3]{cohen}, the Hilbert metric dominates all $\phi$-divergences, thus we expect rather strong conditions to ensure the contraction of Sinkhorn transitions w.r.t.~the Hilbert metric, see for instance~\cite{borwein,chen,deligiannidis,franklin}. Thus contraction estimates
w.r.t.~the Hilbert metric yield exponential decays w.r.t.~$\phi$-divergences, including the total variation norm, the relative entropy and Wasserstein distances of any order.  Nevertheless for the same reason (i.e.~the domination of all $\phi$-divergences) they
can rarely be used to derive contraction estimates with respect to some $\phi$-divergence criteria (unless the criteria is equivalent to the Hilbert metric).

To extend our discussion,  we observe that for any $n\geq 0$ we have
\begin{equation}\label{ref-min-cond-s}
\begin{array}{l}
\delta_{x_1}\Ka_{2n}(dy)~\delta_{x_2}\Ka_{2n}(dz)\geq \epsilon_W~\delta_{x_2}\Ka_{2n}(dy)~\delta_{x_1}\Ka_{2n}(dz)\\
\\
\Longrightarrow \delta_{x_1}\Ka_{2n}\geq 
\epsilon_W~\delta_{x_2}\Ka_{2n}\Longrightarrow
{\sf dob}(\Ka_{2n})\leq 1-\epsilon_W.
\end{array}
 \end{equation}
In a similar manner, for any $n\geq 0$ we check that 
\begin{equation}\label{ref-min-cond-s-y}
 \delta_{y_1}\Ka_{2n+1}\geq 
\epsilon_W~\delta_{y_2}\Ka_{2n+1}\Longrightarrow
{\sf dob}(\Ka_{2n+1})\leq 1-\epsilon_W.
 \end{equation}
Without further work, the Markov transport equations 
presented in (\ref{decrease2n}) and (\ref{decrease2n1}) yield the following theorem.
\begin{theo}\label{tx-st}
For any convex function $\Phi$ satisfying (\ref{Phi-ref}) we have the contraction estimates
\begin{equation}\label{ref-contract-even-i}
\begin{array}{rcccl}
 \Ha_{\Phi}\left(\pi_{2(n+1)},\eta\right)&\leq& (1-\epsilon_W)~  \Ha_{\Phi}\left(\mu ,\pi_{2n+1}\right)&\leq&  (1-\epsilon_W)^2~\Ha_{\Phi}\left(\pi_{2n},\eta \right),\\
 &&&&\\
 \Ha_{\Phi}\left(\eta,\pi_{2(n+1)}\right)&\leq&(1-\epsilon_W)~  \Ha_{\Phi}\left(\pi_{2n+1},\mu \right)&\leq &(1-\epsilon_W)^2~ \Ha_{\Phi}\left(\eta ,\pi_{2n}\right).
 \end{array}
 \end{equation}
The above  inequalities remain valid when we replace $(\eta,\pi_{2n})$ by  $(\mu,\pi_{2n+1})$. 
\end{theo}

\begin{prop}\label{prop-series-radon}
There exists some $c>0$ such that for any $n\geq 1$ we have
$$
\left\Vert
\frac{d\pi_{2n}}{d\eta}-1\right\Vert\vee \left\Vert
\frac{d\eta}{d\pi_{2n}}-1\right\Vert\vee \left\Vert\log{\frac{d\eta}{d\pi_{2n}}}\right\Vert
\leq c~(1-\epsilon_W)^{2n}.
$$
The above  estimates remain valid when we replace $(\pi_{2n},\eta)$ by 
$(\pi_{2n+1},\mu)$. 
\end{prop}
\proof
By (\ref{ref-min-cond-s-y}) for any $n\geq 1$ and $x\in\XX$ and $y_1,y_2\in \YY$ 
  we have
$$
\epsilon_W~ \delta_{y_1}\Sa_{2n}\leq\delta_{y_2}\Sa_{2n}=\int~\delta_{y_2}\Ka_{2n-1}(dz)~\delta_z\Ka_{2n}\leq 
\epsilon_W^{-1}~ \delta_{y_1}\Sa_{2n}.
$$
Using (\ref{pre-gibbs}) and (\ref{fb-sinkhorn}), 
integrating $y_1$ w.r.t. $\eta$
and $y_2$ w.r.t. $\pi_{2(n-1)}$  this implies that
$$
\begin{array}{l}
\epsilon_W~\eta\leq \pi_{2n}\leq \epsilon_W^{-1}~\eta
\\
\\
\displaystyle\Longrightarrow\quad
\frac{d\eta}{d\pi_{2n}}\vee \frac{d\pi_{2n}}{d\eta}\leq \epsilon_W^{-1}
\quad \mbox{\rm and similarly}\quad
\frac{d\mu}{d\pi_{2n+1}}\vee \frac{d\pi_{2n+1}}{d\mu}\leq \epsilon_W^{-1}.
\end{array}$$
By (\ref{pre-gibbs}) and Proposition~\ref{prop-commute-S} for any $k,l\geq 0$ we have
$$
\frac{d\pi_{2 (l+k)}}{d\eta}(y)-1= \left(\delta_y-\eta\right)\Sa_{2(l+k)}\ldots \Sa_{2(l+1)}
\left(\frac{d\pi_{2l}}{d\eta}\right).
$$
On the other hand, (\ref{ref-min-cond-s-y}) and (\ref{ref-min-cond-s}) ensure that for any $n\geq 1$ we have ${\sf dob}(\Sa_{n})\leq (1-\epsilon_W)^2$.
Using (\ref{tv-osc-norm}) and (\ref{defi-beta-Phi}) for 	any $l\geq 1$ and $k\geq 0$ we check that
$$
\left\vert
\frac{d\pi_{2 (l+k)}}{d\eta}(y)-1\right\vert\leq  2\epsilon_W^{-1}~
\left\Vert \left(\delta_y-\eta\right)\Sa_{2(l+k)}\ldots \Sa_{2(l+1)}\right\Vert_{\sf tv}
\leq  2\epsilon_W^{-1}~(1-\epsilon_W)^{2k}.
$$
In the same vein, we have
$$
\begin{array}{l}
\displaystyle\frac{d\pi_{2(l+k)+1}}{d\mu}(x)-1=
(\delta_x-\mu)\Sa_{2(l+k)+1}\ldots \Sa_{2(l+1)+1}\left(\frac{d\pi_{2l+1}}{d\mu}\right)
\\
\\
\displaystyle\Longrightarrow\quad
\left\vert \frac{d\pi_{2(l+k)+1}}{d\mu}(x)-1\right\vert\leq  2\epsilon_W^{-1}~(1-\epsilon_W)^{2k}.
\end{array}$$
By Proposition~\ref{prop-commute} we have
$$
\begin{array}{l}
\displaystyle
\frac{d\eta}{d\pi_{2l}}-1=
\Ka_{2l+1}\left(\frac{d\pi_{2l-1}}{d\mu}-1\right)
\quad\mbox{\rm and}\quad
\frac{d\mu}{d\pi_{2l+1}}-1=
\Ka_{2(l+1)}\left(\frac{d\pi_{2l}}{d\eta}-1\right).
\end{array}$$
Finally, note that for any function $f(z)>0$ we have
\begin{eqnarray}
\left\vert\log{f(z)}\right\vert&=& \log{\left(1+(f(z)-1)\right)}~1_{ f(z)\geq 1}+
\log{\left(1+\frac{1}{ f(z)}-1\right)}~1_{ {1}/{ f(z)}\geq 1}\nonumber\\
&\leq&\left(f(z)-1\right)~1_{ f(z)\geq 1}+\left(\frac{1}{ f(z)}-1\right)~1_{ {1}/{ f(z)}\geq 1}.\label{log-ratio}
\end{eqnarray}
The last assertion follows from the fact that $ \log{(1+u)}\leq u$,  for any $1+u>0$.
\cqfd

\subsubsection{Schr\"odinger bridges}\label{sec-series-UU}
By Proposition~\ref{UVn-prop-series} and
Proposition~\ref{prop-series-radon} there exists some constant $c>0$ such that for any $n\geq 0$ we have
$$
\Vert V_{2n}-\VV\Vert\vee \Vert U_{2n}-\UU\Vert\leq c~(1-\epsilon_W)^{2n}
$$
with the  uniformly absolutely-convergent series
\begin{equation}\label{UU-VV}
\VV=V_{0}-\sum_{l\geq 0 }\log{\left(\frac{d\eta}{d\pi_{2l}}\right)}\quad \mbox{and}\quad
\UU=U_{0}-\sum_{l\geq 0 }\log{\left(\frac{d\mu}{d\pi_{2l+1}}\right)}.
\end{equation}
In addition, we have
$$
\PP(d(x,y)):=e^{-\UU(x)}~e^{-W(x,y)}~e^{-\VV(x)}~
\gamma(d(x,y))\Longrightarrow
\Vert \Pa_{2n}-\PP\Vert_{\sf tv} 
\leq c~(1-\epsilon_W)^{2n}.
$$
We check this claim using the fact that for any probability measures of the form
$$
p_u(dz):=e^{-u(z)}~\gamma(dz)\quad\mbox{\rm
and}\quad p_v(dz):=e^{-v(z)}~\gamma(dz)
$$
for some real valued functions $u(z),v(z)$ we have
\begin{eqnarray*}
\Vert p_u-p_v\Vert_{\sf tv} &=&
\frac{1}{2}\int~\vert e^{-u(z)}-e^{-v(z)}\vert~\gamma(dz)\\
&\leq& \frac{1}{2}~\Vert u-v\Vert
\int~(e^{-u(z)}+e^{-v(z)})~\gamma(dz)= \Vert u-v\Vert.
\end{eqnarray*}
In addition, using (\ref{UVn}) one has
$$
\VV=V+\log{K^{\flat}(e^{-\UU})}\quad\mbox{\rm and}\quad
\UU=U+\log{K(e^{-\VV})}\Longrightarrow\PP=P_{\mu,\eta}.
$$
Also note that
\begin{eqnarray*}
\mathcal{H}(P_{\mu,\eta}~|~\Pa_{2n})  &=&\mu(U_{2n}-\UU)+\eta(V_{2n}-\VV).
\end{eqnarray*}
We summarize the above discussion with the following theorem.
\begin{theo}
There exists some $c>0$ such that for any $n\geq 0$ we have
$$
\mathcal{H}(P_{\mu,\eta}~|~\Pa_{2n})  \vee
  \Vert P_{\mu,\eta}-\Pa_{2n}\Vert_{\sf tv} 
\leq c~(1-\epsilon_W)^{2n}.
$$
In addition, we have the entropic series expansion
\begin{eqnarray*}
\mathcal{H}(P_{\mu,\eta}~|~\Pa_{2n})  
  &=&\sum_{l\geq  n}\left[\Ha(\mu|\pi_{2l+1})+\Ha(\eta|\pi_{2l})\right].
\end{eqnarray*}
\end{theo}

\subsection{Gaussian bridges}\label{Gauss-review-sec}

\subsubsection{Linear Gaussian models}
Consider the state spaces $\XX=\YY=\RR^d$, for some integer $d\geq 1$ and the Lebesgue measure $\lambda(dz)=\nu(dz)=dz$. Hereafter, points $x$ in the Euclidean space $\XX=\RR^d$ are represented by $d$-dimensional column vectors (or, equivalently, by $d \times 1$ matrices). When there is no possible confusion, we use the notation $\Vert\cdot\Vert$ for any equivalent matrix or vector norm.

 Consider the target Gaussian probability measures $
(\mu,\eta)=(\lambda_U,\lambda_V)=(\nu_{m,\sigma},\nu_{\overline{m},\overline{\sigma}})$
discussed in (\ref{KL-def}).
Equivalently, the log-density functions $(U,V)$ are defined by 
\begin{equation}\label{def-U-V}
U(x)=-\log{{\sf gauss}_{\sigma}(x-m)}
\quad \mbox{\rm and}\quad V(y):=-\log{{\sf gauss}_{\overline{\sigma}}(y-\overline{m})}.
\end{equation}
Denote by $G$ a $d$-dimensional centred Gaussian random variable (r.v.) with unit covariance.  Then we have that the linear Gaussian transition $K$ is defined by  (\ref{def-W}) takes the form
\begin{equation}\label{ref-lin-gauss}
K(x,dy):=e^{-W(x,y)}~dy=\Ka_0(x,dy):=\PP(\Za_0(x)\in dy)
\end{equation}
with the random map
\begin{equation}\label{ref-lin-gauss-sigma0}
\begin{array}{l}
\Za_0(x):=\alpha+\beta~x+\tau^{1/2}~G\\
\\
\Longrightarrow \quad \pi_0:=\mu K=\nu_{m_0,\sigma_0}
\quad\mbox{with}\quad m_0:=\alpha+\beta m \quad\mbox{\rm and}\quad
\sigma_0:=\beta\sigma\beta^{\prime}+\tau.
\end{array}
\end{equation}

\subsubsection{Sinkhorn bridges}

We can rewrite $ \Za_0(x)$ into the following alternative form
$$
\Za_0(x):=m_0+\beta_0~(x-m)+\tau_0^{1/2}~G
 \quad\mbox{\rm with}\quad (\beta_0,\tau_0):=(\beta,\tau).
$$
The conjugate formula (\ref{Bayes1}) is defined via
\begin{equation}\label{bayes-update}
\begin{array}{l}
K^{\ast}_{\mu}(y,dx)=\Ka_1(y,dx):=\PP(\Za_1(y)\in dx)\\
\\
\mbox{with the random map}\quad
 \Za_1(y):=m+\beta_1~(y-m_0)+\tau_1^{1/2}~G
 \end{array}
\end{equation}
defined in terms of the gain matrix
$$
\beta_1:=\sigma\beta^{\prime}\sigma_0^{-1}
\quad \mbox{\rm and}\quad
\tau_1^{-1}:=\sigma^{-1}+\beta^{\prime}\tau^{-1}\beta,
$$
In signal processing and filtering literature, this explicit formula is often referred as the Bayesian update formula for the Kalman filter.

Applying the matrix inversion lemma to the covariance matrix $\sigma_0$ defined in (\ref{ref-lin-gauss-sigma0}) we have
$$
\begin{array}{l}
\sigma^{-1}_0=\tau^{-1}-\tau^{-1}\beta~\left(\sigma^{-1}+\beta^{\prime}\tau^{-1}\beta\right)^{-1}~\beta^{\prime}\tau^{-1}=
\tau^{-1}-\tau^{-1}\beta~\tau_1~\beta^{\prime}\tau^{-1}\\
\\
\Longrightarrow
\beta_1=\sigma\beta^{\prime}
\tau^{-1}-\sigma~(\beta^{\prime}\tau^{-1}\beta)~\tau_1\beta^{\prime}\tau^{-1}=
\sigma\beta^{\prime}
\tau^{-1}-\sigma~(\tau_1^{-1}-\sigma^{-1})~\tau_1\beta^{\prime}\tau^{-1}\\
\\
\Longrightarrow
\beta_1=\tau_1~\beta^{\prime}\tau^{-1}.
\end{array}
$$
Then one can obtain $ \Za_1(y)$ as
$$
\Za_{1}(y)
=m_{1}+\beta_{1}~(y-\overline{m})+\tau_{1}^{1/2}~G
\quad \mbox{\rm with}\quad
m_{1}=m+
\beta_{1} (\overline{m}-m_{0})
$$
This implies that $\pi_1=\nu_{\overline{m},\overline{\sigma}}\Ka_1=\nu_{m_1,\sigma_1}$
with the covariance matrix
\begin{equation}\label{ref-lin-gauss-sigma1}
\sigma_{1}:=\beta_{1}~\overline{\sigma}~\beta_{1}^{\prime}+\tau_{1}.
\end{equation}

The term (\ref{s-2}) can be expressed using the regression formula
\begin{eqnarray}
(\Ka_1)^{\ast}_{\eta}(x,dy)=\Ka_2(x,dy):=\PP(\Za_2(x)\in dy) \label{bayes-update-sec}\\
 \Za_2(x):=\overline{m}+\beta_2~(x-m_1)+\tau_2^{1/2}~G \label{bayes-update-sec_1}
\end{eqnarray}
where
$$
\beta_2:=\overline{\sigma}\beta_1^{\prime}\sigma_1^{-1}
\quad \mbox{\rm and}\quad
\tau_2^{-1}:=\overline{\sigma}^{-1}+\beta_1^{\prime}\tau_1^{-1}\beta_1.
$$
Here $ \Za_2(x)$ is given by
$$
\Za_{2}(x)
=m_{2}+\beta_{2}~(x-m)+\tau_{2}^{1/2}~G
\quad \mbox{\rm with}\quad
m_{2}=\overline{m}+
\beta_{2} (m-m_{1}).
$$
Applying the matrix inversion lemma to the covariance matrix $\sigma_1$ defined in  (\ref{ref-lin-gauss-sigma1}) we have
$$
\begin{array}{l}
\sigma^{-1}_1=\tau^{-1}_1-\tau_1^{-1}\beta_1~\left(\overline{\sigma}^{-1}+\beta_1^{\prime}\tau_1^{-1}\beta_1\right)^{-1}~\beta_1^{\prime}\tau_1^{-1}=
\tau^{-1}_1-\tau^{-1}_1\beta_1~\tau_2~\beta_1^{\prime}\tau_1^{-1}\\
\\
\Longrightarrow
\beta_2=\overline{\sigma}\beta_1^{\prime}
\tau_1^{-1}-\overline{\sigma}~(\beta_1^{\prime}\tau_1^{-1}\beta_1)~\tau_2\beta_1^{\prime}\tau_1^{-1}=
\overline{\sigma}\beta_1^{\prime}
\tau^{-1}_1-\overline{\sigma}~(\tau_2^{-1}-\overline{\sigma}^{-1})~\tau_2\beta_1^{\prime}\tau^{-1}_1\\
\\
\Longrightarrow
\beta_2=\tau_2~\beta_1^{\prime}~\tau_1^{-1}=\tau_2~\tau^{-1}\beta.
\end{array}
$$

Iterating the Bayes update formulae (\ref{bayes-update}),  \eqref{bayes-update-sec}, \eqref{bayes-update-sec_1}  
(see section 4 and \cite[Appendix C]{adm-24}) the Sinkhorn transitions  can be expressed as
$$
\Ka_{2n}(x,dy)=\PP\left(\Za_{2n}(x)\in dy\right)\quad\mbox{\rm and}\quad
\Ka_{2n+1}(y,dx)=\PP\left(\Za_{2n+1}(y)\in dx\right)
$$
with linear-Gaussian random maps defined, for any $n\geq 0$, by
\begin{eqnarray}
\Za_{2n}(x)&=&m_{2n}+\beta_{2n}(x-m)+\tau_{2n}^{1/2}~G \nonumber\\
\Za_{2n+1}(y)
&=&m_{2n+1}+\beta_{2n+1}(y-\overline{m})+\tau_{2n+1}^{1/2}~G.\label{sinkhorn-maps}
\end{eqnarray}
The parameters in (\ref{sinkhorn-maps}) are described sequentially by the following theorem. 
\begin{theo}[\cite{adm-24}]
The gain $\beta_n$ and covariance matrices $\tau_n$ in the Bayes update formulae (\ref{sinkhorn-maps}) are defined for any $n\geq 0$ by the recursion
\begin{equation}\label{tau-rec-0}
\begin{array}{rclcrcl}
\beta_{2n}&=&\tau_{2n}~\tau^{-1}\beta&\quad&
\beta_{2n+1}&=&\tau_{2n+1}~\beta^{\prime}\tau^{-1}.
\\
&&&&&&\\
\tau_{2n+1}^{-1}&=&\sigma^{-1}+\beta_{2n}^{\prime}~\tau_{2n}^{-1}~\beta_{2n}
&\quad&
\tau_{2(n+1)}^{-1}&=&\overline{\sigma}^{-1}+\beta_{2n+1}^{\prime}\tau_{2n+1}^{-1}\beta_{2n+1}
\end{array}
\end{equation}
The parameters $m_n$ in (\ref{sinkhorn-maps}) are also defined, for any $n\geq 0$, by the recursion
\begin{eqnarray}
m_{2n+1}&=&m+
\beta_{2n+1} (\overline{m}-m_{2n}), \nonumber\\
m_{2(n+1)}&=&\overline{m}+\beta_{2(n+1)} ~(m-m_{2n+1})\label{sinkhorn-mean-maps}
\end{eqnarray}
with initial condition
$
m_0=\alpha_0+\beta_0 m$, where, again, $(\alpha_0,\beta_0,\tau_0)=(\alpha,\beta,\tau)$ are parameters of the reference linear Gaussian transition (\ref{ref-lin-gauss}). 
\end{theo}

Note that
$
\pi_n=\nu_{m_{n},\sigma_n}
$
with the covariance matrices
$$\sigma_{2n}:=
\beta_{2n}~\sigma~\beta_{2n}^{\prime}+\tau_{2n}
\quad \mbox{and}\quad\sigma_{2n+1}:=\beta_{2n+1}~\overline{\sigma}~\beta_{2n+1}^{\prime}+\tau_{2n+1}.
$$

\begin{examp}
In the context of stochastic flows
discussed in (\ref{def-diff}) with a linear drift $ b_s(x)=Ax$
 and an homogeneous diffusion function
 $ \sigma_s(x)=\Sigma
 $,
for some $A\in \RR^{d\times d}$ and $\Sigma\in \Sa_d^+$, for any given time horizon $t>0$ one has
\begin{equation}\label{lin-gauss}
T_t(x,dy)=\PP(X_{t}(x)\in dy)=e^{-\Wa_t(x,y)}~dy
\quad \mbox{\rm with}\quad
\Wa_t(x,y):=-\log{{\sf gauss}_{\tau(t)}(y-\beta(t) x)}
\end{equation}
and the parameters
$$
(\beta(t),\tau(t))=\left(e^{tA}, \int_0^{t}e^{sA}~\Sigma~e^{sA^{\prime}}~ds\right).
$$
In this context we can choose the reference transition (\ref{ref-lin-gauss}) with
$W=\Wa_t$.
\end{examp}

\subsubsection{Riccati equations}\label{ricc-equations-sec}

Consider the rescaled covariance matrices
\begin{equation}\label{rescaled-tau}
\upsilon_{2n}:=
\overline{\sigma}^{-1/2}~\tau_{2n}~\overline{\sigma}^{-1/2}
\quad\mbox{and}\quad
\upsilon_{2n+1}:=
\sigma^{-1/2}~\tau_{2n+1}~\sigma^{-1/2}
\end{equation}
Uisng (\ref{tau-rec-0})  we check that
\begin{eqnarray*}
\tau_{2n+1}^{-1}&=&\sigma^{-1}+(\tau_{2n}~\tau^{-1}\beta)^{\prime}~\tau_{2n}^{-1}~(\tau_{2n}~\tau^{-1}\beta)=\sigma^{-1}+(\tau^{-1}\beta)^{\prime}~\tau_{2n}^{-1}(\tau^{-1}\beta)\\
\tau_{2(n+1)}^{-1}&=&\overline{\sigma}^{-1}+(\tau_{2n+1}~\beta^{\prime}\tau^{-1})^{\prime}~\tau_{2n+1}^{-1}~(\tau_{2n+1}~\beta^{\prime}\tau^{-1})=\overline{\sigma}^{-1}+(\tau^{-1}\beta)~
\tau_{2n+1}^{-1}~(\tau^{-1}\beta)^{\prime}
\end{eqnarray*}
This yields the recursions
\begin{equation}\label{k-ricc-0}
\upsilon_{2n+1}^{-1}=
I+  \gamma^{\prime}~\upsilon_{2n}~ \gamma
\quad \mbox{and}\quad
\upsilon_{2(n+1)}^{-1}=I+ \gamma~\upsilon_{2n+1}~  \gamma^{\prime}
\end{equation}
with
$$
\gamma:=\overline{\sigma}^{1/2}~\tau^{-1}\beta~\sigma^{1/2}.
$$

Next, observe that
\begin{equation}\label{proof-ricc}
u^{-1}=I+\gamma^{\prime}v\gamma
\Longrightarrow
(I+\gamma u\gamma^{\prime})^{-1}=\left(I+\left(\varpi+v\right)^{-1}\right)^{-1},~~\text{with}\quad\varpi:=\left(\gamma\gamma^{\prime}\right)^{-1}.
\end{equation}
To verify this claim, we use the  matrix inversion lemma to prove that 
$$
\begin{array}{l}
u=I-\gamma^{\prime}~(v^{-1}+\gamma\gamma^{\prime})^{-1}\gamma.
\\
\\
\Longrightarrow
\gamma u \gamma^{\prime}=\gamma\gamma^{\prime}-\gamma\gamma^{\prime}~(v^{-1}+\gamma\gamma^{\prime})^{-1}\gamma
\gamma^{\prime}=\left(\left(\gamma\gamma^{\prime}\right)^{-1}+v\right)^{-1}=
\left(\varpi+v\right)^{-1}
\end{array}$$
Applying the above to $(u,v)=(\upsilon_{2n+1},\upsilon_{2n})$ and then to $(u,v)=(\upsilon_{2n},\upsilon_{2n-1})$ one has.

\begin{theo}[\cite{adm-24}]
For any $n\geq 1$ we have
\begin{equation}\label{k-ricc}
\upsilon_{2n}=\mbox{\rm Ricc}_{\varpi}\left(\upsilon_{2(n-1)}\right)\quad
\mbox{and}\quad
\upsilon_{2n+1}=\mbox{\rm Ricc}_{ \overline{\varpi}}\left(\upsilon_{2n-1}\right),
\end{equation} 
with the Riccati maps $(\mbox{\rm Ricc}_{\varpi},\mbox{\rm Ricc}_{ \overline{\varpi}})$ defined in (\ref{ricc-maps-def}) and the matrices
\begin{equation}\label{def-varpi}
 ( \varpi^{-1},\overline{\varpi}^{-1}):=(\gamma\,\gamma^{\prime},\gamma^{\prime}\gamma).
\end{equation} 
\end{theo}
The exponential decays of Riccati matrices towards the positive fixed points of the Riccati maps follows the contraction analysis discussed in (\ref{cv-ricc-intro}).
In addition, Corollary 4.4 in~\cite{adm-24} yields for any $n\geq 1$ the uniform estimates 
\begin{equation}\label{ue-ref-locc}
\overline{\sigma}^{1/2}(I+\varpi^{-1})^{-1}\overline{\sigma}^{1/2}\leq \tau_{2n}\leq \overline{\sigma}\quad\mbox{and}\quad
\sigma^{1/2}(I+\overline{\varpi}^{-1})^{-1}\sigma^{1/2}\leq \tau_{2n+1}\leq \sigma.
\end{equation}

\subsubsection{Schr\" odinger bridges}

Note that  the transitions of  the Gibbs-loop process  introduced in (\ref{fb-sinkhorn})
for even indices can be rewritten as
$$
 \Sa_{2n}(y_1,dy_2):=\PP\left(\Za^{\circ}_{2n}(y_1)\in dy_2\right)
$$
with the random maps
\begin{eqnarray*}
\Za^{\circ}_{2n}(y)&:=&\overline{m}+\beta^{\circ}_{2n}(y-\overline{m})+\left(\tau^{\circ}_{2n}\right)^{1/2}~G
\end{eqnarray*}
defined in terms of the parameters
$$
\begin{array}{rclcrcl}
\beta^{\circ}_{2n}&:=&\beta_{2n}~\beta_{2n-1},&&
\text{and}~~\tau^{\circ}_{2n}&:=&\tau_{2n}+
\beta_{2n}~\tau_{2n-1}~\beta_{2n}^{\prime}.
\end{array}$$
Consider the directed matrix products
\begin{equation}\label{directed-prod}
\beta^{\circ}_{2n,0}:=\beta^{\circ}_{2n}\beta^{\circ}_{2(n-1)}\ldots \beta^{\circ}_{2}\end{equation}
Then we have
\begin{equation}\label{directed-prod-sig}
m_{2n}-\overline{m}=\beta^{\circ}_{2n} ~(m_{2(n-1)}-\overline{m})
\quad\mbox{\rm and}\quad
\sigma_{2n}-\overline{\sigma}=\beta^{\circ}_{2n,0}(\sigma_0-\overline{\sigma})\left(\beta^{\circ}_{2n,0}\right)^{\prime}.
\end{equation}
The l.h.s.~assertion in (\ref{directed-prod-sig}) is a direct consequence of (\ref{sinkhorn-mean-maps}) and the
detailed proof of the  covariance formula (\ref{directed-prod-sig}) is provided in Appendix C in~\cite{adm-24}. By \cite[Lemma 4.9]{adm-24}, we have
$$
0\leq \overline{\sigma}^{-1/2}~\beta^{\circ}_{2n}~\overline{\sigma}^{1/2}=I-\upsilon_{2n}\leq (I+\varpi)^{-1}.
$$
This shows that
$$
\lim_{n\rightarrow\infty}\Vert \beta^{\circ}_{2n,0}\Vert=0
\quad\mbox{\rm and therefore}\quad
(m_{2n},\sigma_{2n})\longrightarrow_{n\rightarrow\infty}(\overline{m},\overline{\sigma})
$$
More refined estimates are provided in Section~\ref{cv-theo-sink-gauss}.

In terms of the rescaled matrices (\ref{rescaled-tau})
combining (\ref{sinkhorn-maps}) with (\ref{tau-rec-0}) we check that
$$
\Za_{2n}(x)=m_{2n}+\tau_{2n}~\tau^{-1}\beta~(x-m)+\tau_{2n}^{1/2}~G
\quad\mbox{\rm with}\quad
\tau_{2n}=\overline{\sigma}^{1/2}~\upsilon_{2n}~
\overline{\sigma}^{1/2}
$$
with the solution $\upsilon_{2n}$ of the Riccati equation (\ref{k-ricc}). Taking the limit as $n\rightarrow\infty$ one has the following theorem.

\begin{theo}[\cite{adm-24}]\label{Th1}
The Schr\"odinger bridge from $\mu=\nu_{m,\sigma}$ to $\eta=\nu_{\overline{m},\overline{\sigma}}$ with respect to the reference measure 
 $P:=(\mu\times K)$ is given by the formula
 $$
 P_{\mu,\eta}(d(x,y))=\mu(dx)~\PP\left(Z(x)\in dy\right)
 $$
 with the random map
\begin{equation}\label{bridge-map-i}
Z(x)=\overline{m}+\varsigma~\tau^{-1}~\beta~\left(x-m\right)+\varsigma^{1/2}~G
\end{equation}
and the parameters $\varsigma:=\overline{\sigma}^{1/2}~r~\overline{\sigma}^{1/2}$  and 
the matrix $r=\mbox{\rm Ricc}_{\varpi}(r)$ as in (\ref{def-fix-ricc-1}).
\end{theo}

The transport property stems from the equivalences
\begin{eqnarray*}
(\varsigma~\tau^{-1}~\beta)~\sigma~(\varsigma~\tau^{-1}~\beta)^{\prime}+\varsigma_{\theta}=\overline{\sigma}
&\Longleftrightarrow & \varsigma+\varsigma~\left(\overline{\sigma}^{1/2}~\varpi~\overline{\sigma}^{1/2}\right)^{-1}~\varsigma=\overline{\sigma}\nonumber\\
&\Longleftrightarrow & r~\varpi^{-1}~r+ r=I\Longleftrightarrow (\ref{equiv-r}).
\end{eqnarray*}
When $\beta=I$ and $\tau$ is diagonal the above formula reduces to \cite[Eq (2)]{bunne2023schrodinger}. 
Theorem~\ref{Th1} can be seen as a simplification and a generalization of 
Schr\" odinger bridge formulae recently presented in the series of articles~\cite{agueh,bojitov,loubes,janadi,mallasto} when the drift matrix $\beta$ is arbitrary and $\tau$ is an arbitrary positive definite matrix.

\subsubsection{A convergence theorem}\label{cv-theo-sink-gauss}

We only consider Sinkhorn bridges $\Pa_{2n}$ associated  with even indices, the analysis of the bridges associated with odd indices follows the same line of arguments, thus it is skipped
(further details can be found in Section 5 in~\cite{adm-24}). 

Recall that Sinkhorn random maps $\Za_{2n}(x)$ defined in (\ref{sinkhorn-maps})
and the random map $Z(x)$ of the Schr\"odinger bridge defined in (\ref{bridge-map-i})
take the form
\begin{equation}\label{bridge-map-ib}
\begin{array}{rclcrcl}
\Za_{2n}(x)&=&m_{2n}+\tau_{2n}~\tau^{-1}\beta~(x-m)+\tau_{2n}^{1/2}~G
& \mbox{\rm with}&
\tau_{2n}&:=&\overline{\sigma}^{1/2}~\upsilon_{2n}~\overline{\sigma}^{1/2}
\\
&&&&&&\\
Z(x)&=&\overline{m}+\varsigma~\tau^{-1}~\beta~\left(x-m\right)+\varsigma^{1/2}~G
& \mbox{\rm with}&
\varsigma&:=&\overline{\sigma}^{1/2}~r~\overline{\sigma}^{1/2}.
\end{array}
\end{equation}
In addition, the matrix flow $\upsilon_{2n}$ satisfies the Riccati equation 
\begin{equation}\label{ricc-lin-gauss}
\upsilon_{2n}=\mbox{\rm Ricc}_{\varpi}\left(\upsilon_{2(n-1)}\right)\longrightarrow_{n\rightarrow\infty}
\mbox{\rm Ricc}_{\varpi}\left(r\right)=r.
\end{equation}
with the fixed point $\mbox{\rm Ricc}_{\varpi}\left(r\right)=r$ of the Riccati map given in closed form in (\ref{def-fix-ricc-1}).  The next result is a consequence of (\ref{cv-ricc-intro}) and 
the Ando-Hemmen inequality (\ref{square-root-key-estimate}).

\begin{theo}[\cite{adm-24}]\label{theo-qs} 
There exists some  $c<\infty$ such that for any $n\geq 0$ we have
\begin{eqnarray*}
\Vert \tau_{2n}-\varsigma\Vert\vee \Vert  \tau_{2n}^{1/2}-\varsigma^{1/2}\Vert
&\leq &c~(1+\lambda_{\text{\rm min}}(r+\varpi))^{-2n}~
\Vert \tau-\varsigma\Vert
\end{eqnarray*}
with  $(\varpi,r,\varsigma)$ as in (\ref{def-varpi}) and (\ref{bridge-map-i}). In addition, we have
$$
\Vert\overline{\sigma}^{-1/2}~\beta^{\circ}_{2n,0}~\overline{\sigma}^{1/2}\Vert\leq  c~(1+\lambda_{\text{\rm min}}(r+\varpi))^{-n}.
$$
\end{theo}
Theorem~\ref{theo-qs} can be used to derive a variety of quantitative estimates. For instance, using (\ref{KL-def}) we check the  entropy formula
\begin{equation}\label{ent-intro-even}
\begin{array}{l}
2\mathcal{H}\left(\Pa_{2n}~|~P_{\mu,\eta}\right) \\
\\
=  D_{\sf burg}(\tau_{2n}~|~\varsigma)+\Vert\varsigma^{-1/2}\left(m_{2n}-\overline{m}\right)\Vert_2^2+\Vert\varsigma^{-1/2}
(\tau_{2n}-\varsigma)~(\tau^{-1}\beta\sigma^{1/2})\Vert_F^2
\end{array}
\end{equation}
with the Burg distance $D_{\sf burg}$ defined in (\ref{burg-def}) and the Schr\"odinger bridge $P_{\mu,\eta}$ defined in (\ref{bridge-map-i}).  For a more detailed discussion we refer the reader to \cite[Section 5.2]{adm-24}.
We also note that by \cite[Lemma B.1]{adm-24} (see also \cite[Section 11]{dmt-18}) we have
$$
 \displaystyle\Vert \tau_{2n}-\varsigma\Vert_F
\leq \frac{1}{2 \Vert \varsigma^{-1}\Vert_F}\Longrightarrow
D_{\sf burg}(\tau_{2n}~|~\varsigma)\leq \frac{5}{2}~
  \left\Vert \varsigma^{-1}\right\Vert_F~\left\Vert \tau_{2n}-\varsigma\right\Vert_F.
$$
\begin{cor}
There exists some finite $c<\infty$ and some $n_0\geq 0$ such that for any $n\geq n_0$ we have the entropy estimates
$$
\mathcal{H}\left(\Pa_{2n}~|~P_{\mu,\eta}\right)  \leq c~(1+\lambda_{\text{\rm min}}(r+\varpi))^{-2n}~\left(
\Vert \tau-\varsigma\Vert+\Vert m_{0}-\overline{m}\Vert^2\right).
$$
\end{cor}
\begin{cor}\label{cor-was-sinkhorn}
 For any $p\geq 1$ there exists some finite $c<\infty$ such that for any $n\geq 0$ we have
 the $p$-Wasserstein estimates
$$
\begin{array}{l}
\Da_p\left(\Pa_{2n}~|~P_{\mu,\eta}\right)
\leq 
c~(1+\lambda_{\text{\rm min}}(r+\varpi))^{-2n}
\Vert \tau-\varsigma\Vert_F+c~(1+\lambda_{\text{\rm min}}(r+\varpi))^{-n}~\Vert m_{0}-\overline{m}\Vert_F.
\end{array}
$$
\end{cor}

\subsection{Log-concave-at-infinity models}
Consider the marginal measures $(\mu,\eta)$ the Markov transition $K$ discussed in Section~\ref{sch-pot} on the state spaces  $\XX=\RR^d=\YY$ equipped with the Lebesgue measure $\lambda(dz)=dz=\nu(dz)$. Also assume that the log density $W$ is differentiable. In this case,  
by (\ref{UV-Wn}) for any $n\geq 0$ we have
\begin{eqnarray}
  \nabla_2W_{2n}(x_1,y)-  \nabla_2W_{2n}(x_2,y)&=&
   \nabla_2W(x_1,y)- \nabla_2W(x_2,y)\nonumber
\\
  \nabla_2W^{\flat}_{2n+1}(y_1,x)-  \nabla_2W^{\flat}_{2n+1}(y_2,x)&=&
   \nabla_2W^{\flat}(y_1,x)- \nabla_2W^{\flat}(y_2,x).
\label{nabla1-WN}
\end{eqnarray}
To simplify the presentation, we restrict our attention to linear Gaussian transition potential of the form (\ref{def-W}) and assume that 
 $U$ and $V$ satisfy the strongly convex outside a ball condition (\ref{ex-lg-cg-v2}).
 This condition ensures that the probability measures  $(\mu,\eta)$ 
 satisfy  the $LS(\rho)$ inequality and thus the $\TT_2(\rho)$ inequality 
for some $\rho>0$. More general situations are discussed in~\cite{dp-25}.
 
We have that
\begin{eqnarray*}
   \nabla_2W(x_1,y)- \nabla_2W(x_2,y)&=&\tau^{-1}\beta~(x_2-x_1)
   \\
      \nabla_2W^{\flat}(y_1,x)- \nabla_2W^{\flat}(y_2,x)&=&\beta^{\prime}\tau^{-1}~(y_2-y_1).
\end{eqnarray*}
Thus,  for any $n\geq 0$ we have
\begin{equation}\label{Jlip}
 \Ja(\delta_{z_1}\Ka_{n}~|~\delta_{z_2}\Ka_{n})
 \leq \kappa^2~\Vert z_1-z_2\Vert_2^2\quad \mbox{\rm with}\quad
 \kappa:=\Vert\tau^{-1}\beta\Vert_2.
\end{equation}
 This shows that Sinkhorn transitions $\Ka_n$ satisfy the $\JJ(\Vert\tau^{-1}\beta\Vert_2)$ inequality.

\begin{examp}\label{ex-BB-OU}
In the context  of the linear-Gaussian transitions 
discussed in (\ref{def-diff}):\\
\noindent $\bullet$ When $A=0$ we have 
\begin{equation}\label{pure-BB}
\beta(t)=I\quad\mbox{and}\quad\tau(t)=t~\Sigma
\Longrightarrow\kappa(t):=\Vert\tau(t)^{-1}\beta(t)\Vert_2=\frac{1}{t}~\Vert\Sigma^{-1}\Vert_2.
\end{equation}
$\bullet$ Whenever $A$ is Hurwitz there exists some parameters $c,a>0$ such that
\begin{equation}\label{OU-ref}
\begin{array}{l}
\displaystyle\Vert\beta(t)\Vert_2\leq c~e^{-at} \quad\mbox{and}\quad \forall t\geq t_0\quad\tau(t_0)\leq \tau(t)\leq \tau^+ :=\int_0^{\infty}\beta(s)~\Sigma~\beta(s)^{\prime}~ds\\
\\
\Longrightarrow\quad \exists c_0>0~:~\forall t\geq t_0\qquad
\kappa(t):=\Vert\tau(t)^{-1}\beta(t)\Vert_2\leq c_0~e^{-at}.
\end{array}
\end{equation}
 \end{examp}

\subsubsection{Curvature estimates}
Assume that $(U,V)$ are twice differentiable.  Using
(\ref{ref-grad-hessian}) and (\ref{UVn}) for any $n\geq 0$ we also check that
\begin{eqnarray}
\nabla^2 V_{2n+1}(y)&=&
\nabla^2 V(y) -\Ka_{2n+1}(W^{\flat}_{y,y})(y)+  \mbox{\rm cov}_{\Ka_{2n+1},W^{\flat}}(y)
\nonumber\\
 \nabla^2 U_{2(n+1)}(x)&=&\nabla^2U(x)
 -\Ka_{2(n+1)}(W_{x,x})(x)+  \mbox{\rm cov}_{\Ka_{2(n+1)},W}(x).\label{GHn}
\end{eqnarray}
Combining the above with (\ref{UV-Wn}) we have
\begin{eqnarray*}
\nabla^2_2W_{2(n+1)}(x,y)&=&\nabla^2 V(y) +W^{\flat}_{y,y}(x)-\Ka_{2n+1}(W^{\flat}_{y,y})(y)+  \mbox{\rm cov}_{\Ka_{2n+1},W^{\flat}}(y)
\end{eqnarray*}
with the conditional covariance functions $\mbox{\rm cov}_{\Ka_{2n+1},W^{\flat}}$ defined in (\ref{ref-cond-cov}).
One also has
$$
\nabla^2_2W_{1}^{\flat}(y,x)=W_{x,x}(y)+\nabla^2U(x).
$$
In addition, for $n\geq 1$ we have
\begin{eqnarray*}
\nabla^2_2W_{2n+1}^{\flat}(y,x)&=&\nabla^2U(x)+W_{x,x}(y)
 -\Ka_{2n}(W_{x,x})(x)+  \mbox{\rm cov}_{\Ka_{2n},W}(x).
\end{eqnarray*}
For linear Gaussian transition potential of the form (\ref{def-W}) the above decompositions resume to (\ref{HW-UV-0}). In this situation, for any $n\geq 1$ we  have
\begin{eqnarray}
\nabla^2_2W_{2n-1}^{\flat}(y,x)\geq \nabla^2U(x)\quad \mbox{\rm and}\quad
\nabla^2_2W_{2n}(x,y)\geq 
\nabla^2 V(y).\label{HW-UV}
\end{eqnarray}
Thus, for linear Gaussian transition  (\ref{def-W}) potentials  $(U,V)$ satisfying the convex-at-infinity condition (\ref{ex-lg-cg-v2}),
the Markov transitions $\Ka_n$ with $n\geq 1$ satisfy  the $LS(\rho)$ and thus the $\TT_2(\rho)$ inequality. In the context of linear-Gaussian models, compare (\ref{HW-UV}) with (\ref{sinkhorn-maps}) and the uniform estimates  (\ref{ue-ref-locc}).

Applying Lemma~\ref{lem-Int-log-sob} to Sinkhorn transition $\Ka_{2n}=K_{V_{2n}}$ and  the Schr\"odinger bridge transition
$L=K_{\VV}$ defined respectively in (\ref{Kan}) and (\ref{def-KVV-KUU}) one can prove the following result.
\begin{theo}[\cite{del2025entropic}]\label{theo-cov-cm}
For any $n\geq 1$, we have the bias and covariance entropic bounds estimate
 $$
\int~\mu(dx)~  \Vert\mbox{\rm mean}_{K_{V_{2n}}}(x)-\mbox{\rm mean}_{K_{\VV}}(x)\Vert_{2}^2\leq 
 2\rho~\Ha\left(P_{\mu,\eta}~|~\Pa_{2n}\right)
 $$
 as well as the covariance entropic bounds
 $$
 \begin{array}{l}
\displaystyle \int~\mu(dx)~
 \Vert \mbox{\rm cov}_{K_{V_{2n}}}(x)-\mbox{\rm cov}_{K_{\VV}}(x)\Vert_F\\
 \\
\displaystyle\leq 4\rho~\Ha\left(P_{\mu,\eta}~|~\Pa_{2n}\right)+c~\left(8\rho~\Ha\left(P_{\mu,\eta}~|~\Pa_{2n}\right)\right)^{1/2}\quad\mbox{with}\quad
c^2:=\mu\left(\tr(\mbox{\rm cov}_{K_{\VV}})\right)
 \end{array}
 $$
\end{theo}

\subsubsection{Contraction inequalities}

Using  (\ref{def-entropy-note}) and (\ref{bridge-Pn}) for any $n\geq 0$
we check that
$$
P_{\mu,\eta}=(\Pa_{2n})_{\mu,\eta}
 \quad \mbox{\rm and}\quad
P^{\flat}_{\eta,\mu}=
(\Pa_{2n+1}^{\flat})_{\eta,\mu}.
$$
Recalling that
$$
\pi_{2n}\Ka_{2n+1}=\mu \quad \mbox{\rm and}\quad
 \pi_{2n-1} \Ka_{2n}=\eta
$$
we also have the bridge map formulae
$$
\begin{array}{rclcrcl}
\Pa_{2n}&=&\mu\times \Ka_{2n}=(\Pa_{2n})_{\mu, \mu \Ka_{2n}}
&&
P_{\mu,\eta}&=&(\Pa_{2n})_{\mu,\eta}=
(\Pa_{2n})_{\mu, \pi_{2n-1} \Ka_{2n}}
\\
&&&&&&\\
\Pa^{\flat}_{2n+1}&=&\eta\times \Ka_{2n+1}=\left(\Pa^{\flat}_{2n+1}\right)_{\eta,\eta \Ka_{2n+1}}
&&
 P^{\flat}_{\eta,\mu}
 &=& \left(\Pa^{\flat}_{2n+1}\right)_{\eta,\mu}=
 \left(\Pa^{\flat}_{2n+1}\right)_{\eta,\pi_{2n} \Ka_{2n+1}}.
\end{array}$$
Recall that  the probability measures $(\mu,\eta)$ as well as the collection of Sinkhorn transitions $\Ka_n$
 satisfy  for any $n\geq 1$ the $LS(\rho)$ and thus the $\TT_2(\rho)$ inequalities for some common constant $\rho$. However,  by (\ref{Jlip}),
 Sinkhorn transitions $\Ka_n$ also satisfy the
  Fisher-Lipschitz inequality  (\ref{reflocllip}) with  $\kappa:=\Vert\tau^{-1}\beta\Vert_2$. 
  
 Next we present some direct consequences of Theorem~\ref{key-theo}:
\begin{itemize} 
\item{Applying (\ref{H2WmuT2}) to $(P,K,\eta)=(\Pa_{2n},\Ka_{2n},\pi_{2n-1} )$ we check that
\begin{equation}\label{pi2nP-v2}
\Ha(\eta|\pi_{2n})
\leq \Ha\left(P_{\mu,\eta}~|~\Pa_{2n}\right)
\leq  \epsilon_{\kappa}(\rho) ~\Ha(\pi_{2n-1}|\mu)\quad \mbox{with}\quad \epsilon_{\kappa}(\rho):=(\kappa\rho)^{2}
\end{equation}
Thus, by (\ref{pi2nP}) we have
\begin{eqnarray}
\Ha(P_{\mu,\eta}|\Pa_{2n})&\leq& \Ha(P_{\mu,\eta}|\Pa_{2(n-1)})-\Ha(\pi_{2n-1}|\mu)
\nonumber\\
&\leq & \Ha(P_{\mu,\eta}|\Pa_{2(n-1)})-\epsilon_{\kappa}(\rho)^{-1}~\Ha(P_{\mu,\eta}|\Pa_{2n})\label{th-stab-2n}
\end{eqnarray}
Similarly, by (\ref{H2WKT2}) we have the $2$-Wasserstein estimates
\begin{equation}\label{pi2nP-D2}
\Da_2(\pi_{2n},\eta)
\leq  \kappa\rho~\Da_2(\pi_{2n-1},\mu).
\end{equation}}
\item{Applying (\ref{H2WmuT2}) and (\ref{H2WKT2}) to $(P,K,\eta)=(\Pa^{\flat}_{2n+1},\Ka_{2n+1},\pi_{2n})$ we check that
\begin{equation}\label{pi2n1P-v2}
\Ha(\mu|\pi_{2n+1})
\leq \Ha\left( P_{\mu,\eta}~|~\Pa_{2n+1}\right)
\leq  \epsilon_{\kappa}(\rho)  ~\Ha(\pi_{2n}|\eta)
\end{equation}
Thus, by (\ref{pi2n1P}) we have
\begin{eqnarray}
\Ha(P_{\mu,\eta}|\Pa_{2n+1})&\leq & \Ha(P_{\mu,\eta}|\Pa_{2n-1})-\Ha(\pi_{2n}|\eta)\nonumber\\
&\leq & \Ha(P_{\mu,\eta}|\Pa_{2n-1})- \epsilon_{\kappa}(\rho)^{-1} \Ha\left( P_{\mu,\eta}~|~\Pa_{2n+1}\right).\label{th-stab-2n1}
\end{eqnarray}
}
\end{itemize}

By (\ref{H2WKT2}) we have the $2$-Wasserstein estimates
\begin{equation}\label{pi2n1P-D2}
\Da_2(\pi_{2n+1},\mu)
\leq  \kappa\rho~\Da_2(\pi_{2n},\eta).
\end{equation}
Using (\ref{th-stab-2n}) and (\ref{th-stab-2n1})   one has the following theorem.
\begin{theo}\label{theo-stab-rev}
For any $n\geq 0$ we have  the entropic exponential decays
\begin{equation}\label{theo-edecay}
\Ha(P_{\mu,\eta}|\Pa_{2n})\vee 
\Ha(P_{\mu,\eta}|\Pa_{2n+1})\leq \left(1+1/\epsilon_{\kappa}(\rho)\right)^{-n}~\Ha(P_{\mu,\eta}|\Pa_{0})
\end{equation}
with $\epsilon_{\kappa}(\rho)$ as in (\ref{pi2nP-v2}).
In addition, when $\epsilon_{\kappa}(\rho)<1$ we have  Wasserstein exponential decays 
\begin{equation}\label{theo-d2decay}
\Da_2(\pi_{2n},\eta)
\leq  \epsilon_{\kappa}(\rho)^{n}~\Da_2(\pi_{0},\eta)
\quad\mbox{and}\quad
\Da_2(\pi_{2n+1},\mu)
\leq \epsilon_{\kappa}(\rho)^{n}~\Da_2(\pi_{1},\mu).
\end{equation}
\end{theo}
More refined estimates are provided in~\cite{dp-25}, see also Section~\ref{mr-sec}.
\begin{examp}
In the context  of the linear-Gaussian transitions 
discussed in (\ref{def-diff}) and in Example~\ref{ex-BB-OU} we have
$$
(\ref{pure-BB})\Longrightarrow
\epsilon_{\kappa(t)}(\rho)^{1/2}=
\kappa(t)\rho=\frac{\rho}{t}~\Vert\Sigma^{-1}\Vert_2
\quad\mbox{
and}\quad
(\ref{OU-ref})\Longrightarrow
\epsilon_{\kappa(t)}(\rho)^{1/2}\leq c_0\rho~e^{-at}.
$$
\end{examp}

\subsubsection{Some quantitative estimates}\label{mr-sec}

For instance, combining (\ref{ref-convex-UV}) with (\ref{HW-UV}) we also check that
 \begin{equation}\label{ref-convex-UV-app}
\begin{array}{l}
\exists  \sigma,\overline{\sigma}\in  \Sa^+_d\quad\mbox{\rm such that}\quad \forall x\in\RR^d\quad
\nabla^2 U(x)\geq \sigma^{-1}\quad\mbox{\rm and}\quad \nabla^2 V(x)\geq \overline{\sigma}^{-1}\\
\\
\Longrightarrow\quad \left\{\begin{array}{l}
\mbox{\rm $\lambda_U$ and $\Ka_{2n+1}$ satisfy the $LS(\rho)$ inequality with $\rho=\Vert\sigma\Vert_2$}\\ 
\\
\mbox{\rm $\lambda_V$ and $\Ka_{2n}$ satisfy the $LS(\overline{\rho})$ inequality with $\overline{\rho}=\Vert\overline{\sigma}\Vert_2$}.
\end{array}\right.
\end{array}\end{equation}
More generally, assume $(U,V)$ satisfying the convex-at-infinity conditions
 $$
 \nabla^2 U(x)\geq 
a~1_{\Vert x\Vert_2 \geq  l}~I-b 
~1_{\Vert x\Vert_2< l}~I\quad\mbox{\rm and}\quad
\nabla^2 V(x)\geq 
\overline{a}~1_{\Vert x\Vert_2 \geq  \overline{l}}~I-\overline{b} 
~1_{\Vert x\Vert_2< \overline{l}}~I
 $$
 for some parameters $a,\overline{a}>0$, $b,\overline{b}\geq 0$ and $l,\overline{l}\geq 0$.
In this context, combining (\ref{ex-lg-cg-v2}) with (\ref{HW-UV}) there exists some parameters $\rho(a,b,l)\wedge \rho(\overline{a},\overline{b},\overline{l})>0$ such that
 \begin{equation}\label{ex-lg-cg-v2-app}
 \left\{\begin{array}{l}
\mbox{\rm $\lambda_U$ and $\Ka_{2n+1}$ satisfy the $LS(\rho)$ inequality with $\rho=\rho(a,b,l)$}\\ 
\\
\mbox{\rm $\lambda_V$ and $\Ka_{2n}$ satisfy the $LS(\overline{\rho})$ inequality with $\overline{\rho}:=\rho(\overline{a},\overline{b},\overline{l})$.}
\end{array}\right.
\end{equation}
Arguing as above we check that (\ref{pi2nP-v2}) and (\ref{pi2n1P-v2}) and thus (\ref{theo-edecay}) hold with the parameter
\begin{equation}\label{pi2nP-v2-loc}
 \epsilon_{\kappa}(\rho):=\kappa^{2}\rho(a,b,l)\rho(\overline{a},\overline{b},\overline{l}).
\end{equation}
The estimates (\ref{pi2nP-D2}) and resp. (\ref{pi2n1P-D2}) also holds with 
the parameter $\rho=\rho(\overline{a},\overline{b},\overline{l})$ and resp. $\rho=\rho(a,b,l)$, thus (\ref{theo-d2decay}) also holds with $ \epsilon_{\kappa}(\rho)$ as in (\ref{pi2nP-v2-loc}).

One can also show that the entropic exponential decays presented in (\ref{theo-edecay}) can be refined. 
We set 
$$
I_{n}:= \Ha(P_{\mu,\eta}~|~\Pa_{n})\quad\mbox{\rm and}\quad \varepsilon:=
 \epsilon_{\kappa}(\rho).
$$
Combining (\ref{pi2nP}) with Proposition~\ref{prop-decrease2n01} and (\ref{pi2n1P-v2}) we check that
\begin{eqnarray*}
I_{2n}&=&   I_{2n-1}-\Ha(\mu~|~\pi_{2n-1})\leq   I_{2n-1}- \Ha(\pi_{2n}~|~\eta)\leq I_{2n-1}- \varepsilon^{-1}~I_{2n+1}
\end{eqnarray*}
Now, using (\ref{pi2n1P}) along with Proposition~\ref{prop-decrease2n01} and (\ref{pi2nP-v2}) we have
\begin{eqnarray*}
 I_{2n+1}&=& 
  I_{2n}-   \Ha(\eta~|~\pi_{2n})\leq     I_{2n}-     \Ha(\pi_{2n+1}~|~\mu)\leq   I_{2n}-\varepsilon^{-1}    I_{2(n+1)}.
\end{eqnarray*}
In summary we have proved that for any $n\geq 0$ the second order linear homogeneous recurrence relation
\begin{equation}\label{rec-ref}
 I_{n}\geq I_{n+1}+  \varepsilon^{-1} I_{n+2} \geq (1+\varepsilon^{-1})~I_{n+2}.
\end{equation}
These linear recurrence equations are rather well known, for instance
following the proof of \cite[Theorem 6.9]{dp-25}, for any $n\geq 0$ we have
$$
\Ha(P_{\mu,\eta}~|~\Pa_{2(n+1)})\leq \left(1+{1}/{\varepsilon}+{1}/{\overline{\varepsilon}}\right)^{-n}~\Ha(P_{\mu,\eta}~|~\Pa_{0})
$$
with the parameter
$$
 1+1/\overline{\varepsilon}:=\frac{1+\sqrt{1+4/\varepsilon}}{2}\Longrightarrow (1+1/\overline{\varepsilon})^2=1+{1}/{\varepsilon}+{1}/{\overline{\varepsilon}}>1+{1}/{\varepsilon}.
$$

\subsubsection{Strongly convex models}\label{mr-stconv-sec}

Assume there exist some  positive definite matrices $\sigma\geq \sigma_-$ and $\overline{\sigma}\geq \overline{\sigma}_-$ such that
$$
 \sigma^{-1}\leq \nabla^2 U\leq  \sigma^{-1}_-\quad\mbox{\rm and}\quad \overline{\sigma}^{-1}\leq \nabla^2 V \leq \overline{\sigma}^{-1}_-.
$$
This rather strong convexity condition~\cite{chewi2023entropic,conforti2024weak,dp-25,durmus} is often used to improve the curvature estimates in (\ref{HW-UV-0}) and (\ref{HW-UV}). This section underlines some links between these studies and the Riccati matrix difference equations discussed in Section~\ref{ricc-equations-sec} in the context of linear-Gaussian models.

It should be noted that for linear Gaussian transition potential of the form (\ref{def-W}) the decompositions (\ref{HW-UV-0}) yield for any $n\geq 0$ the estimates
$$
\begin{array}{rclcl}
 \sigma^{-1}+ \cchi^{\prime}~
  \mbox{\rm cov}_{\Ka_{2n}}(x)~  \cchi&\leq &
\nabla^2_2W_{2n+1}^{\flat}(y,x)&\leq &  \sigma^{-1}_-+ \cchi^{\prime}~
  \mbox{\rm cov}_{\Ka_{2n}}(x)~  \cchi\\
  \\
\overline{\sigma}^{-1}+  \cchi~  \mbox{\rm cov}_{\Ka_{2n+1}}(y)~\cchi^{\prime}  &\leq &\nabla^2_2W_{2(n+1)}(x,y)&\leq &\overline{\sigma}^{-1}_-+  \cchi~  \mbox{\rm cov}_{\Ka_{2n+1}}(y)~\cchi^{\prime}  .
\end{array}$$
Applying (\ref{bl-cr}) to the Markov transition $\Ka_1$ we have
\begin{eqnarray*}
  \mbox{\rm cov}_{\Ka_{0}}(x)=\tau
 & \Longrightarrow& \tau_1^{-1}:=
  \sigma^{-1}+  \cchi^{\prime}~
 \tau~  \cchi \leq 
 \nabla^2_2W_{1}^{\flat}\leq \tau_{1-}^{-1}:= \sigma^{-1}_-+ \cchi^{\prime}~
 \tau~  \cchi 
\\
&\Longrightarrow&
 \tau_{1-}
\leq \mbox{\rm cov}_{\Ka_{1}}\leq  \tau_{1}
\end{eqnarray*}
with  $\cchi$ as in (\ref{ref-hessian-W-f}).
Applying (\ref{bl-cr}) to the Markov transition $\Ka_2$, this implies that
$$
\begin{array}{l}
\tau^{-1}_2:=
\overline{\sigma}^{-1}+  \cchi~  \tau_{1-}~\cchi^{\prime}\leq 
\nabla^2_2W_{2}\leq \tau^{-1}_{2-}:=\overline{\sigma}^{-1}_-+  \cchi~  \tau_{1}~\cchi^{\prime}
\\
\\
\Longrightarrow \tau_{2-}
\leq \mbox{\rm cov}_{\Ka_{2}}\leq\tau_2.
\end{array}$$
Iterating the procedure, we have the covariance estimates
$$
 \tau_{n-}
\leq \mbox{\rm cov}_{\Ka_{n}}\leq  \tau_{n}
$$
as well as the Riccati difference matrix inequalities
$$
\begin{array}{rclcl}
  \tau_{2n+1}^{-1}:= \sigma^{-1}+ \cchi^{\prime}~
 \tau_{(2n)-}~  \cchi&\leq &
\nabla^2_2W_{2n+1}^{\flat}&\leq &  \tau_{(2n+1)-}^{-1}:= \sigma^{-1}_-+ \cchi^{\prime}~
  \tau_{2n}~  \cchi\\
  \\
\tau_{2(n+1)}^{-1}:=\overline{\sigma}^{-1}+  \cchi~ \tau_{(2n+1)-}~\cchi^{\prime}  &\leq &\nabla^2_2W_{2(n+1)}&\leq &\tau_{2(n+1)-}^{-1}:=\overline{\sigma}^{-1}_-+  \cchi~ \tau_{2n+1}~\cchi^{\prime}  .
\end{array}
$$
For instance, if we set
$$
\overline{\tau}_{2n}:=\overline{\sigma}^{-1/2}~\tau_{2n}~ \overline{\sigma}^{-1/2}
\quad\mbox{\rm and}\quad
\overline{ \tau}_{(2n+1)-}:=
\sigma^{-1/2}_- \tau_{(2n+1)-} \sigma^{-1/2}_-
$$
we obtain
\begin{eqnarray*}
\overline{\tau}_{2(n+1)}&=&
\left(I+\gamma_-~
\overline{ \tau}_{(2n+1)-}~ \gamma_-^{\prime}\right)^{-1} \leq I\\
\overline{\tau}_{(2n+1)-}&=&
\left( I+  \gamma_-^{\prime}~
 \overline{\tau}_{2n}~ \gamma_-\right)^{-1}\quad\mbox{\rm with}\quad
\gamma_-:= \overline{\sigma}^{1/2}~ \cchi~ \sigma^{1/2}_-.
\end{eqnarray*}
Arguing as in (\ref{proof-ricc}) we conclude that
$$
\overline{\tau}_{2(n+1)}=
\mbox{\rm Ricc}_{\varpi_-}\left( \overline{\tau}_{2n}\right)\quad\mbox{\rm with}\quad
\varpi_-^{-1}:=\gamma_-\gamma_-^{\prime}.
$$
The monotone properties (\ref{ricc-maps-incr}) of Riccati maps provide several uniform upper bounds of the flow of covariance matrices $\overline{\tau}_{2n}$. These upper bounds also yield curvature estimates for the Hessian matrices $\nabla^2_2W_{2n}$. 
For a more detailed discussion one can see \cite[Appendix D]{dp-25}.

\subsection{Lyapunov techniques}

\subsubsection{Regularity conditions}

Consider the marginal measures $(\mu,\eta)=(\lambda_U,\nu_V)$ and the integral operators  $K$ and $K^{\flat}$ defined in
(\ref{ref-intro-UV}) and (\ref{def-Qa}) and (\ref{def-Ra}). We assume $U,V,W$ are locally bounded and 
$\Vert e^{-W}\Vert<\infty$. In addition $U,V$ have compact sub-level sets.
These conditions ensures that $(U,V,W)$ are lower bounded by real numbers $(U_{\star},V_{\star},W_{\star})$ so that
$(e^{-U},e^{-V})$ are uniformly bounded. 
We again recall the notation $\Ba_{\infty}(\XX)$ and $\Ba_{0}(\XX)$ from section \ref{sec:notation}.
We remark that for any $\delta>0$ we have
$
(e^{\delta U},e^{\delta V})\in (\Ba_{\infty}(\XX)\times \Ba_{\infty}(\YY))
$ and  set
\begin{equation}\label{UV-delta}
\Ua_{\delta}:=\delta U-W^{V}\quad\mbox{\rm and}\quad
\Va_{\delta}:=\delta V-W_{U}
\end{equation}
with the functions
$$
W^V(x):=\int~\nu_V(dy)~W(x,y)
 \quad \mbox{\rm and}\quad
W_U(y):=\int~\lambda_U(dx)~W(x,y).
$$ 
Consider the following condition:\\
{\it $(\mathtt{A}):$ There exist $\delta\in ]0,1/2[$ and  functions $(\theta_{u},\theta_{v})\in ( \Ba_0(\XX)\times \Ba_0(\YY))$ such that
\begin{equation}\label{UV-delta-theta}
\lambda_{(1-2\delta)U}(1)\vee \nu_{(1-2\delta)V}(1)<\infty\quad \mbox{as well as}\quad e^{-\Ua_{\delta}}\leq \theta_{u}
\quad \mbox{and}\quad
 e^{-\Va_{\delta}}\leq \theta_{v}.
\end{equation}}

Condition $(\mathtt{A})$ is slightly weaker than condition $\Ha^{\prime}_{\delta}(U,V)$ discussed in~\cite{adm-25}. Note that  $(\mathtt{A})$ is met with $(\theta_{u},\theta_{v})=(e^{-\Ua_{\delta}},e^{-\Va_{\delta}})$ as soon as the functions $ (\Ua_{\delta},\Va_{\delta})$  are locally bounded with compact sub-level sets.
Intuitively,  condition $(\mathtt{A})$ ensure that marginal-potentials  $(U,V)$ dominate at infinity the (integrated) transition-potentials $(W_U,W^V)$. In terms of potential energy conversion, the marginal potentials should be large enough to compensate high energy transition losses.

\begin{rmk}
Assume that $(\Ua_{\delta},\Va_{\delta})$ are lower bounded by functions $(\Ua^-_{\delta},\Va_{\delta}^-)$. Note that in this case for any $l>0$ we have
$$
\{\Ua_{\delta}\leq l\}\subset \{\Ua^-_{\delta}\leq l\}\quad\mbox{and}\quad
\{\Va_{\delta}\leq l\}\subset \{\Va^-_{\delta}\leq l\}.
$$
In this scenario,  we can choose $(\theta_{u},\theta_{v}):=(e^{-\Ua^-_{\delta}},e^{-\Va^-_{\delta}})\in (\Ba_0(\XX)\times \Ba_0(\YY))$ as soon as the functions $(\Ua^-_{\delta},\Va_{\delta}^-)$ are locally bounded with compact sub-level sets.
\end{rmk}

Assume that $W(x,y)\leq a_W(x)+b_W(y)$ for some locally bounded functions $(a_W,b_W)$ such that
$\lambda_U(a_W)\vee\nu_V(b_W)<\infty$.  Then, we have
$$
\Ua^-_{\delta}(x):=\delta U(x)-a_W(x)-\nu_V(b_V)\leq \Ua_{\delta}(x)\leq \delta U(x)-W_{\star}.
$$
This shows that $(\Ua_{\delta},\Ua^-_{\delta})$ and similarly $(\Va_{\delta},\Va^-_{\delta})$ 
with $\Va^-_{\delta}:=\delta V-b_W$ are locally bounded.  Finally, $(\Ua^-_{\delta},\Va^-_{\delta})$ have compact 
compact sub-level sets if and only if $(\delta U-a_W,\delta V-b_W)$ have compact 
compact sub-level sets.

\begin{rmk}\label{rmk-W-bounded}
When $\Vert W\Vert<\infty$ we have $w:=\Vert W^V\Vert\vee \Vert W_U\Vert<\infty$. In this case, we have
$$
e^{-\Ua_{\delta}}\leq \theta_{u}=e^{w}~e^{-\delta U}\in \Ba_{0}(\XX)\quad
\mbox{and}\quad 
e^{-\Va_{\delta}}\leq \theta_{v}=e^{w}~e^{-\delta V}\in \Ba_{0}(\YY)
$$
\end{rmk}

\noindent The following examples are taken from~\cite{adm-25,dmg-25}.
\begin{examp}\label{exam-beta}
Consider a Beta distribution on the open unit interval  $\XX:=]0,1[$ with shape parameters $a_u,b_u>0$ defined by
$$
\lambda_U(dx)=c_u~x^{a_u}(1-x)^{b_u}~\lambda(dx)
\quad
\mbox{with}\quad 
\lambda(dx)=1_{]0,1[}(x)~dx
$$
and some normalizing constant $c_u>0$. In this situation, for any $\delta>0$ we have
$$
\displaystyle\Longrightarrow U(x)=a_u~\log{\frac{1}{x}}+b_u~\log{\frac{1}{1-x}}+\log{\frac{1}{c_u}}
\quad
\mbox{and}\quad
e^{-\delta U(x)}=c_u^{\delta}~x^{\delta a_u}~(1-x)^{\delta b_u}
$$
By Remark~\ref{rmk-W-bounded}, when $W$ is bounded we have
$w:=\Vert W^V\Vert<\infty$ as well as the estimate
$$
e^{-\Ua_{\delta}}\leq \theta_{u}=c_{u,w}~x^{\delta a_u}~(1-x)^{\delta b_u}\in \Ba_{0}(\XX)
\quad
\mbox{with}\quad
c_{u,w}:=e^{w} c_u^{\delta}
$$
Thus, condition 
$(\mathtt{A})$ is met for any $\delta\in ]0,1/2[$
as soon as $W$ is bounded  and $(\lambda_U,\nu_V)$ are  both Beta distributions on  $\XX=]0,1[=\YY$ with shape parameters $a_u,b_u>0$ and $a_v,b_v>0$.
\end{examp}

\begin{examp}\label{exam-poly}
Assume that there exist  some parameters $a_1,b_1\in\RR$,  $a_2,b_2> 0$ and $q\geq p>0$ such that for any $x,y\in\XX=\YY=\RR^d$ we have
$$
W(x,y)\leq a_1+a_2~(\Vert x\Vert^{p}+\Vert y\Vert^{p})\quad \mbox{and}\quad
  U(x)\wedge V(x)\geq b_1+b_2~\Vert x\Vert^{q}.
  $$
  In this case, there exists some constant $c\in\RR$ such that
  $$
\Ua_{\delta}(x)\wedge \Va_{\delta}(x)\geq c+\delta b_2~\Vert x\Vert^{q}-
a_2~\Vert x\Vert^{p}.
  $$
 We conclude that (\ref{UV-delta-theta}) is satisfied
 for any $\delta\in ]0,1/2[$  as soon as $q>p$, and for any $\delta\in ]a_2/b_2,1/2[$ as soon as $b_2>2a_2$ and $q=p$. Indeed, in both situations we have
 $$
e^{-\Ua_{\delta}(x)}\vee e^{-\Va_{\delta}(x)}\leq \theta(x):=e^{-c}~\exp{\left(- \delta b_2~\Vert x\Vert^{q}
\left(1-
\frac{a_2}{\delta b_2}~\frac{1}{\Vert x\Vert^{q-p}}\right)\right)}
\longrightarrow_{\Vert x\Vert\rightarrow\infty}~0
 $$
  \end{examp}

A wealth of examples satisfying these regularity conditions  is presented in Section 5 in~\cite{adm-25} ranging from polynomial growth potentials and heavy tailed marginals on general normed spaces to more sophisticated boundary state space models, including semi-circle transitions, Beta, Weibull, exponential marginals as well as semi-compact models. This framework also applies to statistical finite mixture of the above models, including kernel-type density estimators of complex data distributions arising in generative modeling.

\subsubsection{A technical lemma}
Consider the integral operators
$$
K_V(x,dy):=e^{-W(x,y)}~\nu_V(dy)
\quad\mbox{\rm and}\quad
K^{\flat}_U(y,dx):=e^{-W^{\flat}(y,x)}~\lambda_U(dx)
$$
and set
$$
W_{V,U}:=\log{K_V(\exp{(W_U)}})
 \quad \mbox{\rm and}\quad
W^{U,V}:=\log{K^{\flat}_U(\exp{(W^{V}})}).
$$

\begin{lem} Assume $(\mathtt{A})$.  Then we have that
\begin{equation}\label{eWUVb}
\Vert \exp{W_{V,U}}\Vert\vee \Vert \exp{W^{U,V}}\Vert<\infty\quad \mbox{and}\quad
\Vert \exp{W_{(1-\delta)V,U}}\Vert\vee\Vert \exp{W^{(1-\delta)U,V}}\Vert<\infty.
\end{equation} In addition, for any $n\geq 0$ we have
\begin{eqnarray}
-W^{V}+\eta(V)&\leq &\log{
K(\exp{(-V_{2n})})}\nonumber\\
 -W_{U}&\leq &\log{K^{\flat}( \exp{ (-U_{2n} ) } )
}\leq -\eta(V)+W^{U,V}.\label{normct}
\end{eqnarray}
For any $n\geq 1$ we also have
\begin{equation}\label{IUVn}
\begin{array}{rcccl}
U-W^{V}+\eta(V)&\leq& U_{2n} &\leq &U+W_{V,U}\\
V-W_U&\leq& 
V_{2n}&\leq&-\eta(V)+
V+W^{U,V}\quad \mbox{and}\quad \log{K(\exp{(-V_{2n})})}\leq W_{V,U}.
\end{array}
\end{equation}
\end{lem}
\proof
Recalling that $U$ is locally bounded with compact sub-level sets, we have
$$
\lambda_{(1-\delta)U}(1)=\lambda_{(1-2\delta)U}\left(e^{-\delta U}\right)\leq \Vert e^{-\delta U}\Vert~\lambda_{(1-2\delta)U}(1)<\infty
\quad\mbox{\rm and}\quad
 \nu_{(1-\delta)V}(1)<\infty.
$$
When $c:=\Vert e^{-W}\Vert<\infty$ for any $\delta\in ]0,1[$ we have
\begin{eqnarray*}
\exp{(W_{V,U})}\leq c~\nu_V(e^{W_U})=c~\nu_{(1-\delta)V}(e^{-\Va_{\delta}})<\infty
&\Longrightarrow &\Vert e^{W_{V,U}}\Vert<\infty\\
\exp{(W^{U,V})}\leq c~\lambda_U(e^{W^V})=c~\lambda_{(1-\delta)U}(e^{-\Ua_{\delta}})<\infty
&\Longrightarrow& \Vert e^{W^{U,V}}\Vert<\infty.
\end{eqnarray*}
Similiarly, one has
\begin{eqnarray*}
\exp{(W_{(1-\delta)V,U})}=K_{(1-\delta) V}(e^{W_U})&\leq& c~\nu_{(1-\delta)V}(e^{W_U})=
c~\nu_{(1-2\delta)V}(e^{-\Va_{\delta}})\\
\exp{(W^{(1-\delta)U,V})}=K^{\flat}_{(1-\delta) U}(e^{W^V})&\leq& c~ \lambda_{(1-\delta)U}(e^{W^V})=c~
\lambda_{(1-2\delta)U}(e^{-\Ua_{\delta}}).
\end{eqnarray*}

Applying Jensen's inequality and Proposition~\ref{UVn-prop-series}, for any $n\geq 0$ we have
\begin{eqnarray*}
\log{K(\exp{(-V_{2n})})(x)}&=&\log{\int~\eta(dy)~\exp{(-W(x,y)+V(y)-V_{2n}(y))}}\\
&\geq&
-W^V(x)+\eta(V)-\eta(V_{2n})\geq -W^{V}(x)+\eta(V).
\end{eqnarray*}
Note that 
$$
U_{2n}=U+\log{K(e^{-V_{2n}})}\geq U-W^{V}+\eta(V).
$$
In addition,  for any $n\geq 0$ we check that
 \begin{eqnarray*}
\log{K^{\flat}(\exp{(-U_{2n})})(y)}&=&\log{\int~\mu(dx)~\exp{(-W(x,y)+U(x)-U_{2n}(x))}}\\&\geq&
-W_{U}(y)+\mu(U)-\mu(U_{2n})\geq -W_{U}(y)
\end{eqnarray*}
as well as
$$
\log{K^{\flat}(\exp{(-U_{2n})})}\leq-\eta(V)+ 
\log{K^{\flat}(\exp{(-(U-W^{V})})}=-\eta(V)+W^{U,V}.
$$
By (\ref{UVn}) we also have
$$
V-W_U\leq 
V_{2(n+1)}=V+\log{K^{\flat}(e^{-U_{2n}})}\leq
V-\eta(V)+ 
W^{U,V}.
$$
In a similar manner, we have
$$
U-W^{V}+\eta(V)\leq
U_{2(n+1)}=U+\log{K(e^{-V_{2(n+1)}})}\leq U+W_{V,U}.
$$
For any $n\geq 1$ 
$$
\exp{(-W^{V}+\eta(V))}\leq K(\exp{(-V_{2n})})\leq 
K(\exp{(-(V-W_U))})=K_V(\exp{(W_U)}).
$$
This ends the proof of (\ref{IUVn}).
\cqfd

\subsubsection{Some uniform estimates}

In this section we assume $(\mathtt{A})$ is satisfied and we set
\begin{equation}\label{def-gh-pi-2}
(g,h):=(\exp{(\delta U)},\exp{(\delta V)})
\end{equation}

\begin{prop}
There exists some constant $c>0$ such that, for any $n\geq 1$, we have
\begin{equation}\label{W-uvb-2}
\frac{d\pi_{2n}}{d\eta}\vee \frac{d\eta}{d\pi_{2n}}\leq c~\exp{(W_U)}\quad \mbox{and}\quad
\frac{d\pi_{2n+1}}{d\mu}\vee\frac{d\mu}{d\pi_{2n+1}}\leq c~\exp{(W^V)}.
\end{equation}
We also have the uniform estimates
\begin{equation}
\Vert{d\pi_{2n}}/{d\eta}\Vert_{h}\vee \Vert{d\eta}/{d\pi_{2n}}\Vert_h\leq c\quad \mbox{and}\quad
\Vert{d\pi_{2n+1}}/{d\mu}\Vert_g\vee \Vert{d\mu}/{d\pi_{2n+1}}\Vert_g\leq c\label{def-gh-pi}
\end{equation}
\end{prop}
\proof
Using (\ref{ratio-UV}) and (\ref{IUVn}) for any $n\geq 0$ we have
$$
-(W^{V}+W_{V,U}-\eta(V))\leq U_{2(n+1)}-U_{2n}=\log{\left(\frac{d\pi_{2n+1}}{d\mu}\right)}\leq 
W^{V}+W_{V,U}-\eta(V)
$$
as well as
$$
-(
W_U+W^{U,V}-\eta(V))
\leq V_{2(n+1)}-V_{2n}=\log{
\left(
\frac{d\pi_{2n}}{d\eta}
\right)}\leq
W_U+W^{U,V}-\eta(V).
$$
We end the proof of (\ref{W-uvb-2}) taking the exponential  and using (\ref{eWUVb}).
The last assertion follows from the fact that
$$
e^{-\delta V}\left(\frac{d\pi_{2n}}{d\eta}\vee \frac{d\eta}{d\pi_{2n}}\right)\leq c~\exp{(-\Va_{\delta})}\quad
\mbox{\rm and}\quad
e^{-\delta U}\left(\frac{d\pi_{2n+1}}{d\mu}\vee\frac{d\mu}{d\pi_{2n+1}}\right)\leq c~\exp{(-\Ua_{\delta})}.
$$
This ends the proof of the proposition.
\cqfd

\begin{lem}\label{key-lemma-Lyap}
There exists $a>0$ such that for any $n\geq 1$ we have
$$
\Ka_{2n}\left(h\right)/g\leq 
a~e^{-\Ua_{\delta}}~\quad\mbox{and}\quad
\Ka_{2n+1}\left(g\right)/h\leq 
a~e^{-\Va_{\delta}}.~
$$
In addition, there exists some $\imath>0$ such that for any $n\geq 1$ we have
$$
\imath~K_V(x,dy)
\leq
\Ka_{2n}(x,dy)\quad\mbox{and}\quad
\imath~K^{\flat}_U(y,dx)\leq 
\Ka_{2n+1}(y,dx).
$$
\end{lem}

\proof
Combining (\ref{Kan}) with  (\ref{eWUVb}), (\ref{normct}) and (\ref{IUVn}) there exists some
$c_0>0$ such that for any $n\geq 1$ we have
$$
\frac{1}{c_0}~K_V(x,dy)
\leq
\Ka_{2n}(x,dy)=\frac{K_V(x,dy) ~e^{(V-V_{2n})(y)}}{K(e^{-V_{2n}})(x)}\leq c_0~e^{W^V(x)}~K_V(x,dy)~ e^{W_U(y)}.
$$
This implies that
$$
e^{-\delta U}\Ka_{2n}\left(e^{\delta V}\right)\leq 
c_0~e^{-\Ua_{\delta}}~K_{(1-\delta)V}\left(e^{W_U}\right)\leq 
c_1~e^{-\Ua_{\delta}}~
\quad \mbox{\rm with}\quad
c_1:=c_0~\Vert\exp{W_{(1-\delta)V,U}}\Vert.
$$
Then, recalling that $U_{2n+1}=U_{2n}$ we have that
there exists some
$c_0>0$ such that for any $n\geq 1$ we have
$$
\frac{1}{c_0}~K^{\flat}_U(y,dx)\leq 
\Ka_{2n+1}(y,dx)=\frac{K^{\flat}_U(y,dx)~e^{(U-U_{2n})(x)}}{K^{\flat}(e^{-U_{2n}})(y)}\leq c_0~e^{W_U(y)}~K^{\flat}_U(y,dx)~e^{W^V(x)}.
$$
This implies that
$$
e^{-\delta V}\Ka_{2n+1}\left(e^{\delta U}\right)\leq 
c_0~e^{-\Va_{\delta}}~K^{\flat}_{(1-\delta)U}(e^{W^V})\leq 
c_1~e^{-\Va_{\delta}}~
\quad \mbox{\rm with}\quad
c_1:=c_0~\Vert \exp{W^{(1-\delta)U,V}}\Vert.
$$
This ends the proof of the lemma.\cqfd

\subsubsection{A contraction theorem}

\begin{theo}\label{th-lips}
Assume $(\mathtt{A})$.  Then
there exists some $a>0$ and $\varrho\in ]0,1[$
such that for any $n\geq 1$ we have
\begin{equation}\label{lip-sinkhorn}
{\sf lip}_{g_a,h_a}(\Ka_{2n})\vee {\sf lip}_{h_a,g_a}(\Ka_{2n+1})\leq \varrho
\quad\mbox{with}\quad
 g_a:=\frac{1}{2}+a~g
\quad \mbox{and}\quad
h_a:=\frac{1}{2}+a~h
\end{equation}
and $(g,h)$ as in (\ref{def-gh-pi-2}). In addition, we have the uniform estimates
\begin{equation}\label{unif-pi-K}
\sup_{n\geq 1}(\pi_{2n}(h)\vee
 \pi_{2n+1}(g))<\infty\quad\mbox{and}\quad
\sup_{n\geq 1} \left(\Vert \Ka_{2n}\left(h\right)\Vert_g\vee
  \Vert \Ka_{2n+1}\left(g\right)\Vert_h\right)<\infty.
\end{equation}
\end{theo}
\proof

For any $\epsilon\in ]0,1[$, choosing $r>0$ such that
$$
a~e^{-r}\leq \epsilon\quad \mbox{\rm and set}\quad c:=
a~\left(\sup_{\Ua_{\delta}(x)\leq r}\left(e^{-\Ua_{\delta}(x)}~g(x)\right)\vee
\sup_{\Va_{\delta}(y)\leq r}\left(e^{-\Va_{\delta}(y)}~h(y)\right)
\right).
$$
In this notation, we have
$$
\begin{array}{l}
\displaystyle\Ka_{2n}\left(h\right)\leq 
\epsilon~g+c\quad\mbox{and}\quad
\Ka_{2n+1}\left(g\right)\leq 
\epsilon~h+c.
\end{array}
$$
This implies that
\begin{eqnarray}
\pi_{2n}(h)=\lambda_U\Ka_{2n}(h)&\leq& \epsilon~\lambda_{(1-\delta)U}(1)+c\nonumber
\\
 \pi_{2n+1}(g)=\nu_V\Ka_{2n+1}(g)&\leq &\epsilon~\nu_{(1-\delta)V}(1)+c.\label{unif-pin}
\end{eqnarray}
We end the proof of (\ref{unif-pi-K}) recalling that $(U,V)$ are bounded below by a real number.

Since $\XX$ and $\YY$ can be exhausted respectively by the compact sub-levels sets
$\{g\leq l\}$ and $\{h\leq l\}$, there exists some $l_1>1$ such that 
$ \lambda\left(\{g\leq l\}\right)\wedge\nu\left(\{h\leq l\}\right)>0$ for $l=l_1$ and thus for any $l\geq l_1$. We set $\jmath_l:=\sup_{g(x)\vee h(y)\leq l}W(x,y)$.
Then, for any $l\geq l_1$ and $g(x)\leq l$ we have
$$
K_V(x,dy)\geq \iota_V(l)~\frac{\nu_V(dy)~1_{\{h\leq l\}}(y)}{\nu_V(\{h\leq l\})}
\quad \mbox{\rm with}\quad
\iota_V(l):=e^{-\jmath_l}~\nu_V(\{h\leq l\})
$$
Note that
$$
\nu_V(\{h\leq l\})=\nu(e^{-V}~1_{V\leq \log{(l)}/\delta})\geq 
e^{- \log{(l)}/\delta}~\nu\left(\{h\leq l\}\right)>0
$$
In the same vein, for any $l\geq l_1$ and $h(y)\leq l$ we have
$$
K^{\flat}_U(y,dx)
\geq \iota_U(l)~\frac{\lambda_U(dx)~1_{\{g\leq l\}}(x)}{\lambda_U(\{g\leq l\})}
\quad \mbox{\rm with}\quad
\iota_U(l):=e^{-\jmath_l}~\lambda_U(\{g\leq l\}).
$$
This implies that for any $g(x_1)\vee g(x_2)\leq l$ and $h(y_1)\vee h(y_2)\leq l$ we have
$$
\Vert \delta_{x_1} \Ka_{2n}-\delta_{x_2} \Ka_{2n}\Vert_{\tiny tv}\vee \Vert \delta_{y_1} \Ka_{2n+1}-\delta_{y_2} \Ka_{2n+1}\Vert_{\tiny tv}\leq 1-\iota(l)
$$
with $\iota(l):=\imath~(\iota_U(l)\wedge \iota_V(l))$.
By Theorem~\ref{theo-lipa}, the contraction inequalities (\ref{lip-sinkhorn}) are satisfied 
for some $a>0$ and $\varrho\in ]0,1[$ that only depends on the parameters $(\epsilon,c)$ and the function $\iota(l)$ and hence we conclude. \cqfd 

\noindent The next corollary is a consequence of (\ref{lip-sinkhorn}).
\begin{cor}\label{cor-lip-pi}
Under the assumptions of Theorem~\ref{th-lips},
we have the uniform contraction estimates
$$
\sup_{n\geq 1}\left({\sf lip}_{g_a,g_a}(\Sa_{2n+1})\vee
 {\sf lip}_{h_a,h_a}(\Sa_{2(n+1)})\right)
\leq \varrho^2.
$$
with $\varrho\in ]0,1[$ as in (\ref{lip-sinkhorn}).
In addition, for any $n\geq 1$ we have
$$
\vertiii{\pi_{2(n+1)}-\eta}_{h_a}\leq \varrho^2~\vertiii{\pi_{2n}-\eta}_{h_a}
\quad\mbox{\rm and}\quad
\vertiii{\pi_{2n+1}-\mu}_{g_a}\leq \varrho^2~\vertiii{\pi_{2n-1}-\mu}_{g_a}.
$$
\end{cor}
Combining Corollary~\ref{cor-lip-pi} with the uniform estimates (\ref{unif-pin}) there exists some constant $c>0$ such that for any $n\geq 1$ we have
\begin{equation}\label{inter-wn2Dpsi}
\vertiii{\pi_{2n+1}-\mu}_{g_a}\vee\vertiii{\pi_{2(n+1)}-\eta}_{h_a}\leq c~\varrho^{2n}.
\end{equation}
Following the proof of (\ref{wn2Dpsi}) using (\ref{inter-wn2Dpsi}) we obtain the following Kantorovich semi-distance exponential decays.
\begin{cor}
Assume the assumptions of Theorem~\ref{th-lips} are met 
and there exists semi-distance $ \psi_g$ and $\psi_h$ on $\XX$ and $\YY$ such that
$$
\psi_g(x_1,x_2)\leq w_g(x_1,x_2)\quad \mbox{\rm and}\quad
\psi_h(y_1,y_2)\leq w_h(y_1,y_2).
$$
In this situation, for any $n\geq 1$ we have
$$
D_{\psi_g}(\pi_{2n+1},\mu)\vee
D_{\psi_h}(\pi_{2(n+1)},\eta) \leq (c/a)~\varrho^{2n}
$$
with $\varrho\in ]0,1[$ as in (\ref{lip-sinkhorn}) and the Kantorovich semi-distances
$(D_{\psi_g},D_{\psi_h})$ defined in (\ref{ref-cost}).
\end{cor}

\begin{examp}\label{exam-betav2}
Consider the Beta marginal distribution $\lambda_U$ on $\XX:=]0,1[$ discussed in Example~\ref{exam-beta}. In this case, recalling that for any $p>0$ and $u,v\geq 1$ we have
$$
2uv\geq u+v
\quad \mbox{and}\quad
|u-v|^p\leq (u+v)^p\leq 2^{(p-1)_+}~(u^p+v^p)
$$ 
and one has
$$
\displaystyle g(x):=e^{\delta U}=\frac{1}{c^{\delta}}~\left(\frac{1}{x}\right)^{\delta a}~\left(\frac{1}{1-x}\right)^{\delta b}\geq 
\frac{1}{2c^{\delta}}~\left(\left(\frac{1}{x}\right)^{p_a}+\left(\frac{1}{1-x}\right)^{p_b}\right)
$$
with $p_a,p_b):=(\delta a,\delta b)$.  Therefore
$$
w_g(x_1,x_2)\geq \psi_g(x_1,x_2):=c_{a,b}
\left(\left\vert\frac{1}{x_1}-\frac{1}{x_2}\right\vert^{p_a}+\left\vert\frac{1}{1-x_1}-\frac{1}{1-x_2}\right\vert^{p_b}\right)\quad
\mbox{for some}\quad c_{a,b}>0.
$$

\end{examp}

\begin{examp}
Consider the marginal distribution $\lambda_U$ on $\XX:=\RR^d$ discussed in Example~\ref{exam-poly} for some parameters $(b_1,b_2,q)$. In this case we have
$$
\displaystyle g(x):=e^{\delta U}\geq c~\exp{\left(2b\Vert x\Vert^q\right)}
\quad
\mbox{with}\quad (c,b):=\left(e^{b_1},b_2/2\right).
$$
Following Example~\ref{g-lyap-ill} several 
semi-distances $ \psi_g$ can be chosen, for instance we have
\begin{eqnarray*}
w_g(x_1,x_2)&\geq&2c~\left(\exp{\left(b\left(\Vert x_1\Vert^q+\Vert x_2\Vert^q\right)\right)}-1\right)\\
&=& \psi_g(x_1,x_2):=
2c~\left(\exp{\left(\frac{b}{ 2^{(q-1)_+}}~\Vert x_1- x_2\Vert^q\right)}-1\right).
\end{eqnarray*}

\end{examp}

\subsubsection{Entropy estimates}

\begin{prop}\label{prop-series-radon-gh}
Under the assumptions of Theorem~\ref{th-lips}, there exists some $c>0$ such that for any $n\geq 1$ we have
$$
\left\Vert
\frac{d\pi_{2(n+1)}}{d\eta}-1\right\Vert_h\vee \left\Vert
\frac{d\eta}{d\pi_{2(n+1)}}-1\right\Vert_h\vee \left\Vert\log{\frac{d\eta}{d\pi_{2(n+1)}}}\right\Vert_h
\leq c~\varrho^{2(n+1)}
$$
with $\varrho\in ]0,1[$ as in (\ref{lip-sinkhorn}).
The above  estimates remain valid when we replace $(\pi_{2(n+1)},\eta,h)$ by 
$(\pi_{2n+1},\mu,g)$. 
\end{prop}
\proof
Arguing as in the proof of Proposition~\ref{prop-series-radon}
for any $k,l\geq 1$ we have
$$
\left\vert\frac{d\pi_{2 (l+k)}}{d\eta}(y)-1\right\vert
\leq \Vert d\pi_{2l}/{d\eta}\Vert_{h_a}~
 \vertiii{
\left(\delta_y-\eta\right)\Sa_{2(l+k)}\ldots \Sa_{2(l+1)}}_{h_a}.
$$
Note that
$$
\Vert f\Vert_{h_a}\leq (2\Vert f\Vert)\wedge \left(\Vert f\Vert_{h}/a\right).
$$
By (\ref{def-gh-pi}) there exists some constant $c>0$ such that for any $k\geq 1$ we have
$$
\left\vert\frac{d\pi_{2 k}}{d\eta}(y)-1\right\vert
\leq c~\varrho^{2k}~(1+h(y))
\Longrightarrow \
\left\Vert\frac{d\pi_{2k}}{d\eta}-1\right\Vert_{h}
\leq c~\varrho^{2k}.
$$
In addition,  for any $k,l\geq 1$ we have
$$
\begin{array}{l}
\displaystyle\left\vert\frac{d\pi_{2(l+k)+1}}{d\mu}(x)-1\right\vert\leq \Vert {d\pi_{2l+1}}/{d\mu}\Vert_{g_a}~\vertiii{
(\delta_x-\mu)\Sa_{2(l+k)+1}\ldots \Sa_{2(l+1)+1}}_{g_a}
\\
\\
\displaystyle\Longrightarrow\quad \forall k\geq 1\qquad \left\Vert\frac{d\pi_{2k+1}}{d\mu}-1\right\Vert_{g}
\leq c~\varrho^{2k}.
\end{array}$$

By Proposition~\ref{prop-commute} we have
$$
\begin{array}{l}
\displaystyle
\frac{d\eta}{d\pi_{2(k+1)}}-1=
\Ka_{2(k+1)+1}\left(\frac{d\pi_{2k+1}}{d\mu}-1\right)
\quad\mbox{\rm and}\quad
\frac{d\mu}{d\pi_{2k+1}}-1=
\Ka_{2(k+1)}\left(\frac{d\pi_{2k}}{d\eta}-1\right)
\\
\\
\displaystyle\Longrightarrow\quad
\left\Vert \frac{d\eta}{d\pi_{2(k+1)}}-1\right\Vert_h\leq \left\Vert\frac{d\pi_{2k+1}}{d\mu}-1\right\Vert_{g}~\Vert \Ka_{2(k+1)+1}(g)\Vert_h.
\end{array}$$
Using (\ref{log-ratio}) and (\ref{unif-pi-K}) we conclude that for any $k\geq 2$ we have
$$
\left\Vert\frac{d\pi_{2 k}}{d\eta}-1\right\Vert_{h}\vee \left\Vert \frac{d\eta}{d\pi_{2k}}-1\right\Vert_h\vee 
\left\Vert \log{\frac{d\eta}{d\pi_{2k}}}\right\Vert_h
\leq c~\varrho^{2k}.
$$
This ends the proof of the proposition.\cqfd

Proposition~\ref{prop-series-radon-gh} ensures the convergence of the 
series (\ref{UU-VV}) in the weighted Banach spaces $\Ba_g(\XX)$ and  $\Ba_h(\YY)$.
Combining Proposition~\ref{UVn-prop-series} with (\ref{UU-VV}) and using Proposition~\ref{kdec2n-prop} we obtain the following theorem.

\begin{theo}\label{ac-series-theo} 
Under the assumptions of Theorem~\ref{th-lips}, there exists some constant $c>0$ such that  for any $n\geq 1$ we have
\begin{equation}\label{ac-series} 
\|\UU-U_{2n}\|_{g}\vee \|\VV-V_{2n}\|_{h}\vee   \Ha(P_{\mu,\eta}~|~\Pa_{2n}) \leq c~\varrho^{2n}
\end{equation}
with $\varrho\in ]0,1[$ as in (\ref{lip-sinkhorn}) and
with $(U_{2n},V_{2n})$ and $(\UU,\VV)$ as in Proposition~\ref{UVn-prop-series} and (\ref{UU-VV}).
\end{theo}

\end{document}